\documentclass[preprint]{elsarticle}
\usepackage[hidelinks]{hyperref}
\usepackage{hypernat}
\usepackage{lineno} 
\usepackage[T1]{fontenc}
\usepackage{graphicx,wrapfig}
\usepackage{subfiles}
\usepackage{textcomp,gensymb}
\usepackage{color}
\usepackage{caption}
\usepackage{subcaption}
\usepackage{amsmath,amssymb,eucal}
\usepackage{mathtools} 
\usepackage{amsfonts}
\usepackage{upgreek}
\usepackage{mathrsfs}
\usepackage{stmaryrd}
\usepackage{float}
\usepackage{amsthm}
\usepackage{bm} 
\usepackage{colortbl}
\usepackage{algpseudocode,algorithm}
\usepackage{setspace}
\usepackage{lipsum}
\usepackage[normalem]{ulem}
\usepackage[section]{placeins} 
\usepackage[symbol]{footmisc} 
\floatstyle{plaintop}
\newfloat{mybox}{tbhp}{}
\floatname{mybox}{\bf{Box}}
\theoremstyle{remark}

\let\oldequation\equation
\let\oldendequation\endequation
\renewenvironment{equation}
  {\linenomathNonumbers\oldequation}
  {\oldendequation\endlinenomath}

\newcommand\rd{{\rm{d}}}

\def\Int{\int\limits}
\newcommand{\fem}[1]{\textbf{\textsf{#1}}}
\newcommand\mysquare{\mbox{\footnotesize$\square$}}

\DeclareMathOperator{\dive}{div}

\journal{Computer Methods in Applied Mechanics and Engineering}
\begin{document}
\begin{frontmatter}
\title{Towards a sharper phase-field method: a hybrid diffuse--semisharp approach for microstructure evolution problems$^\dagger$}
\author[IPPT]{J\k{e}drzej Dobrza\'{n}ski}
\ead{jdobrz@ippt.pan.pl}
\author[IPPT]{Stanis{\l}aw Stupkiewicz\corref{cor1}}
\ead{sstupkie@ippt.pan.pl}

\cortext[cor1]{Corresponding author}

\address[IPPT]{Institute of Fundamental Technological Research, Polish Academy of Sciences,\\
Pawi\'nskiego 5B, 02-106 Warsaw, Poland.}

\begin{abstract}
A new approach is developed for computational modelling of microstructure evolution problems. The approach combines the phase-field method with the recently-developed laminated element technique (LET) which is a simple and efficient method to model weak discontinuities using non-conforming finite-element meshes. The essence of LET is in treating the elements that are cut by an interface as simple laminates of the two phases, and this idea is here extended to propagating interfaces so that the volume fraction of the phases and the lamination orientation vary accordingly. In the proposed LET-PF approach, the phase-field variable (order parameter), which is governed by an evolution equation of the Ginzburg--Landau type, plays the role of a level-set function that implicitly defines the position of the (sharp) interface. The mechanical equilibrium subproblem is then solved using the semisharp LET technique. Performance of LET-PF is illustrated by numerical examples. 
In particular, it is shown that, for the problems studied, LET-PF exhibits higher accuracy than the conventional phase-field method so that, for instance, qualitatively correct results can be obtained using a significantly coarser mesh, and thus at a lower computational cost.
\footnotetext[2]{Published in \emph{Comput. Methods Appl. Mech. Engrg.}, \textbf{423}, 116841, 2024, doi: 10.1016/j.cma.2024.116841}
\end{abstract}

\end{frontmatter}

\section{Introduction}
    The phase-field method is a well-established, versatile and powerful approach to computational modelling of microstructure evolution problems \citep{Chen2002,Steinbach2009,ProvatasElder2010,WangLi2010,Tourret2022}. 
    Microstructure evolution is inherently associated with propagation of interfaces that separate individual phases, and this constitutes the main challenge in developing computational approaches for the corresponding problems. The essence of the phase-field method is in treating the interfaces as diffuse, rather than sharp, 
    which is achieved by introducing a phase-field variable (so-called order parameter) that differentiates the phases and varies in a continuous manner within the diffuse interfaces. Accordingly, propagation of interfaces can be simulated on a fixed finite-element mesh (this work focuses on the finite element method). Tracking of interfaces is thus avoided, which would otherwise require, for instance, adaptive remeshing or enhancement techniques such as X-FEM.
    
    The phase-field method is highly versatile because the diffuse-interface framework can be combined with virtually any physics and, in fact, has been used in numerous contexts, for instance, solidification \citep{Ode2001}, solid-state transformations \citep{Chen2002,Ubachs2004,GuinKochmann2023}, including martensitic transformation \citep{WangKhachaturyan1997,LevitasPreston2002part1,Bo2020,Tuma2021} and twinning \citep{ClaytonKnap2011,Liu2018,Rezaee2022}, fracture \citep{Bourdin2000,Ambati2014}, corrosion \citep{Cui2021}, and many more.
    
    Another important feature of the phase-field method is that the order parameter carries information about the energy of the interfaces. The interfacial energy is naturally included in the energy balance and thus contributes to the driving force for interface propagation, 
    which is particularly important at smaller scales at which the interfacial energy effects are more pronounced. 
    It is also the interfacial energy term in the free energy function that sets the thickness of the diffuse interfaces through the term involving the gradient of the order parameter. 
    Moreover, considering that the contribution of interfaces to the total energy of a system is size-dependent, the phase-field method is able to model the size effects \citep{Tuma2016,Yeddu2018,Rezaee2020}.
    
    However, the beneficial features of the phase-field method, as discussed above, come at a price of a high computational cost. This is because a sufficiently fine computational grid (mesh) is needed to correctly represent the profile of the order parameter across the diffuse interface and to correctly represent the corresponding interfacial energy. This sets a severe constraint on the maximum allowable element size, and the associated high computational cost limits the maximum physical dimension of the domain that can be considered in the computations.
    
    
    With the aim to overcome the limitations discussed above, the so-called `sharp phase-field method' has been developed by \citet{FinelPRL2018}, see also \citep{Dimokrati2020,Fleck2022,Fleck2023}. The method allows the computational grid to be larger than the theoretical interface thickness, hence significantly coarser meshes can be used compared to the conventional phase-field method, and thus larger physical domains can be effectively simulated. However, the method relies on the notion of the discrete gradient, which is limited to the finite difference method and does not generalize to other discretization techniques, in particular, to the finite-element method.
    
    
    In this work, we address the mentioned limitations of the phase-field method from a different perspective. Focusing on microstructure evolution problems involving elastic interactions and on the finite element method for spatial discretization, 
    we propose a hybrid diffuse--semisharp approach. 
    The microstructure, including diffuse interfaces and their evolution, is described in a phase-field-like fashion using a continuous order parameter, and the order parameter plays the role of a level-set function that implicitly defines the position of the sharp interfaces. Considering that the position of the interfaces is arbitrary (and changes in a continuous manner), a fixed finite-element mesh does not conform to those sharp interfaces. To solve the corresponding mechanical equilibrium subproblem, we employ a recently-developed semisharp approach, namely the laminated element technique (LET) \citep{Dobrzanski2024}. In this approach, the elements that are cut by an interface are treated as laminates of the two involved phases (hence the name of the method), so that the interface-affected elements create a thin (one element thick) layer of elements between the two phases, and the method can thus be classified as semisharp. The method is inspired by the composite voxel technique that has been developed for FFT-based solvers in computational homogenization \citep{Gelebart2015,Kabel2015,Kabel2017,Mareau2017,Keshav2022}. In the present context, it is an important feature of LET that the response is a continuous function of the position and orientation of the interface \citep{Dobrzanski2024}, hence the method is suitable for problems involving moving interfaces.
    
    In the proposed LET-PF method, the two phases are mixed only within a thin layer of laminated elements (of the thickness of one element) along the sharp interface, while in the conventional phase-field method the phases are mixed within the entire diffuse interface that spans several elements (at least 3--4, say). Note that the mixing based on laminated microstructures, thus in a sense similar to LET-PF, is employed in some formulations of the phase-field method \citep{Durga2013,Mosler2014,Schneider2015,Bartels2017}, however, consistent with the phase-field framework, it is applied within the entire diffuse interface. 
    
    The proposed approach bears some similarity to the level-set method \citep{Sethian-book1999} combined with X-FEM \citep{Moes1999,Sukumar2001,Moes2003}. That approach appears attractive because X-FEM is a sharp-interface technique so that the interfaces can be modelled with a high accuracy using a non-matching mesh. 
    However, compared to the phase-field method, the related applications for microstructure evolution problems are much more scarce, see e.g.\ \citep{Ji2002,Duddu2011,Munk2022}. 
    The LET approach used in LET-PF is not as accurate as X-FEM \citep{Dobrzanski2024}, but it is much simpler, in particular, its implementation is carried out solely at the element level and does not require intervention into the structure of the finite-element code. 
    Another pronounced difference between the two approaches is in the way the interfacial energy is treated. The level-set function in the level-set method is a signed distance function and its role is only to (implicitly) define the position of the sharp interface. The interfacial energy (density per unit area) is then assigned directly to the sharp interface. In LET-PF, a diffuse-interface description of the interfacial energy is adopted as in the conventional phase-field method so that the interfacial energy contribution to the total energy is a bulk-like energy (density per unit volume) 
    distributed in the vicinity of the interfaces.
    
    
    In the following, the presentation is limited to the small-strain framework, however, extension to finite deformations is immediate. The implementation and the examples are performed in 2D, 
    {and}, again, extension to 3D is immediate. 
    On the other hand, the proposed LET-PF approach is limited to the case of two phases, but generalization to a multi-phase framework is not immediate, see Remark~2 in Section~\ref{sec:LET}.
    
    The remainder of the paper is organized as follows. In Section~\ref{sec:problem}, the reference sharp-interface microstructure evolution problem is formulated and its conventional phase-field counterpart is described, including the finite-element treatment of the latter. The new hybrid diffuse--semisharp approach (LET-PF) is introduced in Section~\ref{sec:hybrid}. Numerical examples illustrating the performance of LET-PF are reported in Section~\ref{sec:numExamples}. 
\section{Microstructure evolution problem}\label{sec:problem}
    This section briefly summarizes the classical formulation of the microstructure evolution problem in the sharp-interface framework (Section~\ref{sec:sharp}) and in the diffuse-interface framework of the phase-field method (Section~\ref{sec:diffuse}). Also, the finite-element treatment of the simple phase-field model introduced in Section~\ref{sec:diffuse} is discussed in Section~\ref{sec:PFM}. 
    \subsection{Reference sharp-interface problem}\label{sec:sharp}
        Let us consider domain $\Omega$ with two subdomains $\Omega_1$ and $\Omega_2$ separated by an interface $\Gamma$. 
        Each subdomain $\Omega_i$ ($i=1,2$) is occupied by an elastic phase characterized by the (isothermal) Helmholtz free energy density $\psi_i$ of the form
        \begin{equation}
        \label{eq:psiSharp}
        \psi_i (\bm{\varepsilon}) = \psi^0_i + \frac12 (\bm{\varepsilon}-\bm{\varepsilon}^{\rm t}_i) : \mathbb{L}_i : (\bm{\varepsilon}-\bm{\varepsilon}^{\rm t}_i) ,
        \end{equation}
        so that the local constitutive equation reads
        \begin{equation}
        \bm{\sigma} = \frac{\partial\psi_i}{\partial\bm{\varepsilon}} = \mathbb{L}_i : (\bm{\varepsilon}-\bm{\varepsilon}^{\rm t}_i) .
        \end{equation}
        Here, $\bm{\sigma}$ is the Cauchy stress tensor, $\bm{\varepsilon}=\nabla_{\rm s}\bm{u}=\frac12(\nabla\bm{u}+(\nabla\bm{u})^{\rm T})$ is the strain tensor, $\bm{\varepsilon}^{\rm t}_i$ is the transformation strain (known and constant within each phase), $\psi^0_i$ is the free energy in the stress-free state (so-called chemical energy), and $\mathbb{L}_i$ is a positive-definite fourth-order elastic moduli tensor possessing the usual symmetries. 
        
        The displacement field is assumed continuous at $\Gamma$, and $\Gamma$ is assumed to be smooth with $\bm{n}$ denoting the unit normal pointing from phase 1 to phase 2. Then, the classical compatibility conditions at $\Gamma$ read
        \begin{equation}
        \label{eq:compat}
        [\![\bm{\varepsilon}]\!] = \frac12 (\bm{c}\otimes\bm{n} + \bm{n}\otimes\bm{c}) , \qquad
        [\![\bm{\sigma}]\!]\cdot\bm{n} = \bm{0} ,
        \end{equation}
        where $[\![\mysquare]\!]=(\mysquare)_2-(\mysquare)_1$ denotes the jump at the interface and $\bm{c}$ is a vector. 
        The compatibility condition~\eqref{eq:compat}${}_1$ is a consequence of the continuity of the displacement at $\Gamma$. 
        On the other hand, mechanical equilibrium implies that Eq.~\eqref{eq:compat}${}_2$ holds, in addition to the usual equilibrium equation,
        \begin{equation}
        \label{eq:divS}
        \dive\bm{\sigma} = \bm{0} \quad \mbox{in} \;\, \Omega\setminus\Gamma ,
        \end{equation}
        supplemented by the boundary conditions,
        \begin{equation}
        \label{eq:BSs}
        \bm{u} = \bm{u}^\ast \quad \mbox{on} \;\, \partial\Omega_{\bm{u}} , \qquad
        \bm{\sigma} \cdot \bm{\nu} = \bm{t}^\ast \quad \mbox{on} \;\, \partial\Omega_{\bm{t}} ,
        \end{equation}
        where $\bm{u}^\ast$ and $\bm{t}^\ast$ denote the displacement and the surface traction prescribed on $\partial\Omega_{\bm{u}}$ and $\partial\Omega_{\bm{t}}$, respectively, and $\bm{\nu}$ denotes the unit outer normal to $\Omega$.
        
        Consider now $\Gamma$ to be a phase transformation front, so that $\Gamma$ is a moving interface, and assume the following (viscous) kinetic law for the interface \citep{Gurtin2000}
        \begin{equation}
        \label{eq:kinetic}
        \hat{v}_n = \hat{m} \hat{f} , \qquad
        \hat{f} = [\![ \psi ]\!] - \tilde{\bm{\sigma}} : [\![ \bm{\varepsilon} ]\!] + \gamma \kappa ,
        \end{equation}
        where $\hat{v}_n$ is the interface speed in the direction $\bm{n}$, $\hat{m}$ is the interface mobility, $\hat{f}$ is the thermodynamic driving force, $\gamma$ is the interfacial energy (energy density per unit area of the interface), and $\kappa$ is the total curvature (i.e., twice the mean curvature) of $\Gamma$. 
        Moreover, $\tilde{\bm{\sigma}}=\frac12(\bm{\sigma}_1+\bm{\sigma}_2)$ denotes the average stress at the interface; however, in view of the compatibility conditions~\eqref{eq:compat}, we have $\tilde{\bm{\sigma}}:[\![\bm{\varepsilon}]\!]=\bm{\sigma}_1:[\![\bm{\varepsilon}]\!]=\bm{\sigma}_2:[\![\bm{\varepsilon}]\!]$, hence $\tilde{\bm{\sigma}}$ in Eq.~\eqref{eq:kinetic}${}_2$ can be replaced by $\bm{\sigma}_1$ or $\bm{\sigma}_2$.
        Quantities referring to the sharp interface are here denoted by a symbol with a hat.
        
        Direct computational (e.g., finite-element) treatment of the microstructure evolution problem within the sharp-interface framework is not straightforward because it would require remeshing to ensure that the mesh matches the evolving microstructure at each time instant (and possibly at each iteration), see \citep{Moes2023} for a recent development towards this goal. 
        A feasible approach to model moving interfaces is to adopt a diffuse-interface approximation independent of the spatial discretization, which is discussed next.
    \subsection{Diffuse-interface framework: the phase-field method}\label{sec:diffuse}
        The sharp-interface problem discussed in the preceding section can be approximated by adopting a diffuse-interface description using the conventional phase-field method \citep{Chen2002,Steinbach2009,LevitasPreston2002part1,ClaytonKnap2011,HildebrandMiehe2012}. 
        To this end, a continuous order parameter $\phi$, $0\leq\phi\leq1$, is introduced such that $\phi=0$ corresponds to phase 1 and $\phi=1$ to phase~2, while the intermediate values 
        correspond to the diffuse interface. 
        Accordingly, the Helmholtz free energy $\psi$, characterizing both the individual phases and the diffuse interface, is assumed to be a function of the strain $\bm{\varepsilon}$ and of the order parameter $\phi$ and its spatial gradient $\nabla\phi$, and comprises the bulk and interfacial energy contributions, $\psi_{\rm bulk}$ and $\psi_{\rm int}$, respectively,
        \begin{equation}
        \label{eq:psiPF}
        \psi(\bm{\varepsilon},\phi,\nabla\phi) = \psi_{\rm bulk}(\bm{\varepsilon},\phi) + \psi_{\rm int}(\phi,\nabla\phi) .
        \end{equation}
        
        The bulk energy contribution $\psi_{\rm bulk}$ is specified as
        \begin{equation}
        \label{eq:psiBulk}
        \psi_{\rm bulk}(\bm{\varepsilon},\phi) = \psi^0(\phi) +
          \frac12 (\bm{\varepsilon}-\bm{\varepsilon}^{\rm t}(\phi)) : \mathbb{L}(\phi) : (\bm{\varepsilon}-\bm{\varepsilon}^{\rm t}(\phi)) ,
        \end{equation}
        where
        \begin{equation}
        \psi^0(\phi) = (1-h(\phi)) \psi^0_1 + h(\phi) \psi^0_2 , \quad
        \bm{\varepsilon}^{\rm t}(\phi)= (1-h(\phi))\bm{\varepsilon}^{\rm t}_1 + h(\phi)\bm{\varepsilon}^{\rm t}_2 , \quad
        \mathbb{L}(\phi)= (1-h(\phi))\mathbb{L}_1 + h(\phi)\mathbb{L}_2 .
        \end{equation}
        Here, the chemical energy $\psi^0$, the transformation strain $\bm{\varepsilon}^{\rm t}$ and the elastic stiffness tensor $\mathbb{L}$ are interpolated between the respective parameters of the pure phases and depend on the order parameter $\phi$ through the weighting function $h(\phi)$ such that $h(0)=0$, $h(1)=1$, and $h'(0)=h'(1)=0$. 
        We adopt here $h(\phi)=3\phi^2-2\phi^3$, i.e.\ the popular 2-3-4 polynomial (with $\alpha=3$), see \citep{LevitasPreston2002part1,Steinbach2009}. 
        
        The interfacial energy contribution $\psi_{\rm int}$ is adopted in the form of the so-called double-well potential \citep{Steinbach2009},
        \begin{equation}
        \psi_{\rm int}(\phi,\nabla\phi) = \frac{6\gamma}{\ell} \left( \phi^2(1-\phi)^2 + \frac{\ell^2}{4} \nabla\phi \cdot \nabla\phi \right) ,
        \end{equation}
        where $\gamma$ is the interfacial energy, see Eq.~\eqref{eq:kinetic}${}_2$, and $\ell$ is the interface thickness parameter. 
        
        Note that, when $\phi=0$ , the bulk energy $\psi_{\rm bulk}$, Eq.~\eqref{eq:psiBulk}, reduces to $\psi_1$, i.e.\ to the free energy specified by Eq.~\eqref{eq:psiSharp} for the pure phase 1, and likewise for the free energy $\psi_2$ of the pure phase 2 when $\phi=1$.
        
        In the phase-field method, the interfaces are treated as diffuse hence the discontinuities at the interfaces are smeared out, and the compatibility conditions \eqref{eq:compat} need not be considered. Since the strains and stresses are continuous, the mechanical equilibrium equation holds in the entire domain $\Omega$, i.e.\ also within the diffuse interface, 
        \begin{equation}
        \label{eq:equil:pf}
        \dive\bm{\sigma} = \bm{0} \qquad \mbox{in} \;\, \Omega ,
        \end{equation}
        with $\bm{\sigma}=\partial\psi/\partial\bm{\varepsilon}=\mathbb{L}(\phi):(\bm{\varepsilon}-\bm{\varepsilon}^{\rm t}(\phi))$, subject to the boundary conditions~\eqref{eq:BSs}.
        
        Evolution of the microstructure is governed by the time-dependent Ginzburg--Landau equation \citep{PenroseFife1990,Chen2002},
        \begin{equation}
        \label{eq:GL}
        \dot{\phi} = - m \frac{\delta\Psi}{\delta\phi} , \qquad
        \frac{\delta\Psi}{\delta\phi} = \frac{\partial\psi}{\partial\phi} - \nabla \cdot \frac{\partial\psi}{\partial\nabla\phi} ,
        \end{equation}
        where $\delta\Psi/\delta\phi$ is the variational derivative of the total free energy functional $\Psi=\int_\Omega\psi\rd V$, and $m$ is the mobility parameter. 
        Homogeneous Neumann boundary conditions are assumed on the entire boundary,
        \begin{equation}
        \nabla\phi \cdot \bm{\nu} = 0 \quad \mbox{on} \;\, \partial\Omega .
        \end{equation}
        
        It can be shown that the equilibrium profile of the diffuse interface has the following form \citep{HildebrandMiehe2012},
        \begin{equation}
          \label{eq:profile}
          \phi(\xi) = \frac12 \tanh \left( \frac{\xi-\xi_0}{\ell} \right) + \frac12 ,
        \end{equation}
        which can be obtained by minimizing the total interfacial energy $\Psi_{\rm int}=\int_\Omega\psi_{\rm int}\rd V$ of a planar interface located at $\xi=\xi_0$ (i.e.\ $\phi(\xi_0)=\frac12$), where $\xi$ denotes the coordinate in the direction normal to the interface. 
        This profile is only an approximation of the actual profile when the interface is non-planar or when the interface propagates due to a mechanical driving force that is not constant within the interface. In fact, special conditions must be met so that a travelling wave solution is obtained with the interface profile independent of the velocity \citep{Steinbach2009}. 
        
        Adopting the equilibrium profile~\eqref{eq:profile} as an approximation of the actual profile, the mobility parameter $m$ in the evolution equation~\eqref{eq:GL} can be related to the mobility parameter $\hat{m}$ in the sharp-interface description, cf.\ Eq.~\eqref{eq:kinetic}${}_1$. By equating the energy dissipated at the propagating sharp interface to that dissipated within the diffuse interface propagating with the same velocity, thus $\hat{f}\hat{v}_n=\hat{v}_n^2/\hat{m}=\int_{-\infty}^{+\infty}(\dot{\phi}^2/m)\rd\xi$, where $\xi_0=\hat{v}_n t$, the following relationship is found,
        \begin{equation}
          \label{eq:mobility}
          \hat{m} = 3 m \ell .
        \end{equation}
        It follows that when the interface parameter $\ell$ is varied, then the mobility parameter $m$ must be scaled according to Eq.~\eqref{eq:mobility} so that the effective mobility $\hat{m}$ of the interface is not affected (to the first order, in view of the approximate profile of the diffuse interface).
    \subsection{Finite-element treatment of the phase-field model}\label{sec:PFM}
        In order to arrive at the finite-element formulation, the mechanical equilibrium equation~\eqref{eq:equil:pf} is expressed in a weak form by following the standard procedure,
        \begin{equation}
        \label{eq:VWP}
        \int_\Omega \bm{\sigma} : \delta\bm{\varepsilon} \, \rd V
          - \int_{\partial\Omega_{\bm{t}}} \bm{t}^\ast \cdot \delta\bm{u} \, \rd S = 0
          \qquad \forall \, \delta\bm{u} ,
        \end{equation}
        where $\delta\bm{\varepsilon}=\nabla_{\rm s}\delta\bm{u}$, and $\delta\bm{u}$ is the virtual displacement (test function) that vanishes, $\delta\bm{u}=\bm{0}$, on $\partial\Omega_{\bm{u}}$. 
        Likewise, the evolution equation~\eqref{eq:GL} is expressed in a weak form,
        \begin{equation}
        \label{eq:GL:weak}
        \int_\Omega \Big[ \Big( \frac{\dot{\phi}}{m} + \frac{\partial\psi_{\rm bulk}}{\partial\phi} + \frac{\partial\psi_{\rm int}}{\partial\phi} \Big) \delta\phi +
          \frac{\partial\psi_{\rm int}}{\partial\nabla\phi} \cdot \nabla \delta\phi \Big] \rd V = 0
          \qquad \forall \, \delta\phi ,
        \end{equation}
        where $\delta\phi$ is the respective test function and
        \begin{equation}
        \label{eq:dpsi1}
        \frac{\partial\psi_{\rm bulk}}{\partial\phi} = \Big( \Delta\phi^0
          - \bm{\sigma} : \Delta\bm{\varepsilon}^{\rm t}
          + \frac12 (\bm{\varepsilon}-\bm{\varepsilon}^{\rm t}(\phi)) : \Delta\mathbb{L} :
            (\bm{\varepsilon}-\bm{\varepsilon}^{\rm t}(\phi))
          \Big) h'(\phi) ,
        \end{equation}
        \begin{equation}
        \label{eq:dpsi2}
        \frac{\partial\psi_{\rm int}}{\partial\phi} = \frac{12\gamma}{\ell} \phi(1-\phi)(1-2\phi) , \qquad
        \frac{\partial\psi_{\rm int}}{\partial\nabla\phi} = 3 \gamma \ell \, \nabla\phi ,
        \end{equation}
        where $\Delta(\mysquare)=(\mysquare)_2-(\mysquare)_1$. 
        
        The incremental (time-discrete) form of the evolution equation~\eqref{eq:GL:weak} is obtained by applying the implicit backward-Euler time integration scheme. Considering a typical time increment $t_n\to t_{n+1}=t_n+\tau$, the rate of $\phi$ is thus approximated by the finite difference,
        \begin{equation}
        \label{eq:dphi}
        \dot{\phi} \approx \frac{1}{\tau}(\phi-\phi_n) ,
        \end{equation}
        where $\phi_n$ is the known value at the previous time step $t_n$, and $\phi=\phi_{n+1}$ is the unknown value at the current time step. 
        Here and below, to make the notation more compact, the subscript $n+1$ is omitted for the quantities evaluated at the current time step $t_{n+1}$. 
        Adopting this notation, the time-discrete problem at $t_{n+1}$ is given directly by the weak forms~\eqref{eq:VWP} and~\eqref{eq:GL:weak} with $\dot{\phi}$ in Eq.~\ref{eq:GL:weak} replaced by its finite-difference approximation~\eqref{eq:dphi}. 
        
        Spatial discretization is performed using the finite element method. The fields of displacement $\bm{u}$ and order parameter $\phi$ are approximated using the respective basis functions $N_k^{(\bm{u})}$ and $N_k^{(\phi)}$,
        \begin{equation}
        \label{eq:FE:appr}
        \bm{u}^h = \sum_k N_k^{(\bm{u})} \bm{u}_k , \qquad
        \phi^h = \sum_k N_k^{(\phi)} \phi_k ,
        \end{equation}
        where $\bm{u}_k$ and $\phi_k$ denote the respective nodal values with $\fem{p}{}_{\bm{u}}=\{\bm{u}_k\}$ and $\fem{p}{}_{\phi}=\{\phi_k\}$ denoting the global vectors.
        Following the standard approach, $\bm{u}^h$ and $\phi^h$ are introduced into the weak forms~\eqref{eq:VWP} and~\eqref{eq:GL:weak}, and integration is performed over individual finite elements $\omega$ of the triangulation ${\cal T}$ of the domain $\Omega$,
        \begin{equation}
        \label{eq:VWP:FE}
        \sum_{\omega\in{\cal T}} \int_{\omega} \bm{\sigma} : \delta\bm{\varepsilon}^h \, \rd V
          - \sum_{\partial\omega\in{\cal S}_{\bm{t}}} \int_{\partial\omega} \bm{t}^\ast \cdot \delta\bm{u}^h \, \rd S = 0
          \qquad \forall \, \delta\fem{p}{}_{\bm{u}} ,
        \end{equation}
        \begin{equation}
        \label{eq:GL:weak:FE}
        \sum_{\omega\in{\cal T}}
          \int_{\omega} \Big[ \Big( \frac{\phi^h-\phi^h_n}{m\tau} + \frac{\partial\psi_{\rm bulk}}{\partial\phi} + \frac{\partial\psi_{\rm int}}{\partial\phi} \Big) \delta\phi^h +
          \frac{\partial\psi_{\rm int}}{\partial\nabla\phi} \cdot \nabla \delta\phi^h \Big] \rd V = 0
          \qquad \forall \, \delta\fem{p}{}_{\phi} ,
        \end{equation}
        where $\delta\bm{\varepsilon}^h=\nabla_{\rm s}\delta\bm{u}^h$, and ${\cal S}_{\bm{t}}$ denotes the triangulation of the boundary $\partial\Omega_{\bm{t}}$ into surface segments $\partial\omega$, consistent with the triangulation ${\cal T}$ of the bulk.
        
        The discretized weak forms~\eqref{eq:VWP:FE} and~\eqref{eq:GL:weak:FE} define the set of coupled nonlinear equations that can be written in the following residual form,
        \begin{equation}
        \label{eq:RuRphi}
        \fem{R}_{\bm{u}}(\fem{p}{}_{\bm{u}};\fem{p}{}_{\phi}) = \fem{0} , \qquad
        \fem{R}_{\phi}(\fem{p}{}_{\phi};\fem{p}{}_{\bm{u}}) = \fem{0} .
        \end{equation}
        At each time step of the incremental procedure, these equations are solved in a monolithic manner, thus $\fem{R}(\fem{p})=\fem{0}$ with $\fem{R}=\{\fem{R}_{\bm{u}},\fem{R}_{\phi}\}$ and $\fem{p}=\{\fem{p}{}_{\bm{u}},\fem{p}{}_{\phi}\}$, using the Newton method.
\section{LET-PF: hybrid diffuse--semisharp treatment of propagating interfaces}\label{sec:hybrid}
    In Section~\ref{sec:LET}, we briefly present the laminated element technique (LET) that has been recently developed by \citet{Dobrzanski2024} for an efficient treatment of material (i.e.\ non-propagating) interfaces. 
    The presentation is limited to the simplest case of small-strain framework and linear elasticity, but the approach is more general \citep{Dobrzanski2024}.
    Next, in Section~\ref{sec:LET-PF}, LET is combined with the phase-field method thus leading to a hybrid diffuse--semisharp treatment of propagating interfaces. 
    \subsection{Laminated element technique (LET) for material interfaces}\label{sec:LET}
\begin{figure}
            \centerline{
            \begin{tabular}{ccc} 
                \includegraphics[width=0.55\textwidth]{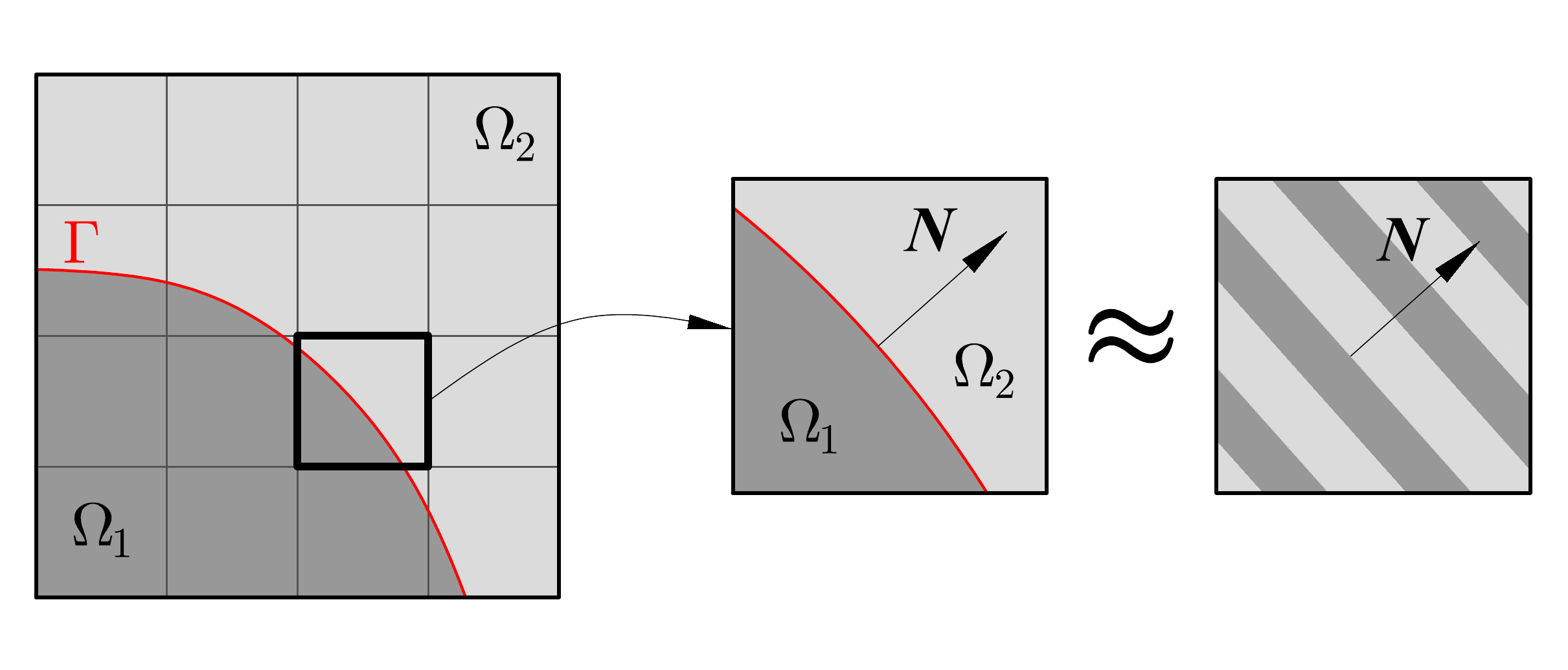} & \hspace{1eM} &
                \raisebox{0.5ex}{\includegraphics[width=0.21\textwidth]{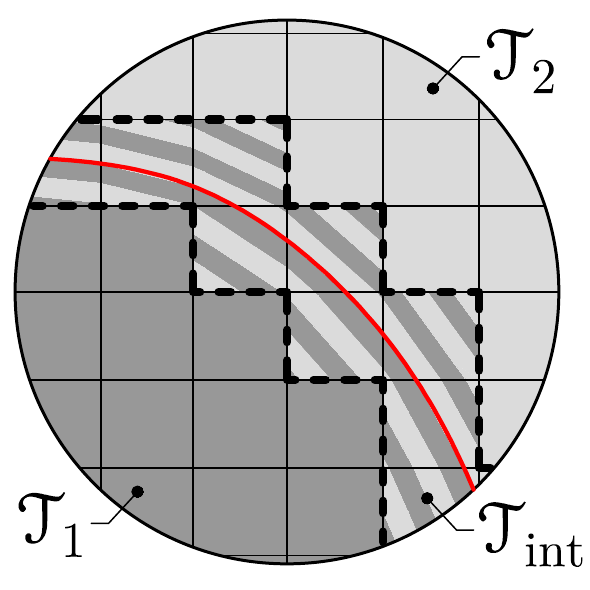}} \\
                {\footnotesize (a)} &  & {\footnotesize (b)}
            \end{tabular}
            }
            \caption{Laminated element technique (LET): (a) the elements cut by an interface are treated as simple laminates; (b)~the laminated elements are grouped in subset ${\cal T}_{\rm int}$, the remaining elements belong to subset ${\cal T}_1$ or ${\cal T}_2$.}
            \label{fig:LET:sketch}
        \end{figure}

        In LET, a non-matching (non-conforming) finite-element discretization is considered, i.e., the mesh is not aligned with the interface. Each element that is cut by the interface is treated as a laminate composed of the two phases involved, and the corresponding laminated microstructure is fully characterized by the volume fraction of the phases and by the lamination orientation which are determined as discussed below. At the same time, no treatment is applied to the remaining elements that belong to one of the phases. 
        Considering that the interface affects (through lamination) only one layer of elements along the interface, see the sketch in Fig.~\ref{fig:LET:sketch}, LET can be classified as
        a \emph{semisharp} approach, as opposed to the sharp-interface approach, characteristic for matching-mesh discretizations and for XFEM-type approaches, and as opposed to the diffuse-interface approach, characteristic for the phase-field method.
        
        The interface is implicitly defined by a level set function $\phi$ defined such that $\phi<\frac12$ corresponds to phase $1$, $\phi>\frac12$ corresponds to phase $2$, and the level set $\phi=\frac12$ specifies the interface. 
        In the discretized setting, the level set function $\phi\approx\phi^h$ is defined in terms of the basis functions $N_i^{(\phi)}$ and respective nodal values $\phi_i$, as in Eq.~\eqref{eq:FE:appr}$_2$. Note, however, that the level set function is here prescribed and does not evolve since it describes a fixed geometry. This assumption will be relaxed in Section~\ref{sec:LET-PF}.
        
        For each element $\omega\in{\cal T}$, the volume fraction $\eta^{(\omega)}=\eta^{(\omega)}_2=1-\eta^{(\omega)}_1$ of phase $2$ within the element is determined according to
        \begin{equation}
        \label{eq:eta:el}
        \eta^{(\omega)} = \frac{\sum_{k=1}^{N_{\rm n}} \langle \phi_k^{(\omega)}-\frac12 \rangle}
          {\sum_{k=1}^{N_{\rm n}} | \phi_k^{(\omega)}-\frac12 |} ,
        \end{equation}
        where $\langle\mysquare\rangle=\frac12(\mysquare+|\mysquare|)$ denotes the Macauley brackets, and $N_n$ denotes the number of nodes per element. 
        The above definition is applicable to four-node quadrilateral elements in 2D and to eight-node hexahedral elements in 3D. 
        For an element cut by the interface, we have $0<\eta^{(\omega)}<1$. Otherwise, for an element that is not cut by the interface and thus fully belongs to one of the phases, we have either $\eta^{(\omega)}=0$ or $\eta^{(\omega)}=1$.
        
        
        The lamination orientation is specified by the normal to the interface and is defined as the normalized gradient of the level-set function evaluated at the element centre $\bm{x}_0^{(\omega)}$,
        \begin{equation}
        \label{eq:nLET}
        \bm{n}^{(\omega)} = \frac{\nabla\phi^h(\bm{x}_0^{(\omega)})}{\|\nabla\phi^h(\bm{x}_0^{(\omega)})\|} .
        \end{equation}
        
        Now, the set of all elements can be split into three disjoint subsets, ${\cal T}={\cal T}_1\cup{\cal T}_2\cup{\cal T}_{\rm int}$, according to the value of the volume fraction $\eta^{(\omega)}$ in each element,
        \begin{equation}
        {\cal T}_1 = \{ \omega: \; \eta^{(\omega)}=0 \} , \qquad
        {\cal T}_2 = \{ \omega: \; \eta^{(\omega)}=1 \} , \qquad
        {\cal T}_{\rm int} = \{ \omega: \; 0<\eta^{(\omega)}<1 \} ,
        \end{equation}
        so that the elements that are cut by the interface are grouped in ${\cal T}_{\rm int}$, while sets ${\cal T}_1$ and ${\cal T}_2$ correspond to pure phases 1 and 2, respectively, see Fig.~\ref{fig:LET:sketch}(b). 
        
        For a fixed microstructure specified by a prescribed level-set function $\phi^h$, the problem at hand is a linear elasticity problem governed by the following weak form,
        \begin{equation}
        \label{eq:VWP:LET}
        \sum_{i=1}^2 \left( \sum_{\omega\in{\cal T}_i} \int_{\omega} \bm{\sigma}_i : \delta\bm{\varepsilon}^h \, \rd V \right)
          + \sum_{\omega\in{\cal T}_{\rm int}} \int_{\omega} \bar{\bm{\sigma}} : \delta\bm{\varepsilon}^h \, \rd V
          - \sum_{\partial\omega\in{\cal S}_{\bm{t}}} \int_{\partial\omega} \bm{t}^\ast \cdot \delta\bm{u}^h \, \rd S = 0
          \qquad \forall \, \delta\fem{p}{}_{\bm{u}} ,
        \end{equation}
        where $\bm{\sigma}_i$ denotes the stress in phase $i$ ($i=1,2$),
        \begin{equation}
        \bm{\sigma}_i = \bm{\sigma}_i(\bm{\varepsilon}^h) = \left. \frac{\partial \psi_i(\bm{\varepsilon})}{\partial \bm{\varepsilon}} \right|_{\bm{\varepsilon}=\bm{\varepsilon}^h}
          = \mathbb{L}_i : (\bm{\varepsilon}^h-\bm{\varepsilon}_i^{\rm t}) ,
        \end{equation}
        and $\bar{\bm{\sigma}}$ is the overall (macroscopic) stress in the laminated microstructure in the elements containing the interface,
        \begin{equation}
        \label{eq:sigmaBar}
        \bar{\bm{\sigma}} = \bar{\bm{\sigma}}(\bm{\varepsilon}^h,\eta^{(\omega)},\bm{n}^{(\omega)})
          =  \left. \frac{\partial \bar{\psi}_{\rm bulk}(\bm{\varepsilon},\eta^{(\omega)},\bm{n}^{(\omega)})}{\partial \bm{\varepsilon}} \right|_{\bm{\varepsilon}=\bm{\varepsilon}^h}
          = \bar{\mathbb{L}} : (\bm{\varepsilon}^h-\bar{\bm{\varepsilon}}^{\rm t}) .
        \end{equation}
        Here, $\bar{\psi}_{\rm bulk}=\bar{\psi}_{\rm bulk}(\bm{\varepsilon},\eta^{(\omega)},\bm{n}^{(\omega)})$ is the overall free energy of the laminate characterized by $\eta^{(\omega)}$ and $\bm{n}^{(\omega)}$, 
        $\bar{\mathbb{L}}=\bar{\mathbb{L}}(\eta^{(\omega)},\bm{n}^{(\omega)})$ is the overall stiffness tensor, and $\bar{\bm{\varepsilon}}^{\rm t}=\bar{\bm{\varepsilon}}^{\rm t}(\eta^{(\omega)},\bm{n}^{(\omega)})$ is the overall inelastic (transformation) strain. 
        Recall that $\eta^{(\omega)}$ and $\bm{n}^{(\omega)}$ depend on the level-set function $\phi^h$, see Eqs.~\eqref{eq:eta:el} and~\eqref{eq:nLET}. 
        The micro-macro transition for the laminated microstructure is discussed in \ref{app:laminate}.
        

        \paragraph{Remark 1}
        The presentation above is limited to the constitutive framework of linear elasticity with eigenstrain. However, the approach is more general, essentially with no restrictions on the constitutive models of the phases. For instance, an application to finite-strain elastoplasticity is presented in~\citep{Dobrzanski2024}.

        \paragraph{Remark 2}
        The LET approach, as presented above, is limited to the case of only two phases and at most one interface within a single finite element. The latter assumption can be relaxed, for instance, by introducing a separate level-set function for each interface, as discussed in Section~4.4 in~\citep{Dobrzanski2024}, although 
        that approach may not directly generalize 
        to the case of propagating interfaces, discussed in the next section. 
        Consideration of several phases would, in particular, require an adequate treatment of triple points. 
        Possible scenarios to treat such cases, for instance, based on higher-rank laminates, are left for future studies.

    \subsection{Combining the phase-field method with the laminated element technique}\label{sec:LET-PF}
\begin{figure}
            \centerline{
            \includegraphics[width=1.1\textwidth]{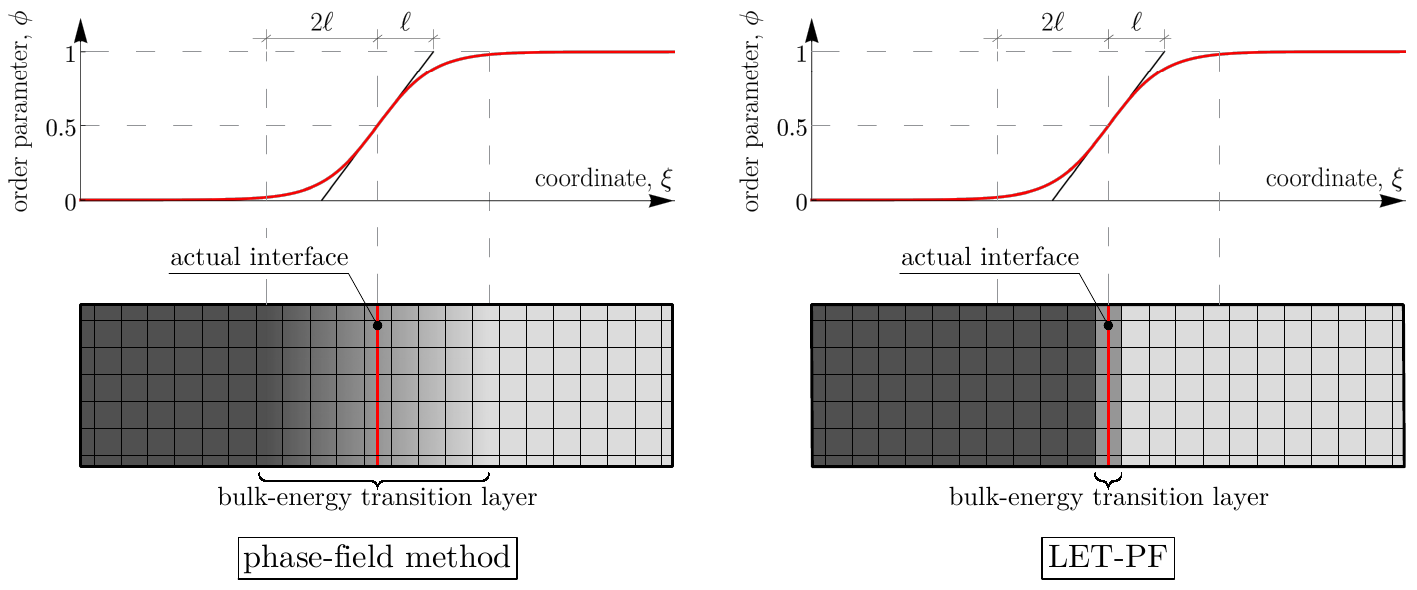}
            }
            \caption{LET-PF compared to the conventional phase-field method. The profile of the order parameter is sketched in the upper figures. Shading of the mesh in the bottom figures indicates the volume fraction of the phases that governs the bulk energy. In the conventional phase-field method, the transition layer is diffuse and spans several elements (left), while in the semisharp LET-PF method it is localized to only one layer of elements (right).
            }
            \label{fig:PFMvsLETPF}
        \end{figure}

        The idea of the proposed hybrid diffuse--semisharp approach is to combine the phase-field method with LET, thus the method will be referred to as LET-PF. 
        Specifically, diffuse interfaces, their propagation and the corresponding microstructure evolution are modelled in a manner similar to the phase-field method. At the same time, the order parameter $\phi$ plays the role of the level-set function in the semisharp LET-based treatment of the mechanical subproblem, see Fig.~\ref{fig:PFMvsLETPF}.
        
        
        As discussed in Section~\ref{sec:PFM}, the weak formulation of the discretized conventional phase-field method is specified by Eqs.~\eqref{eq:VWP:FE} and~\eqref{eq:GL:weak:FE}. 
        In LET-PF, the virtual work principle~\eqref{eq:VWP:FE} is replaced by that of the LET method, Eq.~\eqref{eq:VWP:LET}, which here refers to the current time instant $t=t_{n+1}$ with the microstructure specified by the 
        order parameter field $\phi^h$ at the current time $t_{n+1}$. 
        The time-discrete evolution of the order parameter $\phi^h$ is governed by the following weak form,
        \begin{equation}
        \label{eq:GL:weak:LET}
        \sum_{\omega\in{\cal T}}
          \int_{\omega} \Big[ \Big( \frac{\phi^h-\phi^h_n}{m\tau} + \frac{\partial\psi_{\rm int}}{\partial\phi} \Big) \delta\phi^h +
          \frac{\partial\psi_{\rm int}}{\partial\nabla\phi} \cdot \nabla \delta\phi^h \Big] \rd V +
        \sum_{\omega\in{\cal T}_{\rm int}} \int_{\omega} \sum_{i=1}^{N_{\rm n}}
          \frac{\partial \bar{\psi}_{\rm bulk}}{\partial\phi_i^{(\omega)}} \, \delta\phi_i^{(\omega)} \rd V = 0
          \qquad \forall \, \delta\fem{p}{}_{\phi} ,
        \end{equation}
        which replaces the weak form~\eqref{eq:GL:weak:FE} of the conventional phase-field model. 
        The first term in Eq.~\eqref{eq:GL:weak:LET} includes the interfacial energy contribution and the viscous evolution term and is 
        identical to the respective part of the weak form~\eqref{eq:GL:weak:FE}. 
        The second term in Eq.~\eqref{eq:GL:weak:LET} describes 
        the driving force resulting from the bulk energy $\bar{\psi}_{\rm bulk}$ of the laminated elements along the interface. 
        This term is the counterpart of the respective term in the LET form of the virtual work principle, see Eqs.~\eqref{eq:VWP:LET} and~\eqref{eq:sigmaBar}. 
        Note that, in LET-PF, the bulk energy depends on the order parameter field $\phi^h$ only in the laminated elements $\omega\in{\cal T}_{\rm int}$ (through $\eta^{(\omega)}$ and $\bm{n}^{(\omega)}$). 
        Actually, also the set ${\cal T}_{\rm int}$ itself depends on $\phi^h$ (likewise, ${\cal T}_1$ and ${\cal T}_2$) and may change during a time step and even between the individual iterations of the global Newton scheme. Accordingly, to improve the convergence behaviour, a regularization of the formula for $\eta^{(\omega)}$ in Eq.~\eqref{eq:eta:el} is employed, as discussed in Section~\ref{sec:FE}.
        
        Summarizing, the difference with respect to the conventional phase-field method is that, in LET-PF, the two sub-models are coupled only within the laminated interface elements in ${\cal T}_{\rm int}$, see Figs.~\ref{fig:LET:sketch} and~\ref{fig:PFMvsLETPF}.
        This is, in particular, apparent in Eq.~\eqref{eq:VWP:LET}, where $\bar{\bm{\sigma}}$ is the only term that depends on $\phi^h$, in addition to the set ${\cal T}_{\rm int}$ itself. 
        Likewise, in Eq.~\eqref{eq:GL:weak:LET}, the coupling occurs only through the overall bulk energy $\bar{\psi}_{\rm bulk}$ of the laminates that is defined only within the laminated interface elements in ${\cal T}_{\rm int}$.
        At the same time, in the conventional phase-field method, the coupling occurs within the whole volume of the diffuse interfaces, i.e., whenever $0<\phi<1$, see Fig.~\ref{fig:PFMvsLETPF}. 
        It follows that the mechanical equilibrium subproblem, which is governed by the bulk contribution to the free energy, is treated in LET-PF in a significantly sharper manner than in the conventional phase-field method, and hence coarser meshes can be used in LET-PF to achieve a similar accuracy, as illustrated by the numerical examples below.

    \subsection{Finite-element treatment of the LET-PF model}\label{sec:FE}
        The governing equations of LET-PF, see Eqs.~\eqref{eq:VWP:LET} and~\eqref{eq:GL:weak:LET}, define the set of non-linear equations to be solved at each time step. 
        The unknowns of the problem are the nodal values of the displacement ($\bm{u}_k$) and order parameter ($\phi_k$), both entering the formulation through the finite-element approximation~\eqref{eq:FE:appr} with $\fem{p}{}_{\bm{u}}$ and $\fem{p}{}_{\phi}$ denoting the respective global vectors of unknowns. 
        The general structure of the problem is the same as in the case of the conventional phase-field model, and the residual form of the governing equations of the two subproblems can be written as in Eq.~\eqref{eq:RuRphi}, of course, with due differences in the form of the residual vectors $\fem{R}_{\bm{u}}$ and $\fem{R}_\phi$. 
        A~monolithic Newton-based solution scheme is employed to solve these equations, just like in the case of the conventional phase-field model, see the end of Section~\ref{sec:PFM}.
        
        In the numerical examples, which are limited here to 2D, four-node quadrilateral elements are employed with bilinear shape functions used for both the displacement and the order parameter and with the standard $2\times2$ Gauss quadrature. 
        Computer implementation has been performed using \emph{AceGen}, a symbolic code generation system employing the automatic differentiation (AD) and expression optimization techniques \citep{Korelc2009,KorelcWriggers2016}, thus resulting in a computationally efficient computer code. The actual finite-element computations have been carried out using \emph{AceFEM}, a finite-element system integrated with \emph{AceGen}.
        
        Preliminary computations have revealed that the convergence of the Newton method at the global level is severely affected by the non-differentiable functions (the absolute value and Macauley brackets) that are involved in the formula for the volume fraction $\eta^{(\omega)}$ in element $\omega$, see Eq.~\eqref{eq:eta:el}. In order to improve the robustness of the method, the following regularization is applied to the absolute value function and Macauley brackets used in Eq.~\eqref{eq:eta:el},
        \begin{equation}\label{eq:eta:el:reg}
        |\phi|_{\rm reg} = \left\{
          \begin{array}{ll}
            |\phi| & \mbox{if}\; |\phi|\geq\phi_{\rm reg} , \\
            \frac{1}{2\phi_{\rm reg}}(\phi^2+\phi^2_{\rm reg}) \; & \mbox{if}\; |\phi|<\phi_{\rm reg} ,
          \end{array}
          \right.
          \qquad
          \langle\phi\rangle_{\rm reg} = \frac12 (\phi+|\phi|_{\rm reg}) ,
        \end{equation}
        where $\phi_{\rm reg}>0$ is a regularization parameter with the typical value of 0.1 used in the computations reported in Section~\ref{sec:numExamples}. 
        As illustrated in \ref{sec:regEffect}, the above regularization significantly improves the convergence behaviour so that the analysis can proceed with a larger time step, while it has a small effect on the results.
\section{Illustrative examples}\label{sec:numExamples}
    \subsection{Evolving circular inclusion}\label{sec:circular}
        The aim of this study is to evaluate the overall performance of LET-PF and to compare it to the conventional phase-field method (abbreviated to PFM in the sequel) for a 2D problem that has an analytical solution, to be used as a reference. 
        The sharp-interface problem is specified first, followed by the description of the computational model and by the results of a comprehensive study.
        \subsubsection{Sharp-interface benchmark problem}
\begin{figure}[htbp]
                \centerline{
                \begin{tabular}{ccc} 
                    \includegraphics[width=0.38\textwidth]{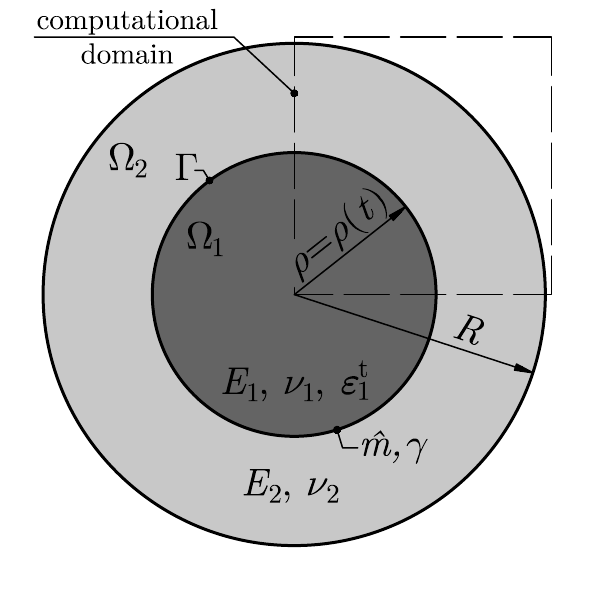} & &
                    \raisebox{.5em}[0pt][0pt]{\includegraphics[width=0.36\textwidth]{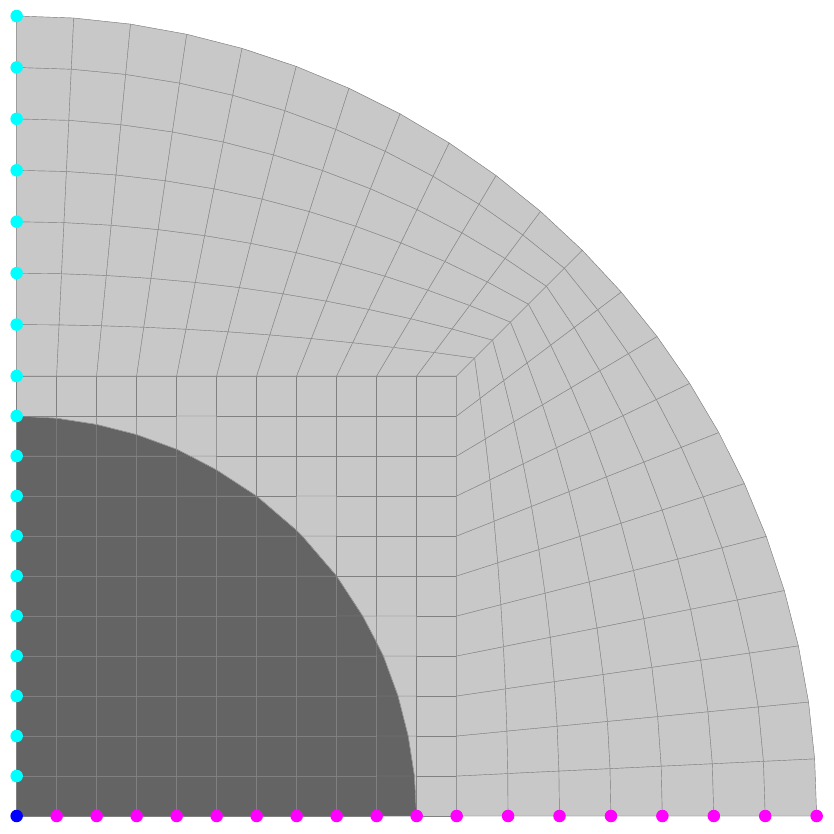}} \\
                    {\footnotesize (a)}~~ & ~~ & {\footnotesize (b)}
                \end{tabular}
                }
                \caption{Evolving circular inclusion: (a) scheme of the problem; (b) computational domain with a regular (non-conforming) mesh of quadrilateral elements. A coarse mesh of $11\times 11$ elements within the internal square part (element size $h=0.1$) is shown, the actual computations are carried out using significantly finer meshes.
                }
                \label{fig:inclusionscheme}
            \end{figure}

             We consider an elastic circular domain $\Omega_2$ (matrix) of the radius $R$ and an evolving elastic circular inclusion $\Omega_1$ of the radius $\rho=\rho(t)$, see Fig.~\ref{fig:inclusionscheme}(a). The phase transformation and thus the evolution of the inclusion are induced by the internal stresses resulting from the volumetric eigenstrain $\bm{\varepsilon}_1^\text{t}=\epsilon\bm{I}$ in the inclusion and by the energy of the interface $\Gamma$ of the density $\gamma$. Propagation of the interface is assumed to be governed by the viscous kinetic law, cf.\ Eq.~\eqref{eq:kinetic}${}_1$,
            \begin{equation}\label{eq:vn:rhodot}
                \hat{v}_n=-\dot{\rho} = \hat{m} \hat{f} , \qquad \hat{f} = \hat{f}_{\rm bulk}+\hat{f}_{\rm int} ,
            \end{equation}
            where $\hat{f}_{\rm bulk}$ and $\hat{f}_{\rm int}$ denote the local thermodynamic driving forces originating from the bulk and interfacial energy, respectively. The outer boundary is assumed to be free.
            %
            %
            %
            %
    
            The sharp-interface problem specified above is a 1D problem and admits an analytical solution that is derived in \ref{app:BimaterialLame}. In the special case of identical elastic properties of the two phases, 
            the evolution of the inclusion radius is governed by the following differential equation, cf.\ Eq.~\eqref{eq:rhodot:app},
            \begin{equation}\label{eq:BimaterialLame:rhodotwA}
                \dot{\rho} = -\frac{\hat{m} \gamma}{\rho}\left( 1 + A\frac{\rho}{\rho_0}\left( \frac{1}{2}-\left( \frac{\rho}{R} \right)^2 \right) \right) , \qquad
                A=\frac{\rho_0 E \epsilon^2}{(1-\nu^2)\gamma} ,
            \end{equation}
            where $\rho_0$ denotes the initial radius of the inclusion, $E$ is the Young's modulus and $\nu$ is the Poisson's ratio. Note that a dimensionless parameter $A$ has been introduced which combines all the problem parameters related to both the elastic strain energy and the interfacial energy (the chemical energy contribution is not considered, hence the elastic strain energy is the only contribution to the bulk energy). For small values of $A$, the evolution is mainly governed by the interfacial energy, and for large values of $A$ by the elastic strain energy.
            Exact integration of Eq.~\eqref{eq:BimaterialLame:rhodotwA} is not possible; integration is thus performed numerically.
    
            The explicit formulae for the driving force contributions $\hat{f}_{\rm bulk}$ and $\hat{f}_{\rm int}$ are given by Eq.~\eqref{eq:rhodot:forces}. Specifically, $\hat{f}_{\rm int}=\gamma/\rho$ is the driving force originating from the interfacial energy of the density $\gamma$, and the interface curvature is $1/\rho$. The driving force $\hat{f}_{\rm bulk}$ originates from the elastic strain energy induced by the inclusion eigenstrain. It is proportional to the elastic modulus $E$ and squared eigenstrain magnitude $\epsilon^2$ and depends on the geometrical factor $\rho/R$. For $\rho_0=R/2$, as assumed in the computations below, the driving force $\hat{f}=\hat{f}_{\rm bulk}+\hat{f}_{\rm int}$ is positive and hence $\dot{\rho}$ is negative, which means that the inclusion decreases in size and ultimately vanishes. 
        \subsubsection{Computational model}
            The two methods (LET-PF and PFM) involve identical material and numerical parameters, including those related to the finite-element model, and the specifications provided below apply to both methods.
            The finite-element computations are carried in 2D on a computational domain encompassing one quarter of the circular domain with adequate symmetry conditions imposed along the horizontal and vertical edges, see Fig.~\ref{fig:inclusionscheme}(b). 
            In the computations, it is assumed that the material properties of the inclusion and matrix are identical with the Young's modulus $E=1$ and Poisson's ratio $\nu=0.25$. 
            The radius of the domain is $R=2$, and the initial radius of the inclusion is $\rho_0=1$. 
            The volumetric eigenstrain of the inclusion is assumed as $\epsilon=0.1$. The effective mobility parameter is set to $\hat{m}=1$, and the actual mobility parameter $m$ used in the computations is then specified as a function of $\ell$ according to Eq.~\eqref{eq:mobility}.
    
\begin{figure}[htbp]
                \centerline{
                    \includegraphics[width=0.6\textwidth]{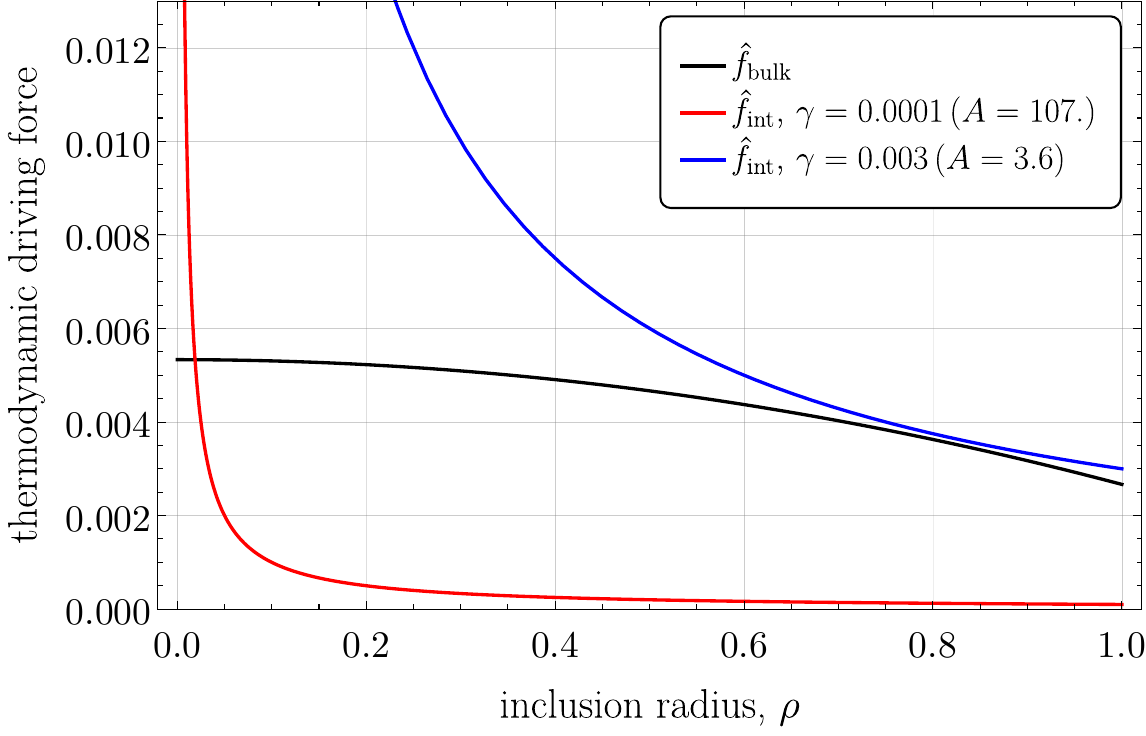} 
                }
                \caption{Thermodynamic driving forces $\hat{f}_{\rm bulk}$ and $\hat{f}_{\rm int}$ as a function of the inclusion radius $\rho$, see Eq.~\eqref{eq:rhodot:forces}.
                }
                \label{fig:regimes}
            \end{figure}

            One of the goals here is to study the performance of LET-PF over a wide range of the values of the interfacial energy density $\gamma$. 
            To this end, the values of $\gamma$ have been chosen from within the range $\gamma\in [0.0001, 0.003]$ (12 values have been chosen), so that the dimensionless parameter $A$ varies within the range ${A}\in [3.6,107.]$. 
            This range has been chosen such that the extreme values correspond to two regimes of interest. For $\gamma=0.0001$ ($A=107.$), the bulk contribution $\hat{f}_{\rm bulk}$ to the driving force is much greater than the interfacial contribution $\hat{f}_{\rm int}$ over nearly the whole range of inclusion radii $\rho$ (note that $\hat{f}_{\rm int}\to\infty$ for $\rho\to0$), see Fig.~\ref{fig:regimes}. Accordingly, in this (elasticity-driven) regime, evolution is mostly governed by the elastic strain (bulk) energy. 
            On the other hand, for $\gamma=0.003$ ($A=3.6$), both contributions are initially of the same order, see Fig.~\ref{fig:regimes}. 
            The third regime, in which $\hat{f}_{\rm int}\gg\hat{f}_{\rm bulk}$, is not interesting here because the evolution is then fully governed by the interfacial energy and the effect of LET-PF is negligible. Moreover, it has been checked that the response and the performance of both methods are then quite similar as for $\hat{f}_{\rm int}\approx\hat{f}_{\rm bulk}$.
            %
            %
            %
            %
            
            
            The impact of the interface thickness parameter $\ell$ on the accuracy of the results is also examined, with the  values of $\ell$ specified relative to the characteristic element size $h$, namely, $\ell/h\in\{0.75,1,1.5,2\}$. 
            As mentioned earlier, for the selected parameter values, the inclusion disappears completely during the evolution process by gradually reducing its radius to zero. Therefore, a regular grid of square finite elements is employed in the central region containing the inclusion, and quadrilateral elements with somewhat less regular shapes are employed in the remaining part of the domain, as depicted in Fig.~\ref{fig:inclusionscheme}(b). The main computations are performed for two mesh sizes, $h=0.02$ and $h=0.01$, which correspond to, respectively, $55\times55$ and $110\times110$ elements in the central region. 
            At the initial time $t=0$, the initial profile of the order parameter $\phi$ is prescribed according to formula \eqref{eq:profile} with $\xi$ replaced by the radius $r$ and $\xi_0$ replaced by the initial inclusion radius $\rho_0$.
    
            In all computations, an adaptive time incrementation scheme has been used, such that the time increment is increased or decreased depending on the current convergence behaviour. In each case, the maximum time increment $\Delta t_{\rm max}$ has been prescribed relative to $T_{\rm exact}$, the time of the complete evolution as resulting from the analytical solution. 
    
            The results will be presented in terms of the mean inclusion radius $\bar{\rho}$ which is defined as the average radius of the inclusion defined by the $\phi=\frac12$ level set. The individual radii are measured along a number of directions taken every 1 degree. 
        \subsubsection{Results}
\begin{figure}[htbp]
                \centerline{\scriptsize
                    \begin{tabular}{cc}
                        \includegraphics[width=0.45\textwidth]{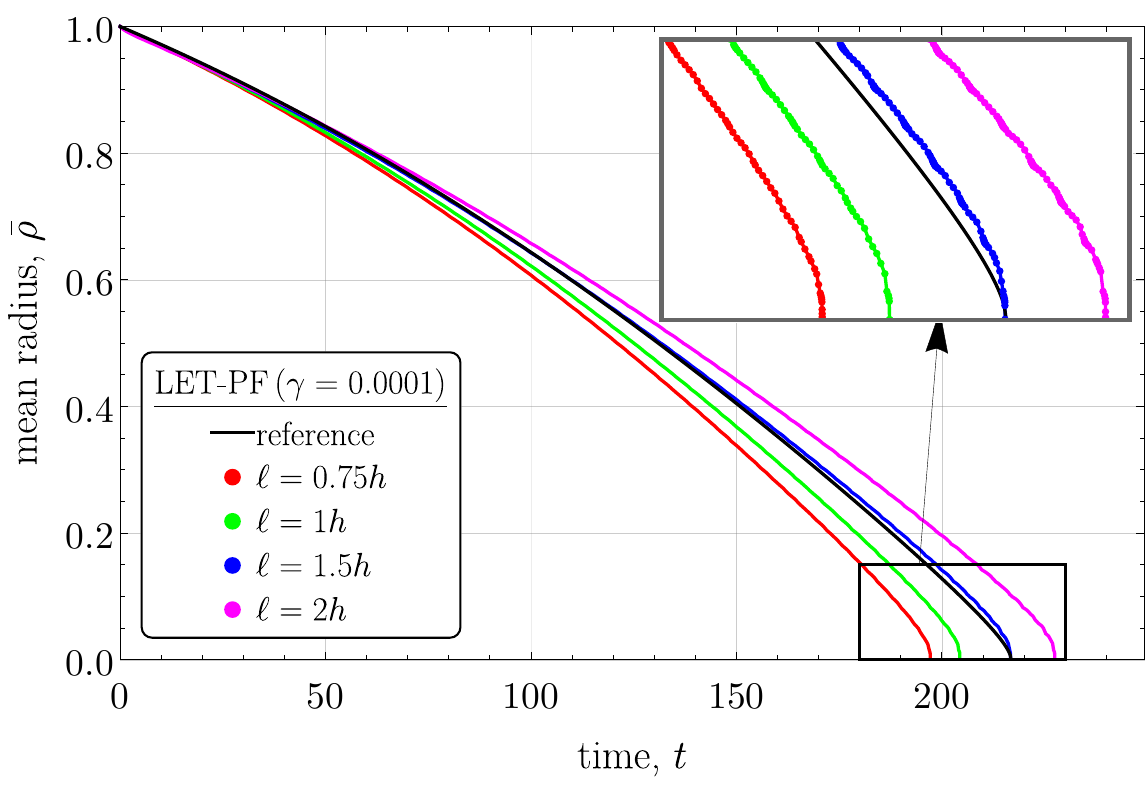} &
                        \includegraphics[width=0.45\textwidth]{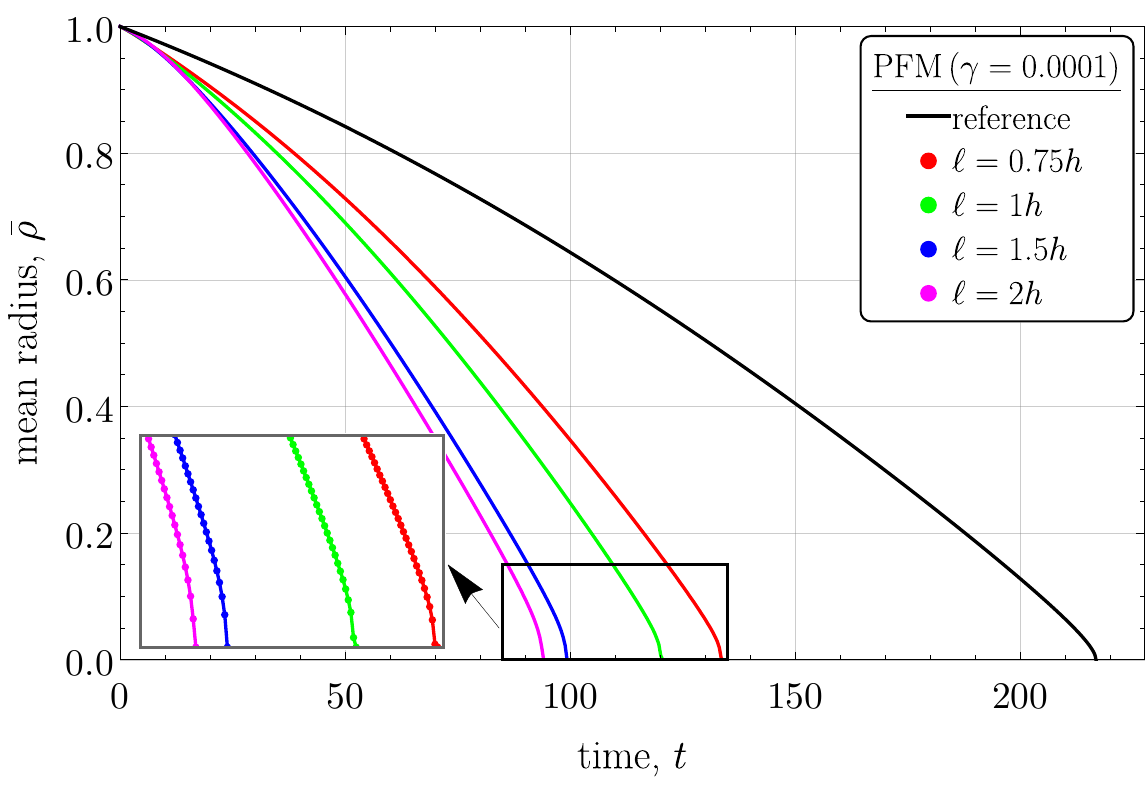} \\[1ex]
                        \hspace*{3em}(a) & \hspace*{3em}(b) \\[2ex]
                        \includegraphics[width=0.45\textwidth]{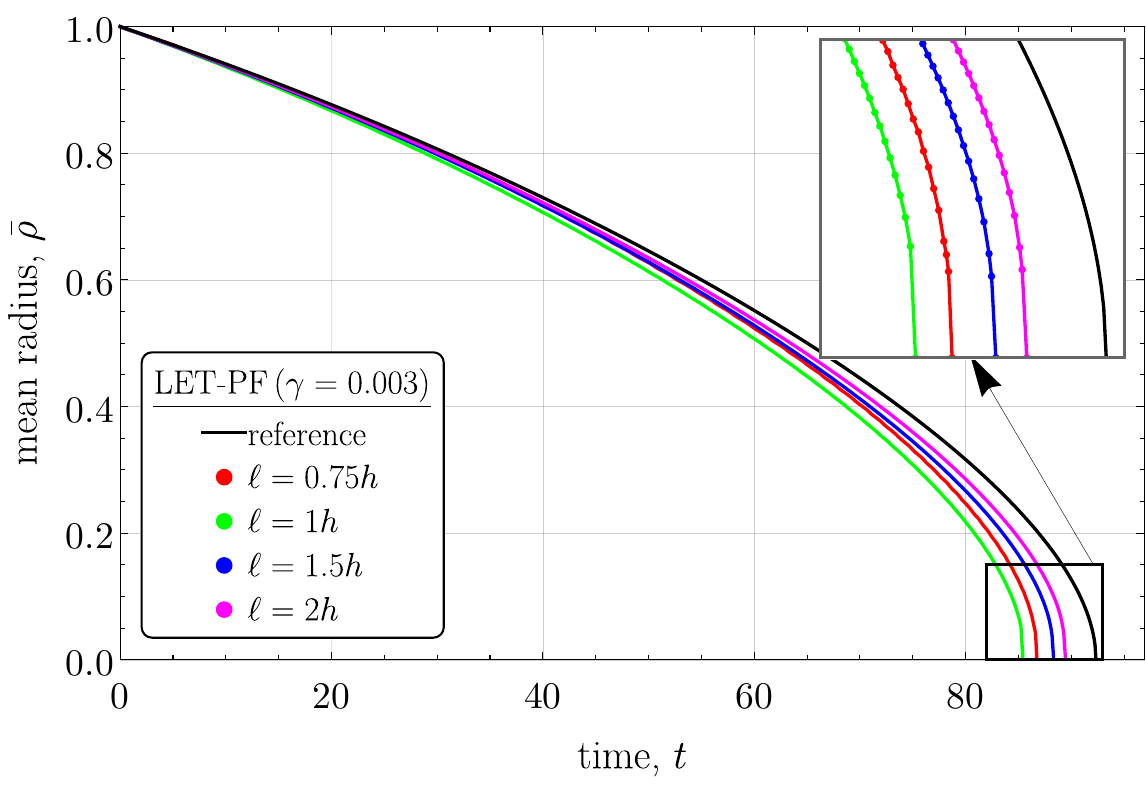} &
                        \includegraphics[width=0.45\textwidth]{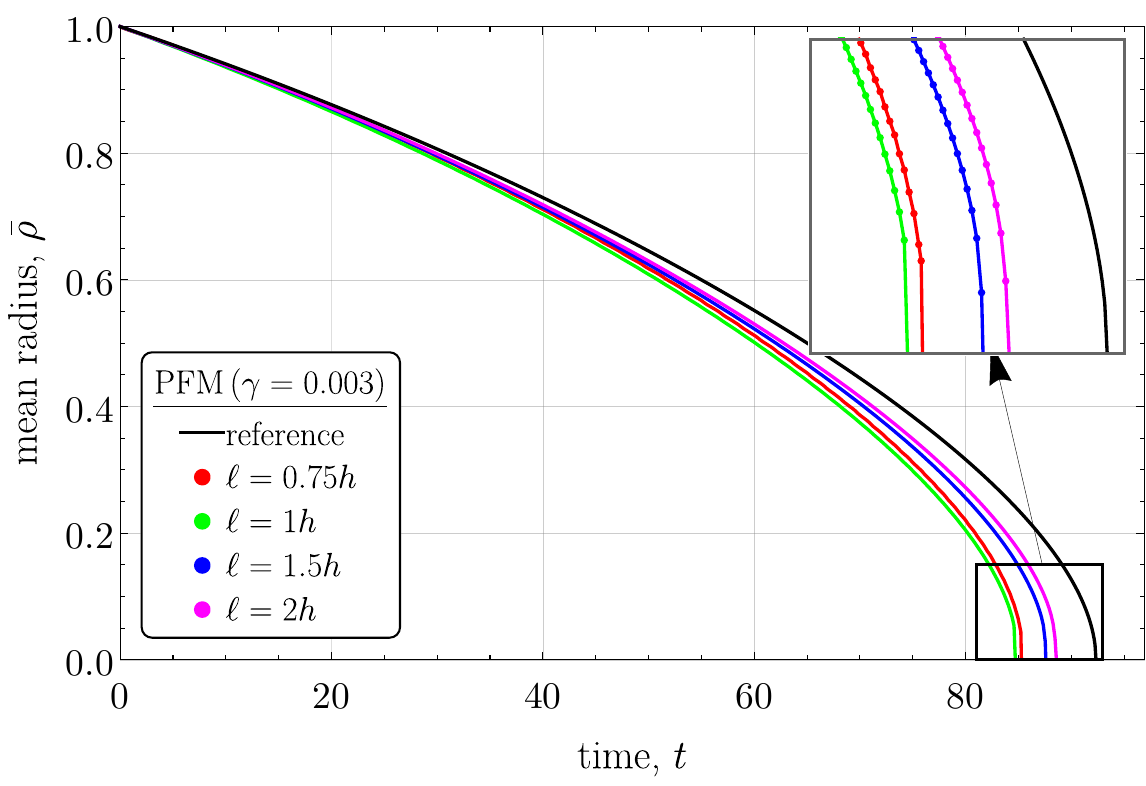} \\[1ex]
                        \hspace*{3em}(c) & \hspace*{3em}(d)
                    \end{tabular}
                    }
                \caption{Time evolution of the mean inclusion radius $\bar{\rho}$ predicted by LET-PF (a,c) and PFM (b,d) for $\gamma=0.0001$ (a,b) and for $\gamma=0.003$ (c,d). The computations are performed using a coarse mesh, $h=0.02$, and small time increment, $\Delta t_{\rm max}=T_{\rm exact}/500$. 
                As a reference, the analytical solution is depicted by a solid black line.}
                \label{fig:RadiusFullEvolution}
            \end{figure}

            The time evolution of the mean inclusion radius $\bar{\rho}$ is shown in Fig.~\ref{fig:RadiusFullEvolution} for the two extreme values of the interfacial energy, $\gamma=0.0001$ and $\gamma=0.03$. As a reference, the analytical solution is also included in Fig.~\ref{fig:RadiusFullEvolution}. 
            The computations have been performed with small time increments, $\Delta t_{\rm max}=T_{\rm exact}/500$, so that the error introduced by the time integration is insignificant. It can be seen that, in the elasticity-driven regime ($\gamma=0.0001$), PFM is highly inaccurate (Fig.~\ref{fig:RadiusFullEvolution}(b)), while LET-PF performs much better (Fig.~\ref{fig:RadiusFullEvolution}(a)). In the interfacial-energy-driven regime ($\gamma=0.003$), both methods perform similarly well. 
            Fig.~\ref{fig:RadiusFullEvolution} illustrates also a significant effect of the interface thickness parameter $\ell$.
            %
            %
            %
            %
    
\begin{figure}[htbp]
                \centerline{\scriptsize
                    \begin{tabular}{cccc}
                        \includegraphics[width=0.25\textwidth]{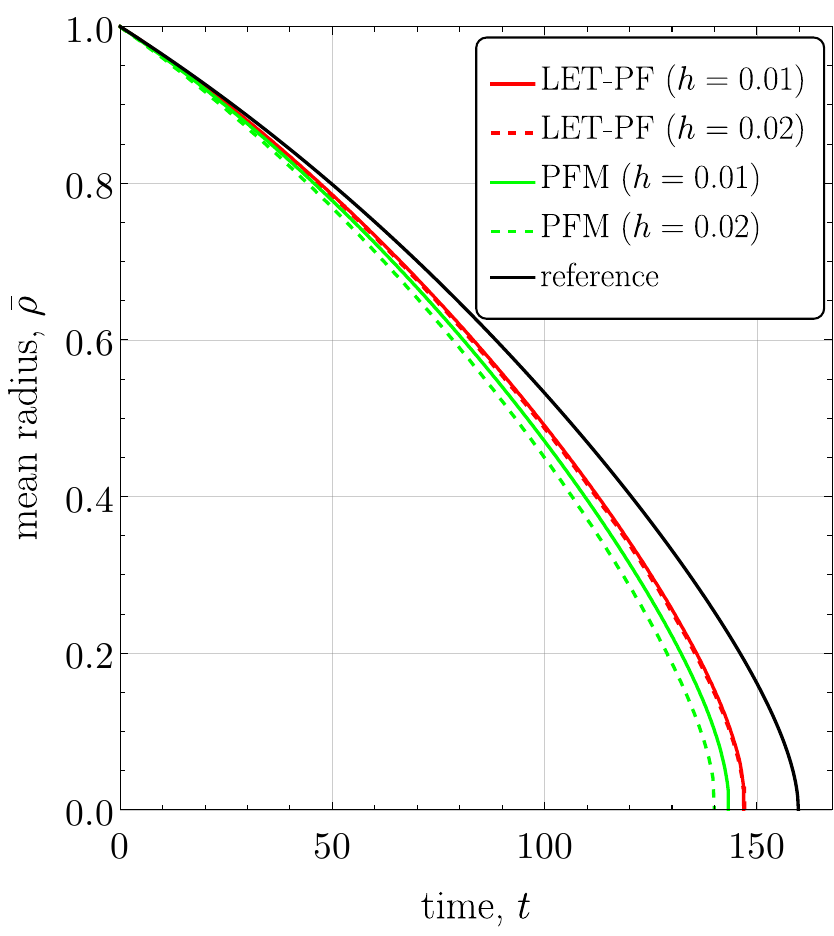} &
                        \includegraphics[width=0.25\textwidth]{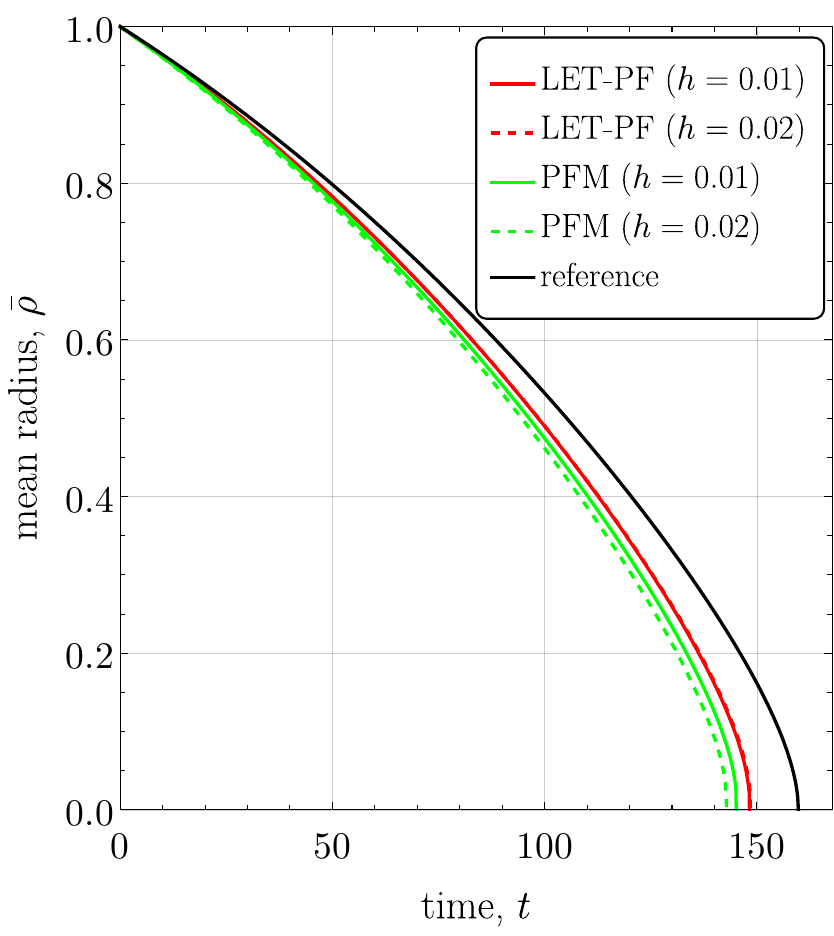} &
                        \includegraphics[width=0.25\textwidth]{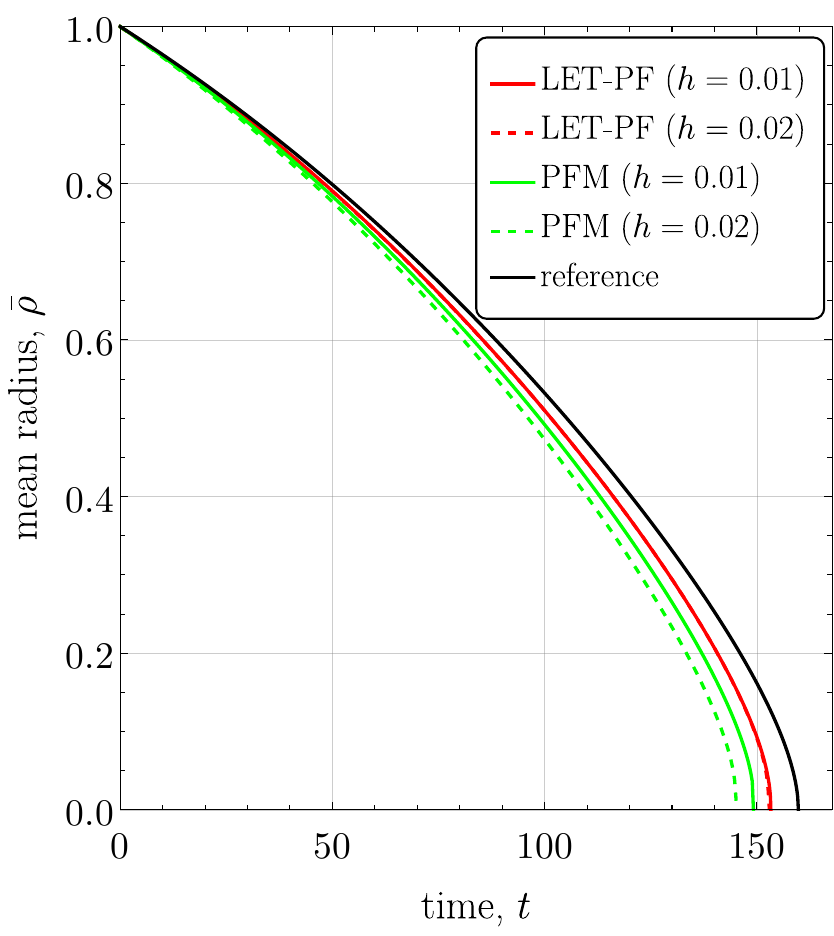} &
                        \includegraphics[width=0.25\textwidth]{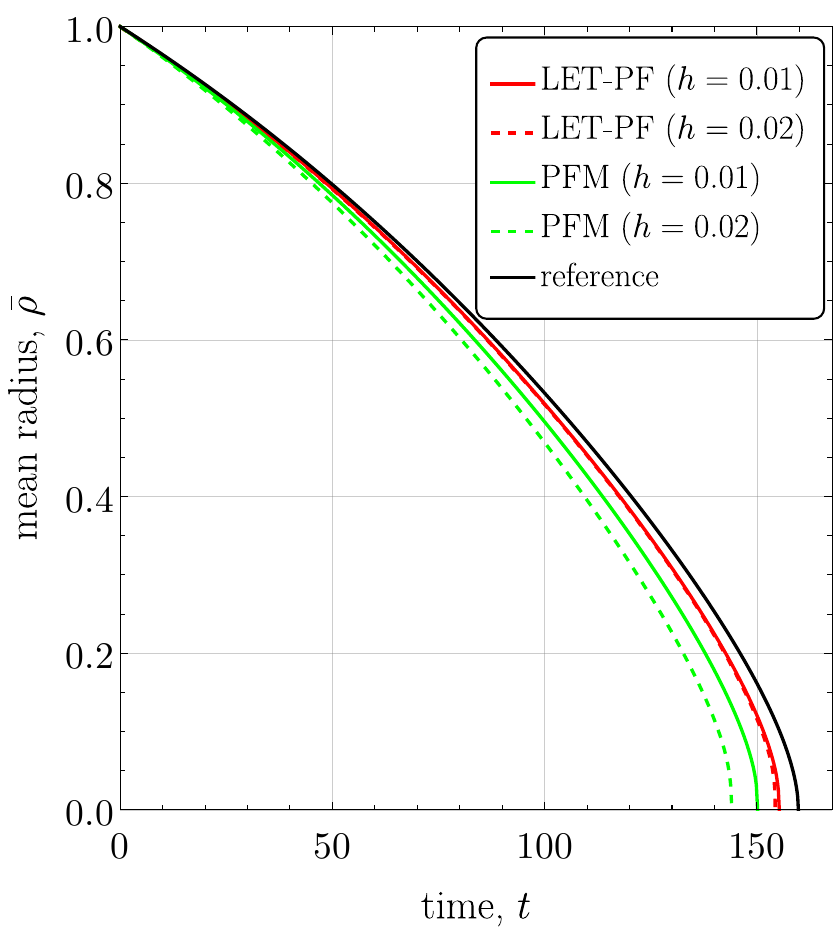} \\[0.5ex]
                        \includegraphics[width=0.25\textwidth]{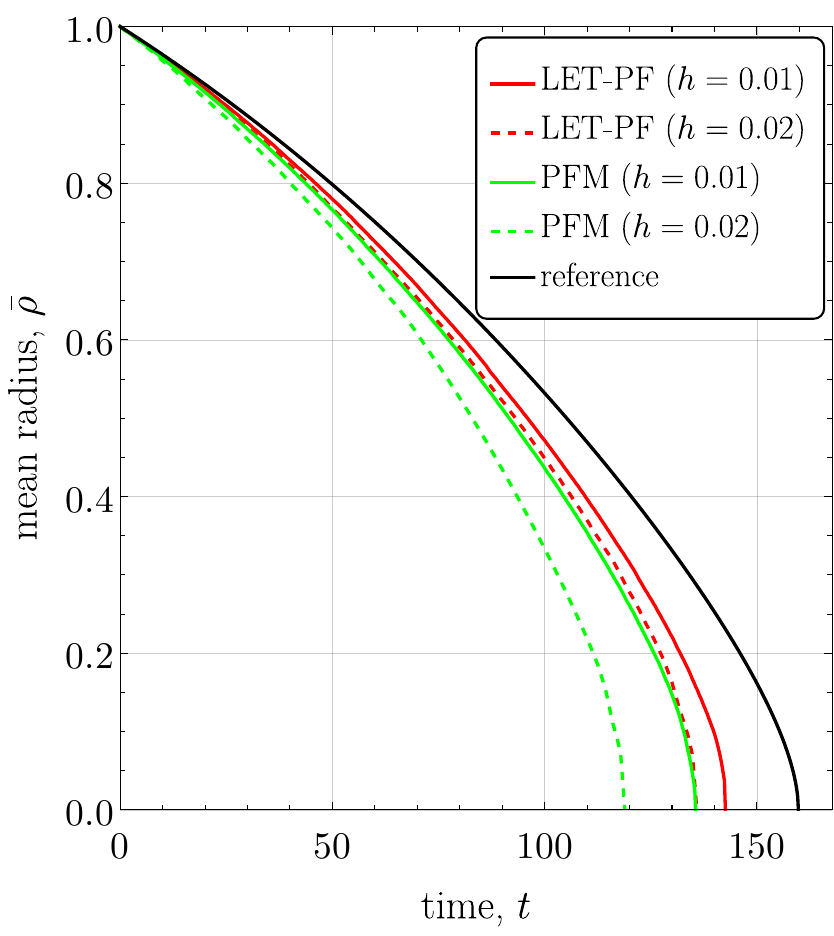} &
                        \includegraphics[width=0.25\textwidth]{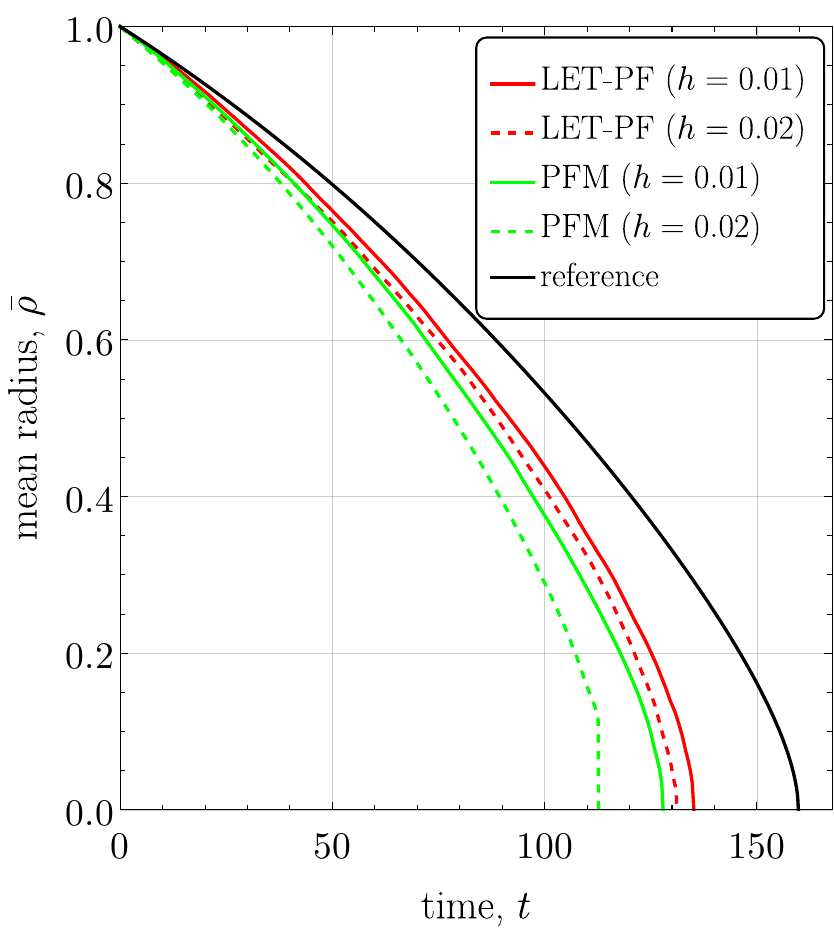} &
                        \includegraphics[width=0.25\textwidth]{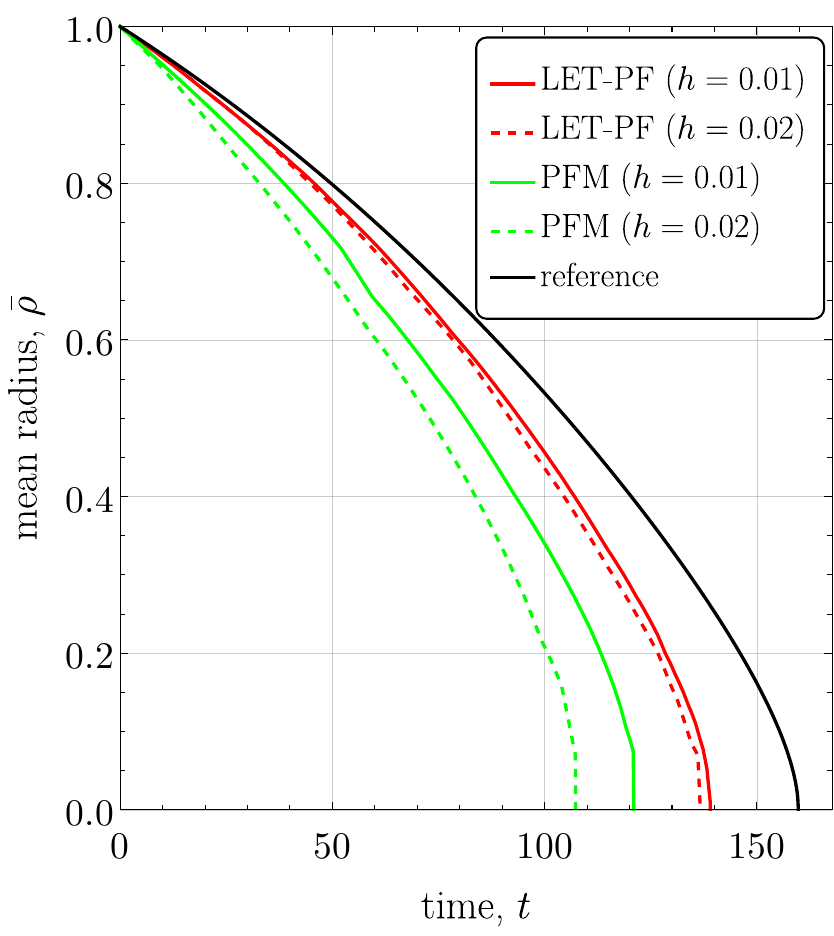} &
                        \includegraphics[width=0.25\textwidth]{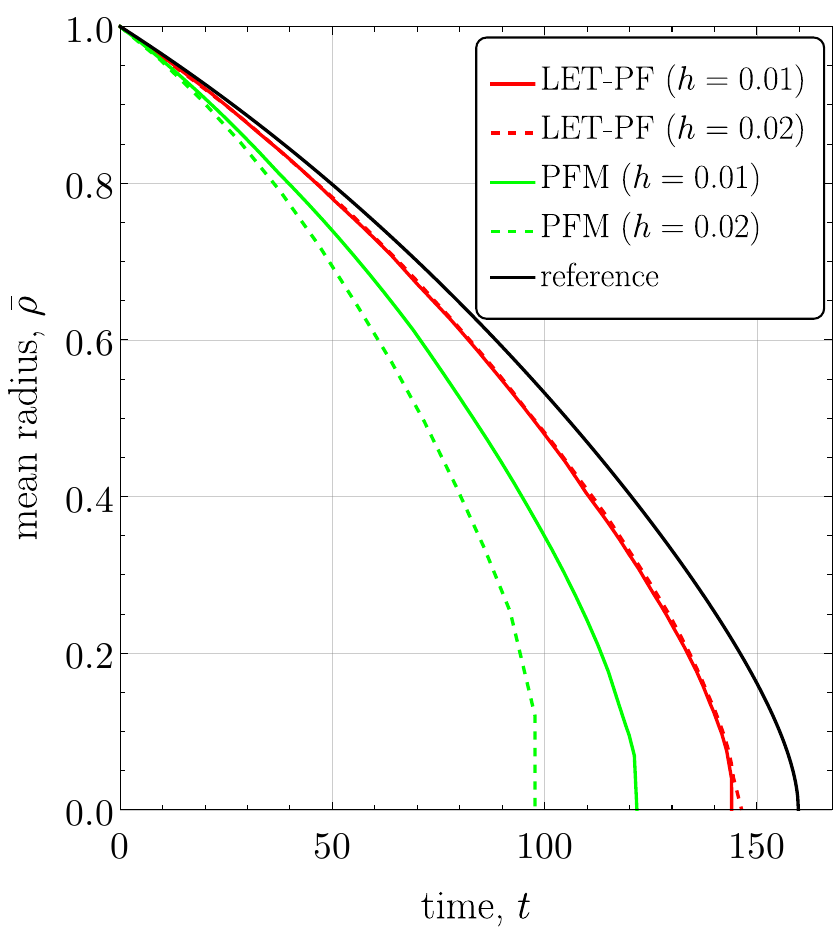} \\[1ex]
                        \hspace*{2em}(a) & \hspace*{2em}(b) & \hspace*{2em}(c) & \hspace*{2em}(d)
                    \end{tabular}
                    }
                \caption{Time evolution of the mean inclusion radius $\bar{\rho}$ for $\gamma=0.0008$: (a) $\ell=0.75h$, (b) $\ell=h$, (c) $\ell=1.5h$, (d) $\ell=2h$. The results shown in the top row are obtained for a small time increment (adaptive time stepping with $\Delta t_{\rm max}=T_{\rm exact}/500$), and those in the bottom row for significantly larger time increments ($\Delta t_{\rm max}=T_{\rm exact}$).
                }
                \label{fig:RadiusEvolutions}
            \end{figure}

            The effect of mesh density and time increment is illustrated in Fig.~\ref{fig:RadiusEvolutions} for an intermediate value of the interfacial energy, $\gamma=0.0008$. In addition to the small time increment ($\Delta t_{\rm max}=T_{\rm exact}/500$) and thus highly accurate time integration, see the top figures in Fig.~\ref{fig:RadiusEvolutions}, the analysis has also been performed with no restriction on the time increment ($\Delta t_{\rm max}=T_{\rm exact}$) so that the actual time increments result from the current convergence behaviour and are significantly larger, see the bottom figures in Fig.~\ref{fig:RadiusEvolutions}. It can be seen that the effect of time step is significant. As discussed later, the number of time steps needed to complete the simulation can be (and will be) used as a measure of the robustness of the model.
            Concerning the effect of mesh density, it is more pronounced in the case of PFM than LET-PF. Again, this effect will be studied in more detail later.
            %
            %
            %
            %
    
\begin{figure}[htbp]
                \centerline{\scriptsize
                    \begin{tabular}{cccc}
                        \rotatebox{90}{\small {\hspace{4eM} \fbox{PFM}\hspace{8.5eM} \fbox{LET-PF}}} &
                        \includegraphics[width=0.28\textwidth]{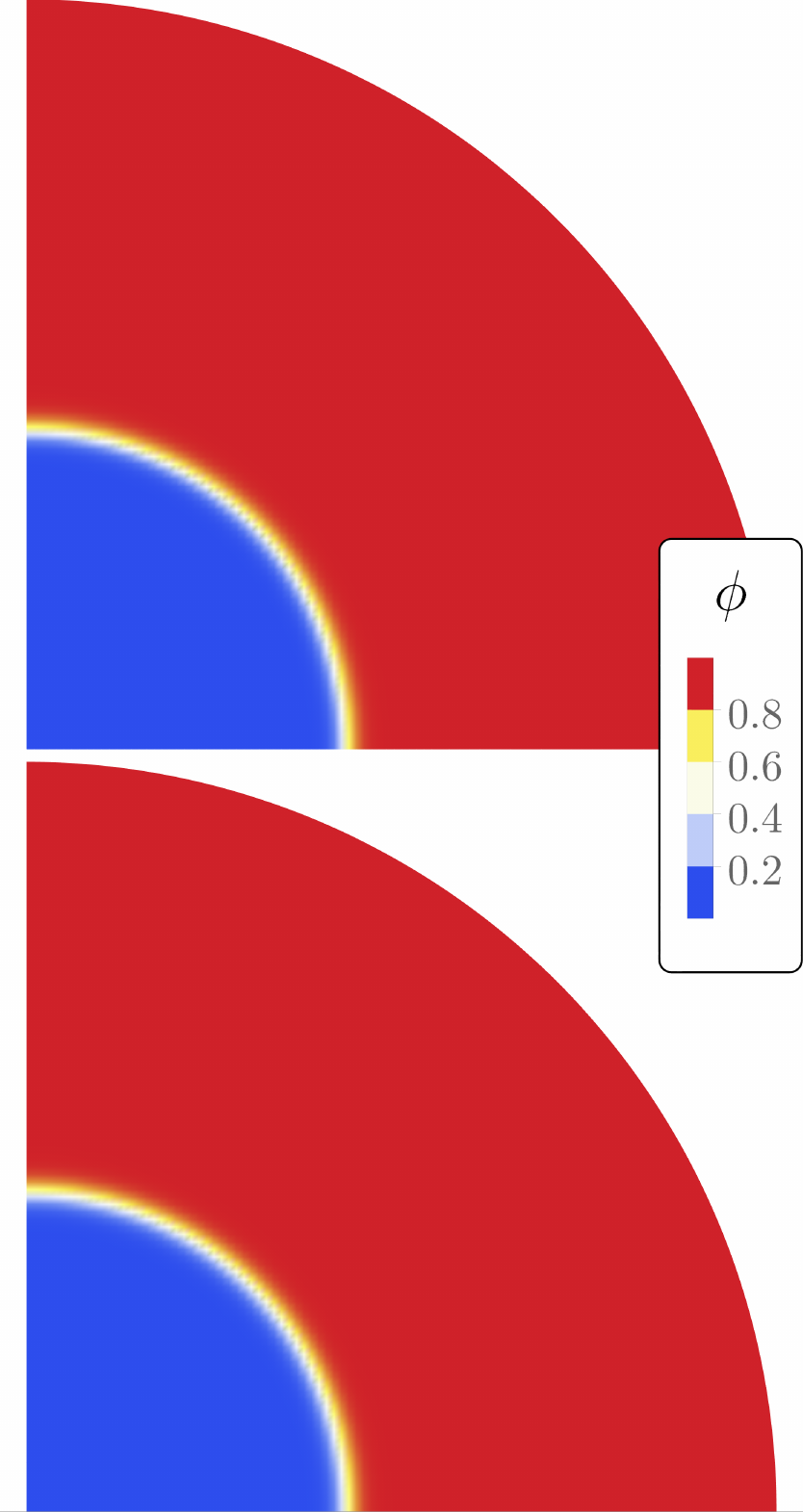} &
                        \includegraphics[width=0.28\textwidth]{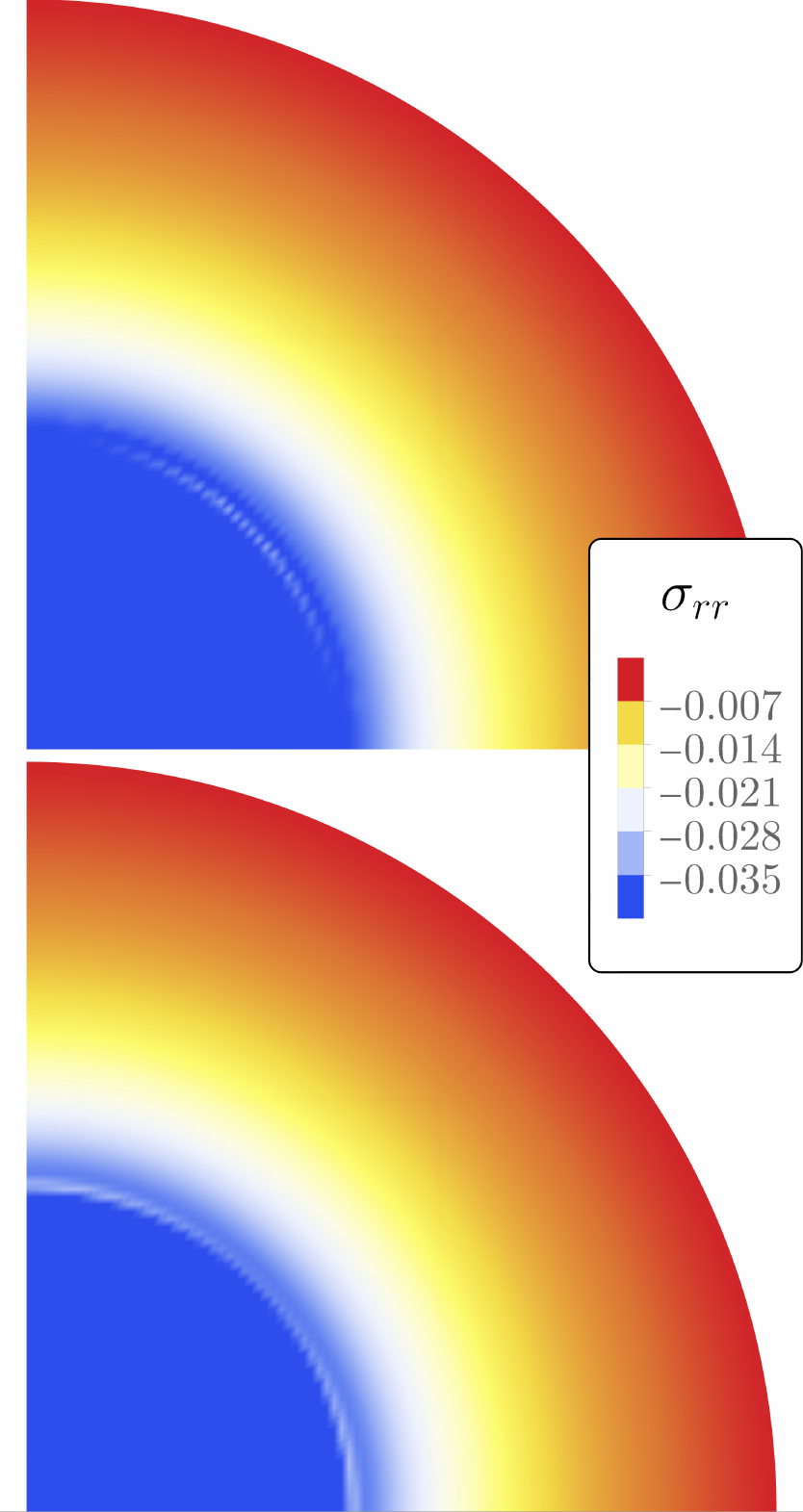} &
                        \includegraphics[width=0.28\textwidth]{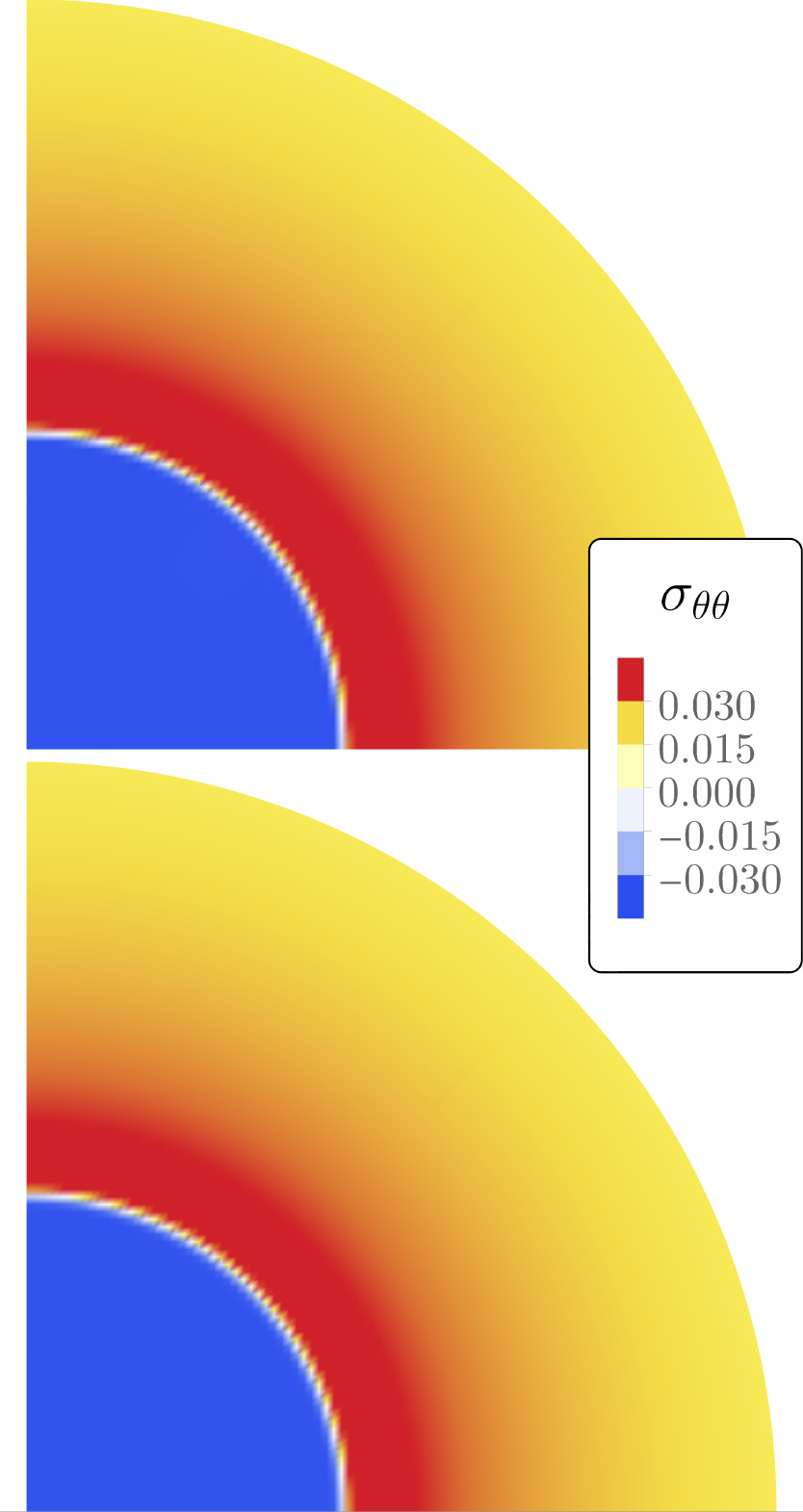} \\[2ex]
                         & (a) & (b) & (c)
                    \end{tabular}
                    }
                \caption{Maps of the order parameter $\phi$ (a), radial stress $\sigma_{rr}$ (b), and hoop stress $\sigma_{\theta\theta}$ (c) obtained using LET-PF (top) and PFM (bottom) at the instant when the inclusion radius reaches $\bar{\rho}\approx 0.85$ ($\gamma=0.0008$, $h=0.02$, $\ell=1.5h$).
                }
                \label{fig:Fields}
            \end{figure}

\begin{figure}[htbp]
                \centerline{\scriptsize
                    \begin{tabular}{cc}
                        \includegraphics[width=0.49\textwidth]{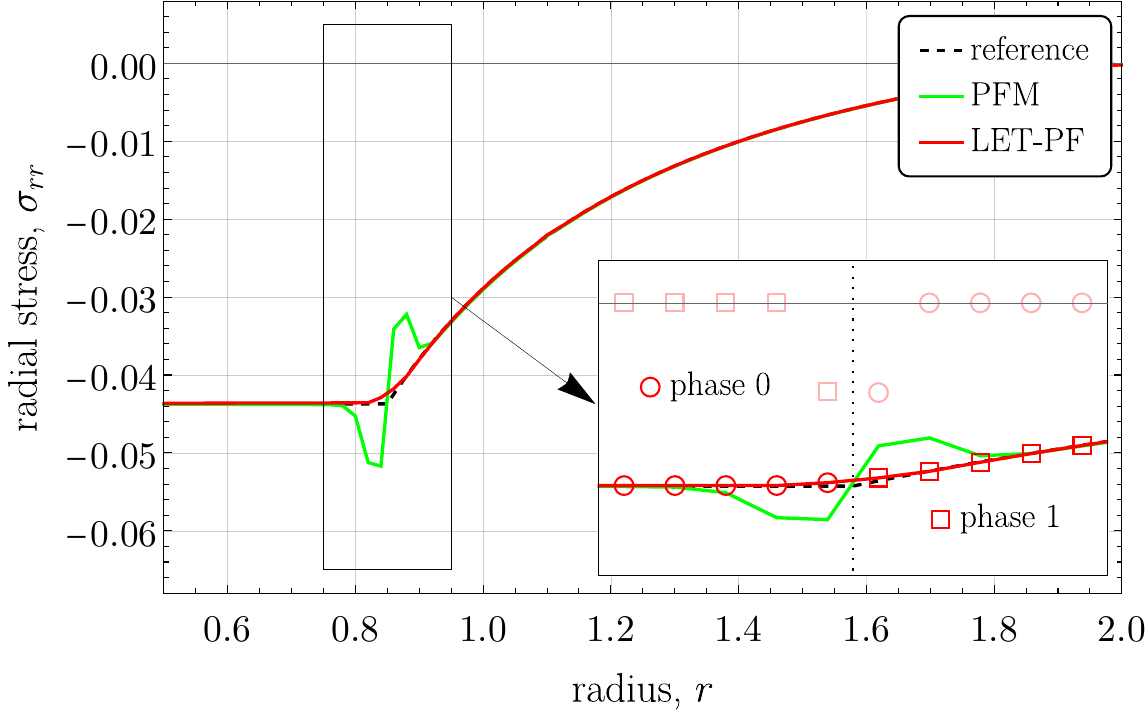} &
                        \includegraphics[width=0.49\textwidth]{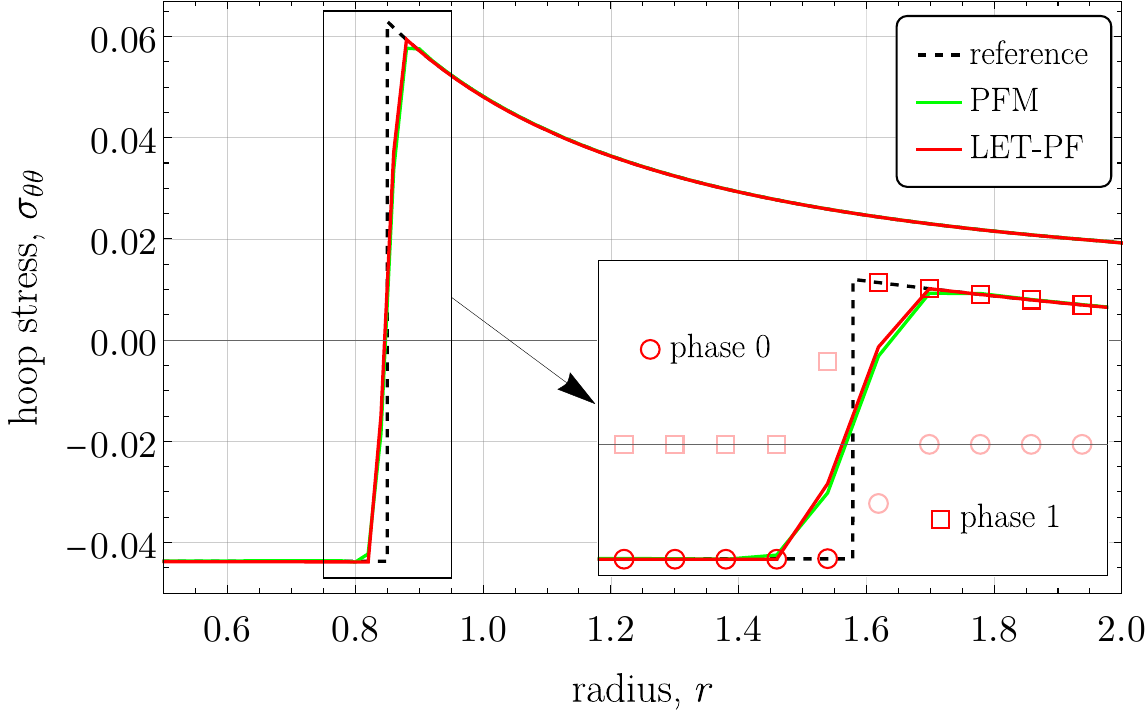}  \\[1ex]
                        \hspace*{3em}(a) & \hspace*{3em}(b) 
                    \end{tabular}
                    }
                \caption{Profiles of the radial stress $\sigma_{rr}$ (a) and hoop stress $\sigma_{\theta\theta}$ (b) along the bottom edge of the computational domain ($\bar{\rho}\approx 0.85$, $\gamma=0.0008$, $h=0.02$, $\ell=1.5h$). 
                The lines depict the overall stresses, while the markers in the insets represent the local stresses within each phase, as predicted by LET-PF.
                }
                \label{fig:StressPlots}
            \end{figure}

            Fig.~\ref{fig:Fields} shows the maps of the order parameter $\phi$, radial stress $\sigma_{rr}$ and hoop stress $\sigma_{\theta\theta}$ at the instant when the inclusion radius reaches $\bar{\rho}\approx0.85$. 
            While the order parameter $\phi$ and the hoop stress $\sigma_{\theta\theta}$ exhibit no substantial differences between LET-PF and PFM, the radial stress $\sigma_{rr}$, which according to the analytical solution should be continuous at the interface, is better represented by LET-PF. 
            This is particularly evident along the horizontal and vertical edges of the computational domain, as the interface 
            is then approximately parallel to the element edges and the LET technique is then highly accurate \citep{Dobrzanski2024}. 
            A more profound insight into these fields behaviours is presented in Fig.~\ref{fig:StressPlots}, showing the profiles of the radial and hoop stresses along the horizontal direction. 
            Regarding the radial stress, as shown in Figure~\ref{fig:StressPlots}(a), a significant distinction between the two methods is readily apparent. PFM exhibits fluctuations within the interface region, while LET-PF produces results much closer to the analytical solution. As for the hoop stress, as illustrated in Figure~\ref{fig:StressPlots}(b), the situation appears somewhat different---both methods seem to yield a similar response. However, it is crucial to note that, when utilizing LET-PF, the information is available also about the local strains and stresses resulting from the homogenization process in the laminated elements. The local stresses are indicated in the insets in Fig.~\ref{fig:StressPlots} by markers, thereby revealing that LET-PF provides also the local response that aligns more closely with the analytical solution.
            %
            %
            %
            %
            %
            %
            %
            %
    
            The results presented in Figs.~\ref{fig:RadiusFullEvolution} and~\ref{fig:RadiusEvolutions} suggest that LET-PF outperforms PFM in terms of accuracy. To examine the accuracy and efficiency of LET-PF in more detail, a comprehensive study has been performed by varying the governing material and numerical parameters. The material parameters are fully characterized by the dimensionless parameter $A$, Eq.~\eqref{eq:BimaterialLame:rhodotwA}${}_2$, and this parameter is controlled by varying the interfacial energy $\gamma$, as described above, while the remaining material parameters are fixed. Concerning the numerical parameters, two mesh densities are considered, $h\in\{0.01,0.02\}$, four values of the interface thickness parameter, $\ell/h\in\{0.75,1,1.5,2\}$, and two values of the maximum time increment, $\Delta t_{\rm max}\in\{T_{\rm exact},T_{\rm exact}/500\}$. 
    
            As shown in Fig.~\ref{fig:regimes}, when the inclusion radius approaches zero, the thermodynamic driving force, and specifically its interfacial contribution $\hat{f}_{\rm int}$, increases to infinity. The evolution is then fully governed by $\hat{f}_{\rm int}$ and the difference between LET-PF and PFM diminishes. 
            Moreover, the radius of curvature of the interface becomes then small compared to the element size and interface thickness so that additional effects influence the solution. 
            Accordingly, in the following, the final part of the evolution is excluded from the analysis, and only the range of the inclusion radii between $\rho_0$ and $0.15\rho_0$ is considered. The corresponding duration of the evolution process, resulting from the analytical solution, depends on $A$ and is denoted by $T_{\rm exact}$. The maximum time increment $\Delta t_{\rm max}$ in the adaptive time incrementation scheme is then expressed in terms of $T_{\rm exact}$, as specified above.

\begin{figure}[htbp]
                \centerline{\scriptsize
                    \begin{tabular}{cc}
                        \includegraphics[width=0.45\textwidth]{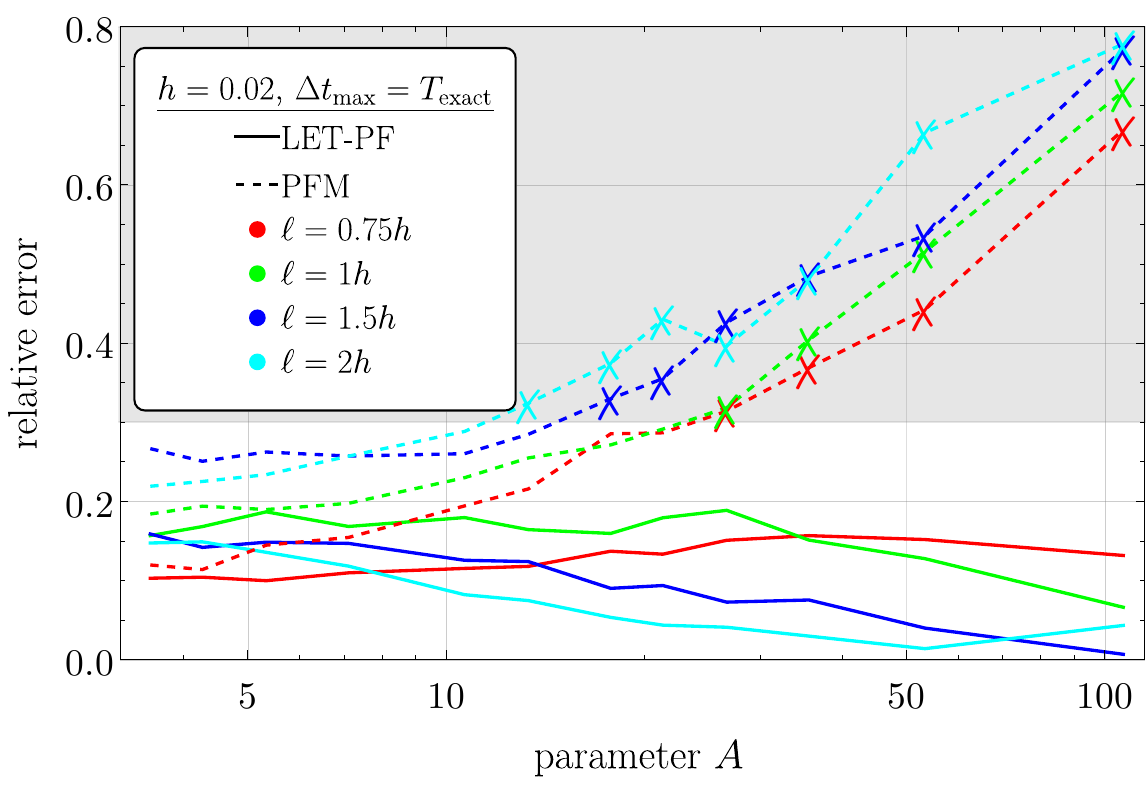} &
                        \includegraphics[width=0.45\textwidth]{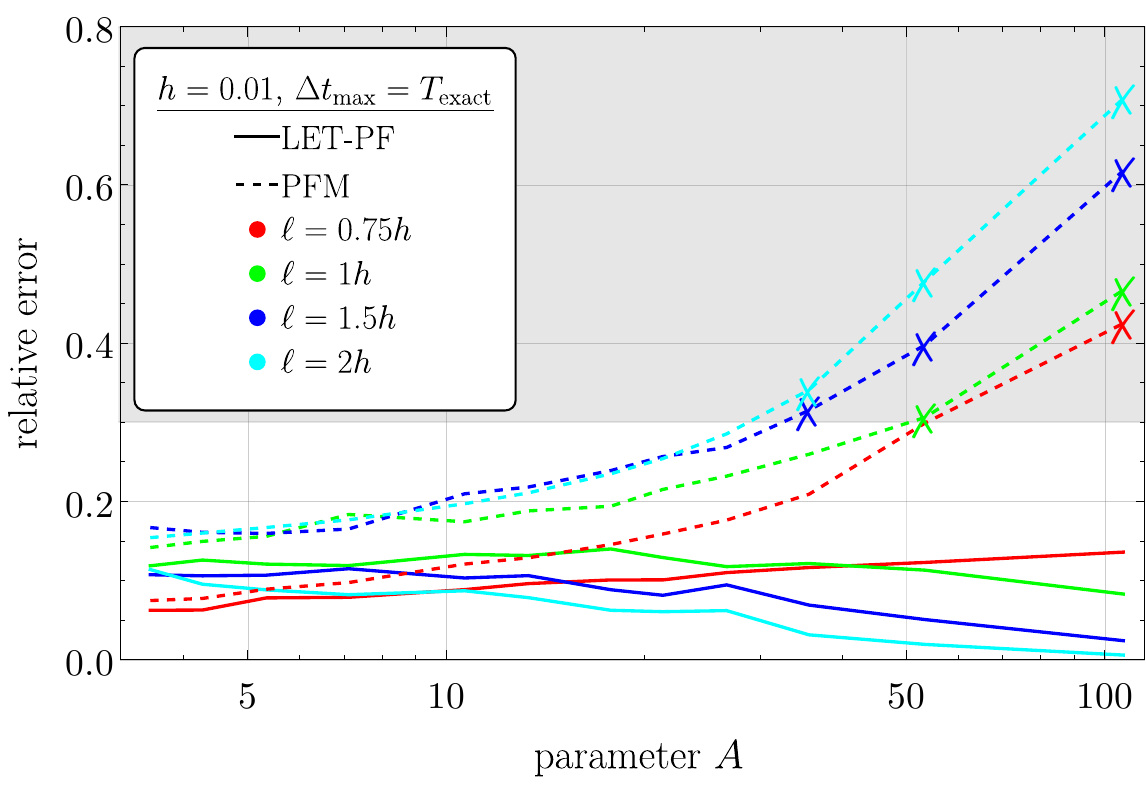} \\[1ex]
                        \hspace*{3em}(a) & \hspace*{3em}(b) \\[2ex]
                        \includegraphics[width=0.45\textwidth]{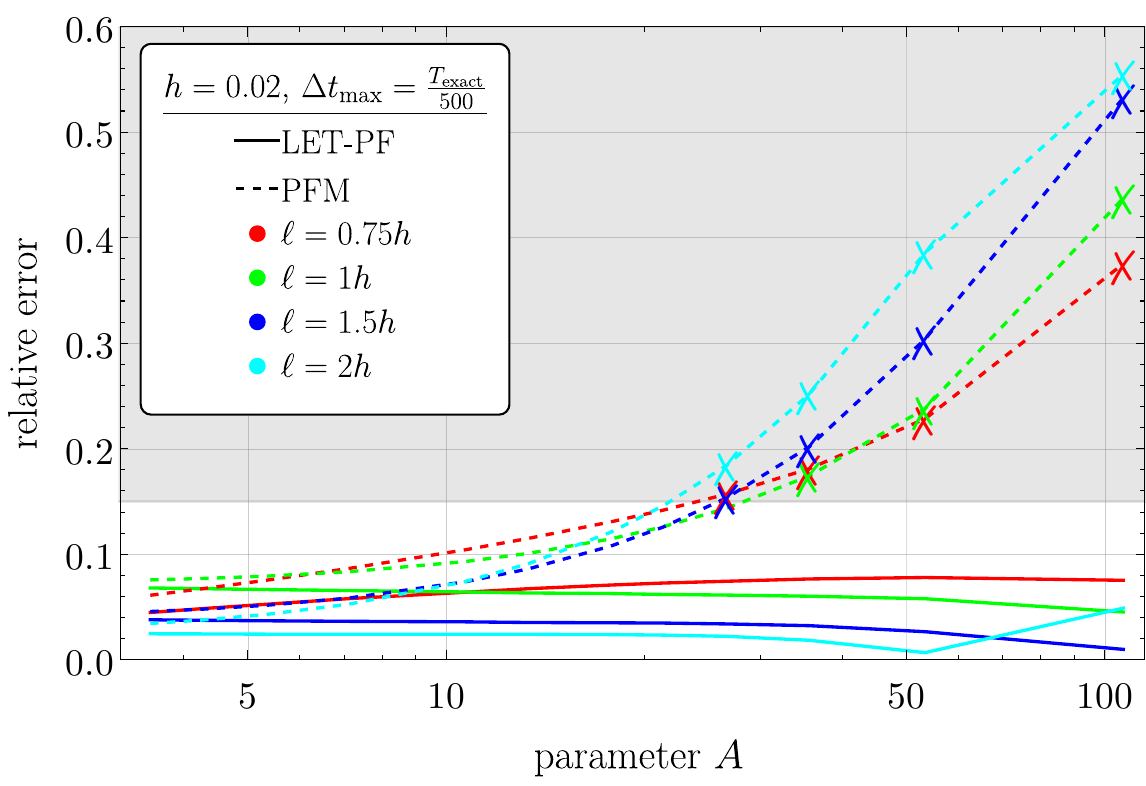} &
                        \includegraphics[width=0.45\textwidth]{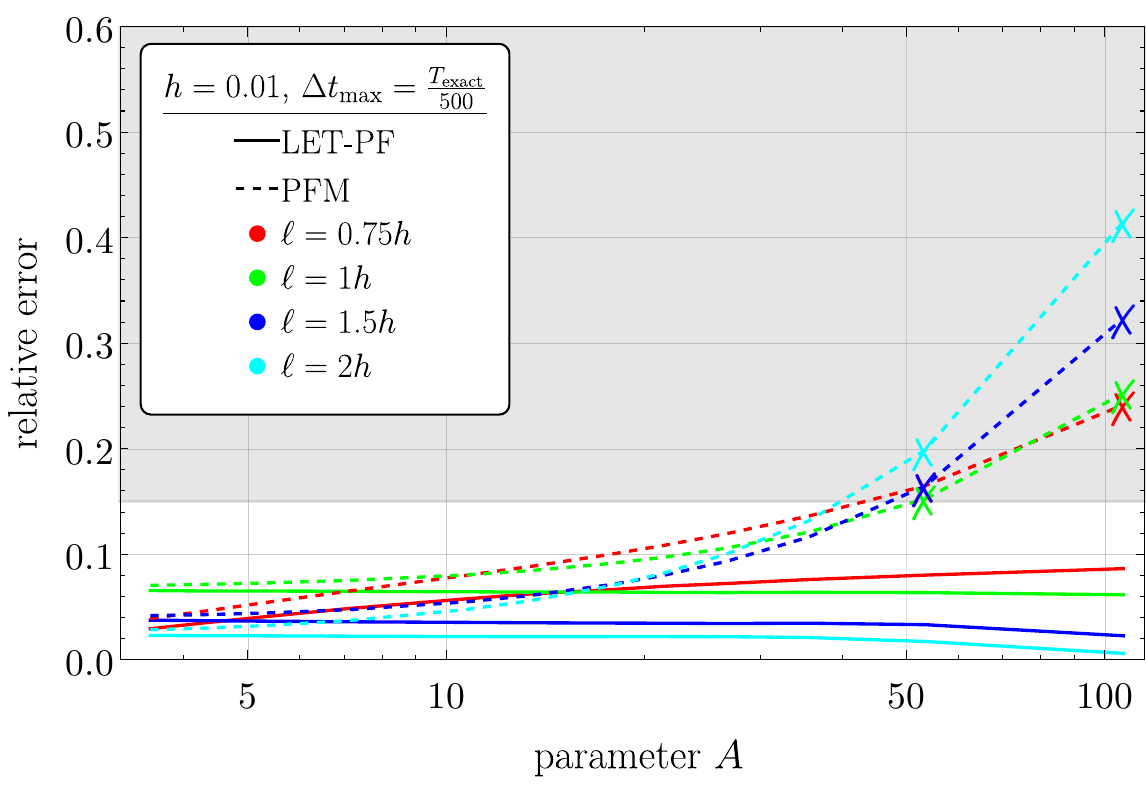} \\[1ex]
                        \hspace*{3em}(c) & \hspace*{3em}(d)
                    \end{tabular}
                    }
                \caption{Relative error, Eq.~\eqref{eq:relerror}, evaluated for all considered cases for LET-PF (solid lines) and PFM (dashed lines).
                The shaded regions correspond to the error that is approximately twice higher than the average error of LET-PF, i.e., exceeding 0.3 for $\Delta t_{\max}=T_{\rm exact}$ (panels (a) and (b)) and exceeding 0.15 for $\Delta t_{\max}=T_{\rm exact}/500$ (panels (c) and (d)). The data points falling within the shaded regions are marked with crosses, see also Fig.~\ref{fig:Cost}.
                }
                \label{fig:RadiiRelErr}
            \end{figure}

            Fig.~\ref{fig:RadiiRelErr} presents the relative error as a function of parameter $A$. The relative error is defined here as
            \begin{equation}\label{eq:relerror}
                {\rm relative\,error}=\dfrac{\int_{0.15\rho_0}^{\rho_0} \lvert \tau^{\rm exact}(r)-\tau^{\rm num}(r) \rvert\,\rd r}{\int_{0.15\rho_0}^{\rho_0} \tau^{\rm exact}(r)\,\rd r},
            \end{equation}
            where $\tau={\rho}^{-1}$ is the inverse function of $\rho(t)$ that describes the time as a function of the inclusion radius, and $\tau^{\rm exact}$ and $\tau^{\rm num}$ correspond to the analytical and numerical solutions, respectively. 
            The use of the inverse function ${\rho}^{-1}$ is here necessary because the duration of the evolution process in each simulation is, in general, different than the duration in the analytical solution so that integration of $\lvert \rho^\text{exact}(t)-\rho^{\rm num}(t) \rvert$ would not be feasible. 
            The integration in Eq.~\eqref{eq:relerror} is performed numerically using the trapezoidal rule. 
            The results in Fig.~\ref{fig:RadiiRelErr} illustrate that, in terms of accuracy, LET-PF outperforms PFM in all cases studied. 
            Although both methods exhibit similar trends for lower values of $A$ (where the interfacial energy has a greater contribution), LET-PF demonstrates superior performance as the contribution of the bulk (elastic) energy becomes more prominent. 
            This is particularly evident in Fig.~\ref{fig:RadiiRelErr}(a,b), where the errors for PFM reach extreme values of nearly 0.8, while LET-PF performs much better. 
            In general, it can be seen that the accuracy of LET-PF only weakly depends on $A$. 
            %
            %
            %
            %
    
\begin{figure}[htbp]
                \centerline{\scriptsize
                    \begin{tabular}{cc}
                        \includegraphics[width=0.45\textwidth]{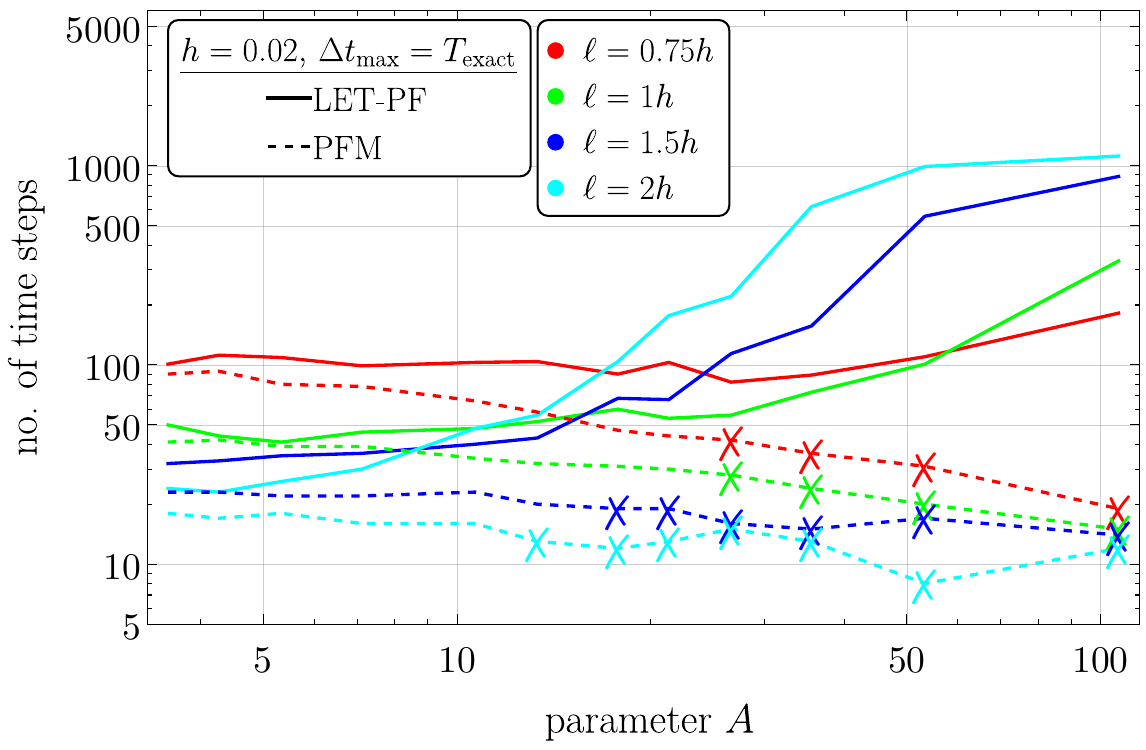} &
                        \includegraphics[width=0.45\textwidth]{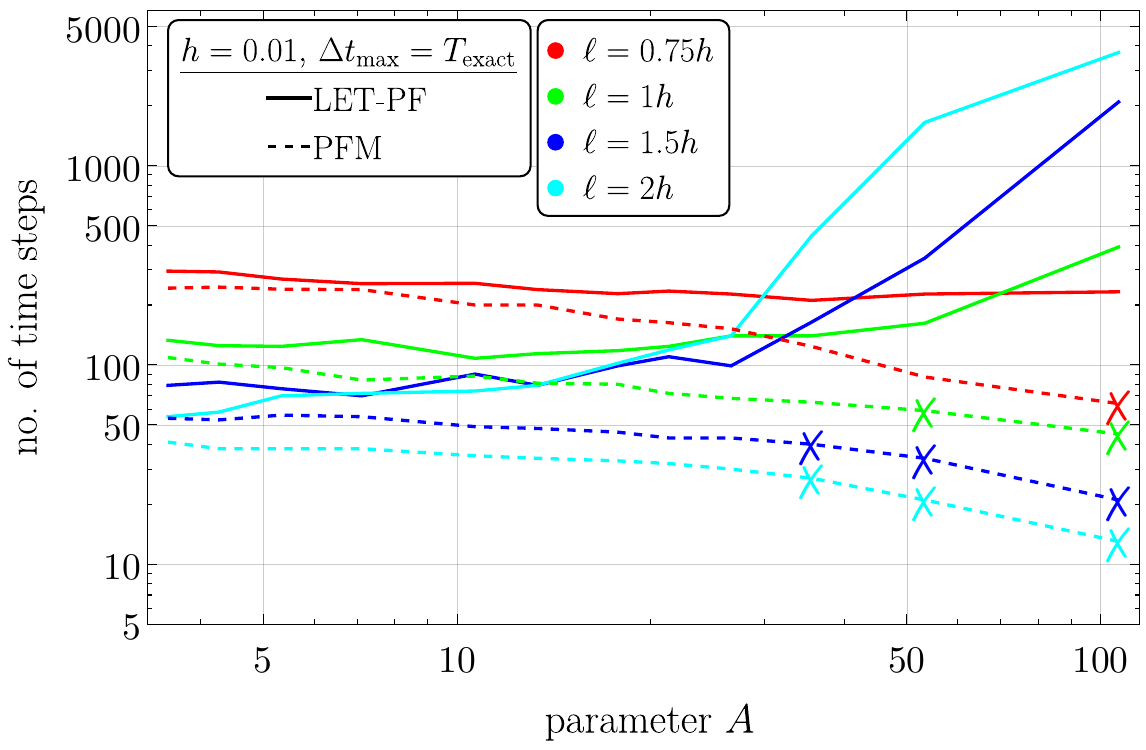} \\[1ex]
                        \hspace*{3em}(a) & \hspace*{3em}(b) \\[2ex]
                        \includegraphics[width=0.45\textwidth]{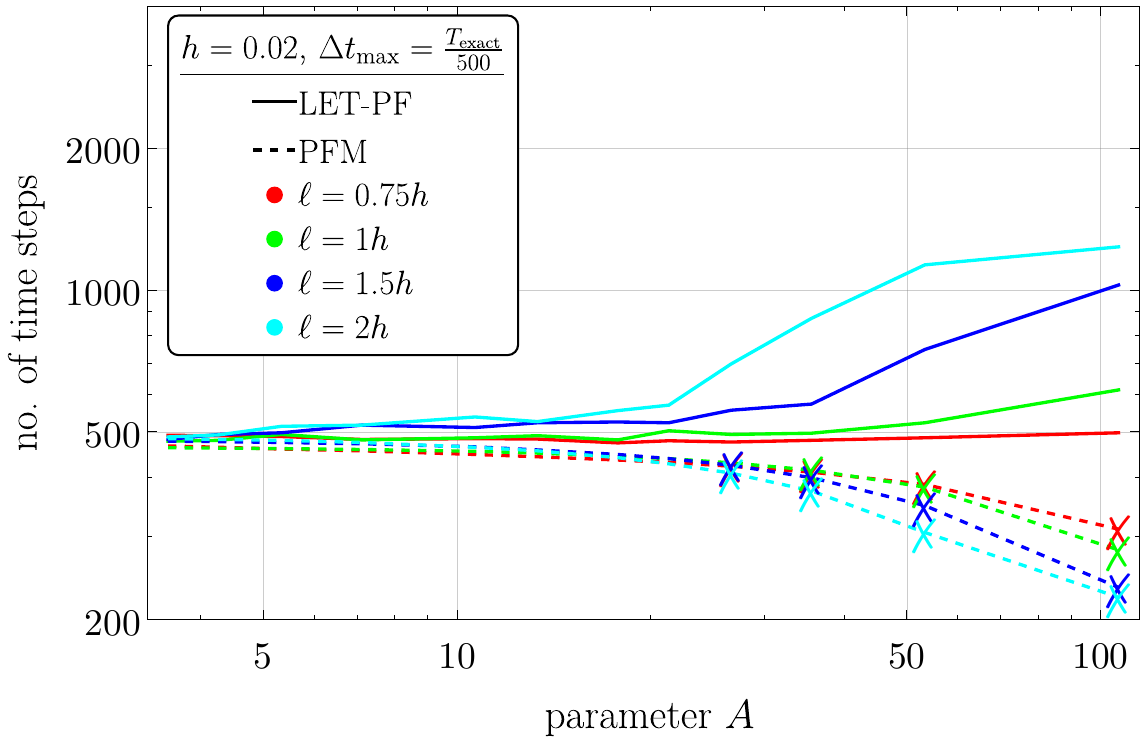} &
                        \includegraphics[width=0.45\textwidth]{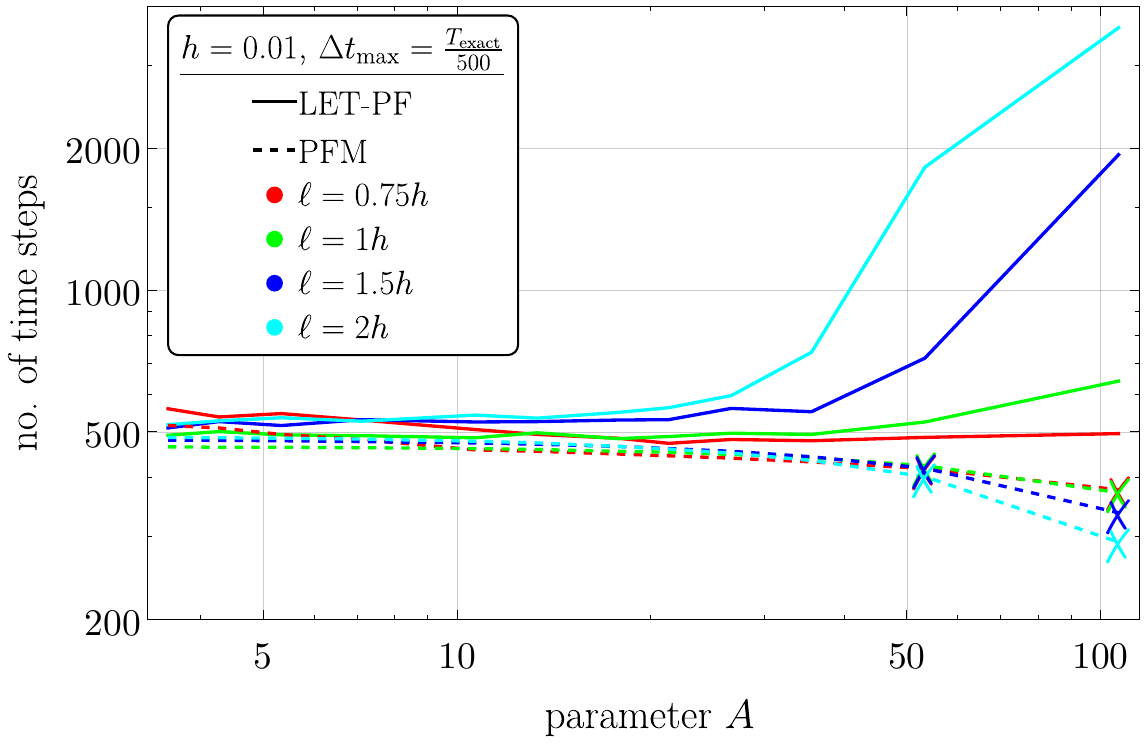} \\[1ex]
                        \hspace*{3em}(c) & \hspace*{3em}(d)
                    \end{tabular}
                    }
                \caption{The number of time steps needed to complete the simulation for $\Delta t_{\rm max} = T_{\rm exact}$ (a,b) and for $\Delta t_{\rm max} = T_{\rm exact}/500$ (c,d). Two mesh sizes are used, namely $h=0.02$ (a,c) and $h=0.01$ (b,d). 
                The results correspond to the cases shown in Fig.~\ref{fig:RadiiRelErr}, and the cross markers indicate the data points of exceedingly high error, see Fig.~\ref{fig:RadiiRelErr}.
                }
                \label{fig:Cost}
            \end{figure}

            Fig.~\ref{fig:Cost} shows the number of time steps needed to complete the simulation for all cases included in Fig.~\ref{fig:RadiiRelErr}. 
            When there is no restriction on the maximum time increment ($\Delta t_{\rm max}=T_{\rm exact}$), so that the simulation proceeds with the time increment dictated by the current convergence behaviour, the number of time steps may serve as an indicator of the robustness of the computational scheme (the smaller the number of time steps, the more robust the scheme). Moreover, the number of time steps quantifies the computational cost since the total computation time is proportional to the number of time steps. However, in terms of the cost, the comparison of LET-PF and PFM is not immediate, because the cost of evaluation of the element quantities (residual vector and tangent matrix) differs for the two methods. 
            The corresponding results are presented in Fig.~\ref{fig:Cost}(a,b). It can be seen that, for small $A$, the number of time steps is similar for both methods, PFM being slightly more efficient. For high values of $A$, the number of time steps significantly increases for LET-PF, while it slightly decreases for PFM. However, caution is needed when interpreting these results. As has been shown in Fig.~\ref{fig:RadiiRelErr}, PFM is highly inaccurate for high $A$, and the corresponding data points are marked with crosses in Fig.~\ref{fig:Cost}, consistent with Fig.~\ref{fig:RadiiRelErr}. It follows that PFM is indeed cheaper and more robust for high $A$, but it is then unacceptable in terms of accuracy. 
            On the other hand, for high $A$, LET-PF can only proceed with small time increments, but it delivers results of similar accuracy as for small $A$.
            %
            %
            %
            %
    
            When the time increment is small ($\Delta t_{\rm max} = T_{\rm exact}/500$), for small and intermediate $A$, the simulation proceeds with the prescribed maximum time increment, and the number of time steps is close to 500 for both LET-PF and PFM, see Fig.~\ref{fig:Cost}(c,d). For high $A$, in the case of LET-PF, the number of time steps increases as in the case of $\Delta t_{\rm max} = T_{\rm exact}$. On the other hand, in the case of PFM, the number of time steps decreases for high $A$, and this is related to the high error of the method and to the related reduction of the duration of the evolution process, see, for instance, Fig.~\ref{fig:RadiusFullEvolution}(b). Note that the corresponding data points are marked by crosses, as in Fig.~\ref{fig:RadiiRelErr}.
    
\begin{figure}[htbp]
                \centerline{\scriptsize
                    \begin{tabular}{cccc}
                        \includegraphics[width=0.25\textwidth]{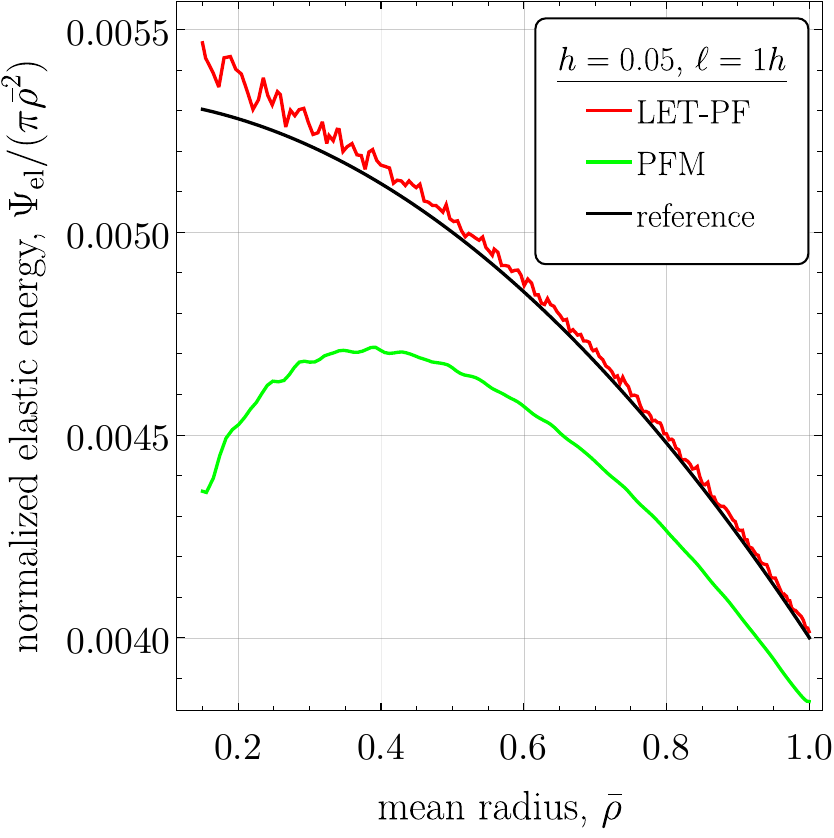} &
                        \includegraphics[width=0.25\textwidth]{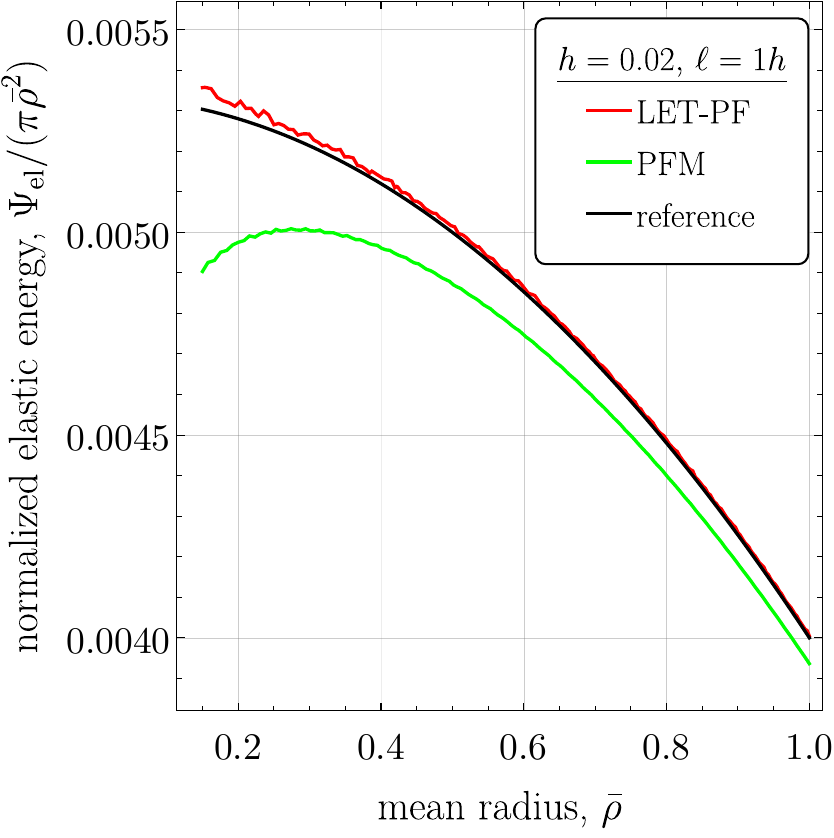} &
                        \includegraphics[width=0.25\textwidth]{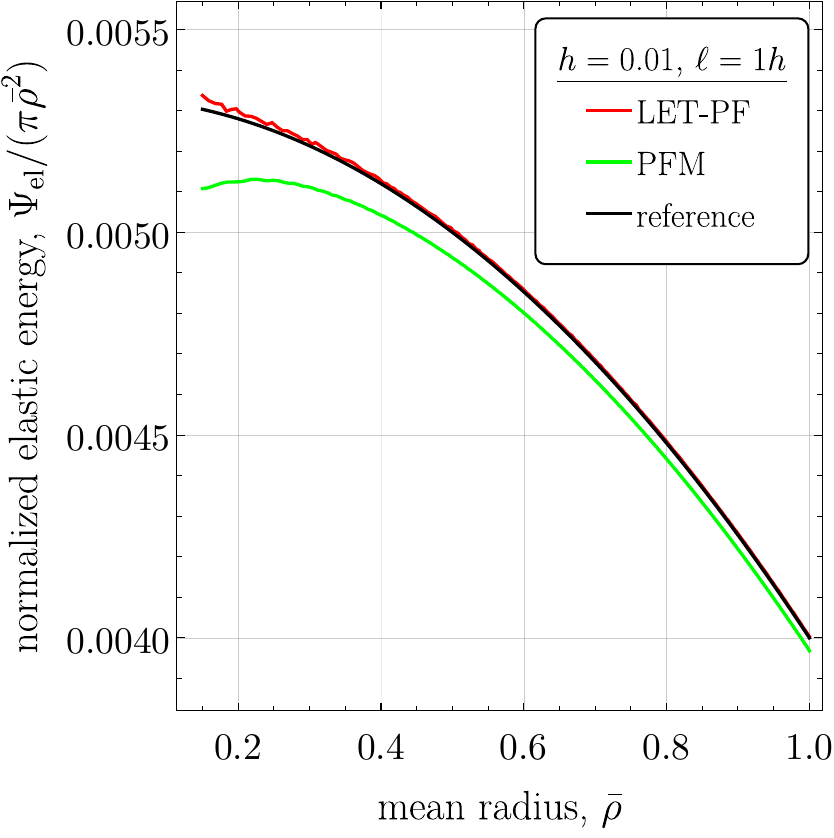} &
                        \includegraphics[width=0.25\textwidth]{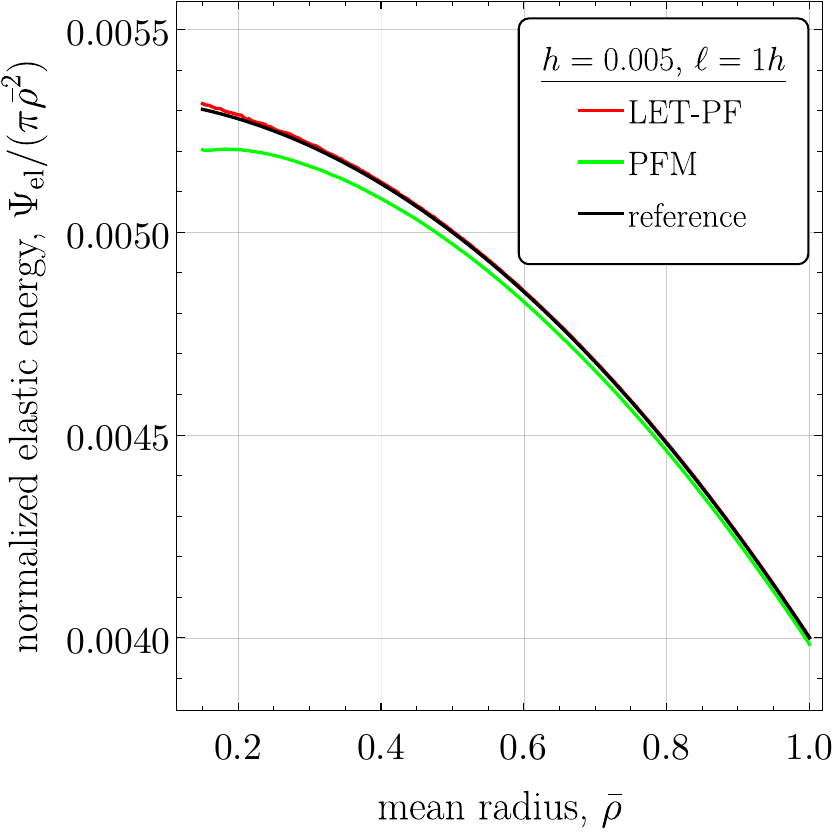} \\[1ex]
                        \multicolumn{4}{c}{(a)}\\[3ex]
                        \includegraphics[width=0.25\textwidth]{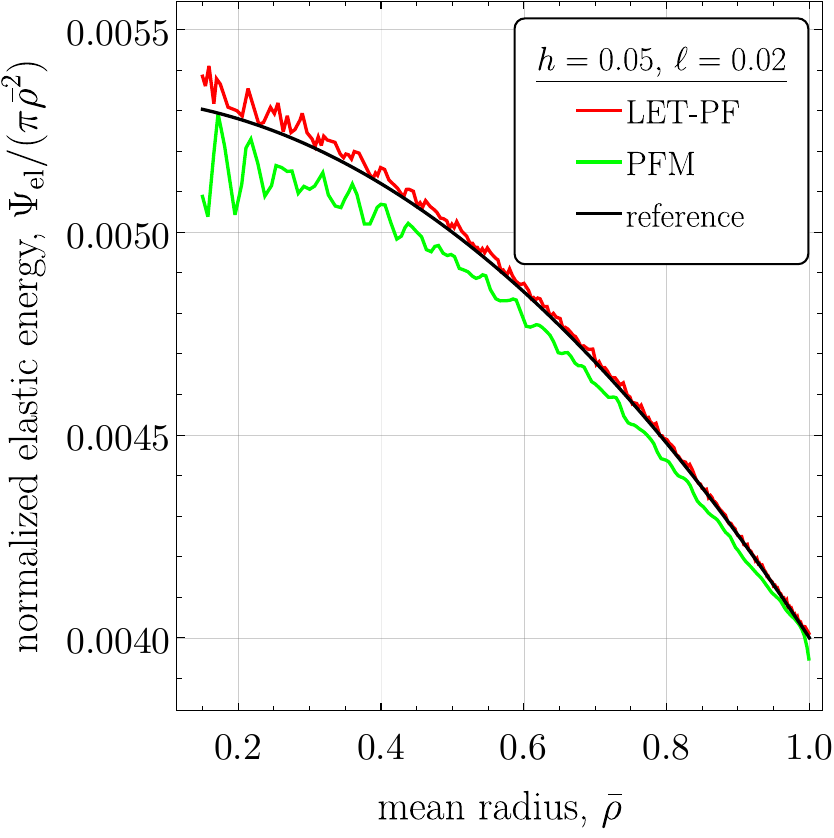} &
                        \includegraphics[width=0.25\textwidth]{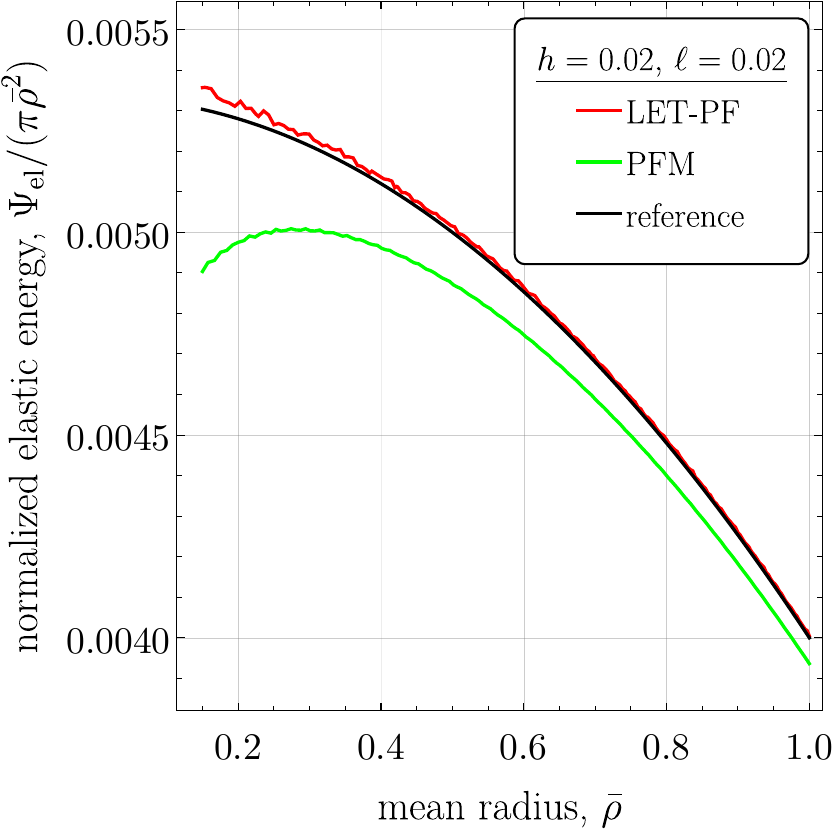} &
                        \includegraphics[width=0.25\textwidth]{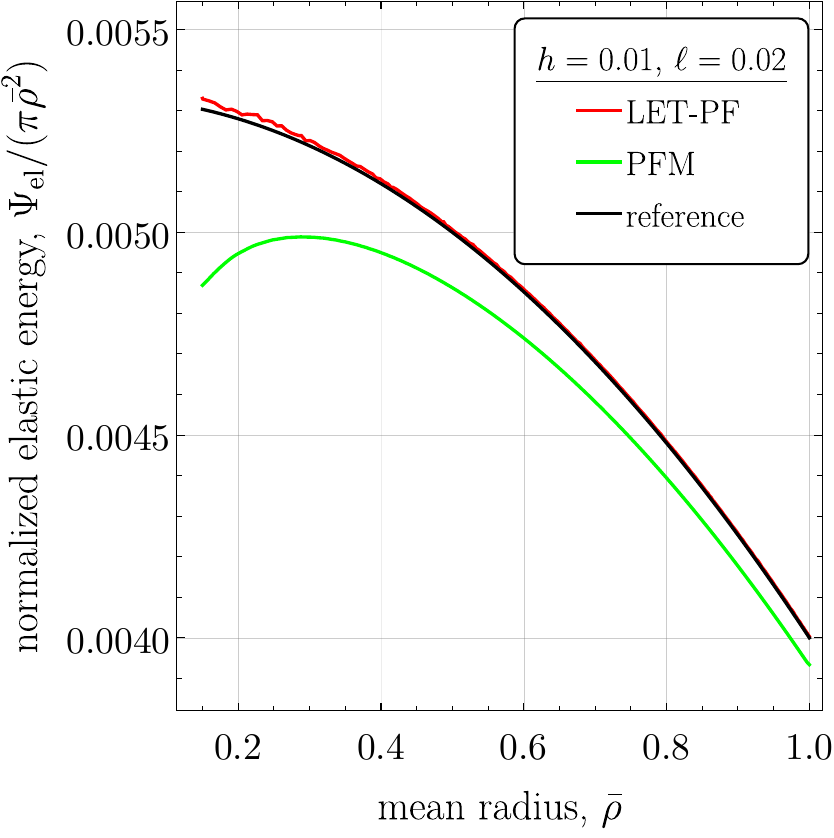} &
                        \includegraphics[width=0.25\textwidth]{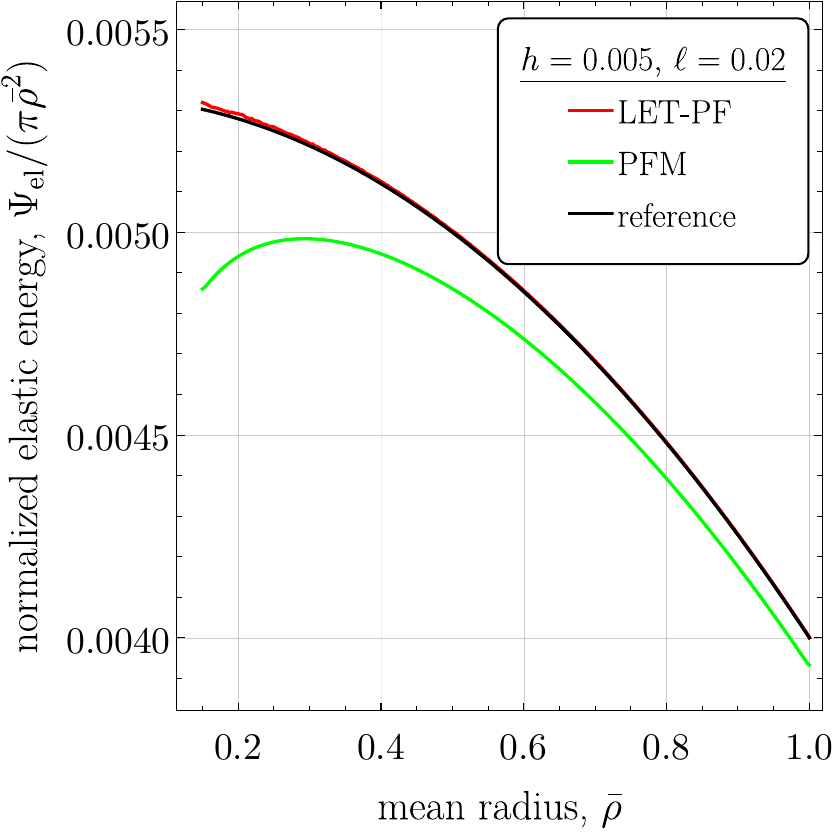} \\[1ex]
                        \multicolumn{4}{c}{(b)}
                    \end{tabular}
                    }
                \caption{Normalized total elastic strain energy $\Psi_{\rm el}/(\pi\bar{\rho}^2)$ as a function of the mean radius $\bar{\rho}$ for increasing mesh density (element size $h$ decreases from left to right) and for two scenarios: (a) $\ell=h$, (b) $\ell=0.02$ ($\gamma=0.0008$, $\Delta t_{\rm max}=T_{\rm exact}/200$).
                }
                \label{fig:ConvRateElEnergy}
            \end{figure}

    
            The results reported in Fig.~\ref{fig:RadiiRelErr} show that the relative error significantly depends on parameter $A$, on the maximum time increment $\Delta t_{\rm max}$, and on the interface thickness parameter $\ell$. The effect of mesh density, although present, is hard to perceive as the other dependencies prevail. Accordingly, additional computations have been performed in which four mesh densities have been used, $h\in\{0.05,0.02,0.01,0.005\}$, and two scenarios concerning the interface thickness parameter $\ell$ have been considered. In the first scenario, $\ell$ is proportional to the mesh size $h$ and is set to $\ell=h$. In the second scenario, $\ell$ is set to a constant value $\ell=0.02$ so that $\ell/h\in\{0.4,1,2,4$). 
            The simulations are performed for $\gamma=0.0008$ and $\Delta t_{\rm max}=T_{\rm exact}/200$. 
            Fig.~\ref{fig:ConvRateElEnergy} shows the total elastic strain energy $\Psi_{\rm el}$ as a function of the mean inclusion radius $\bar{\rho}$. 
            
            The energy is normalized by the area of the inclusion, $\Psi_{\rm el}/(\pi\bar{\rho}^2)$, so that the differences between the individual curves are better visible. 
            It can be seen that LET-PF outperforms PFM in both scenarios. 
            In the fist scenario ($\ell=h$), LET-PF converges to the analytical solution visibly faster than PFM, see Fig.~\ref{fig:ConvRateElEnergy}(a). This is because LET-PF treats the interface in a semisharp manner, and it is the diffuseness of the interface in PFM that is the source of the observed error. 
            This is illustrated in Fig.~\ref{fig:ConvRateElEnergy}(b), where the interface thickness is kept constant ($\ell=0.02$) and PFM converges to a solution that is significantly different from the correct one, particularly at small $\bar{\rho}$, when the interface thickness is large relative to the inclusion radius. 
            At the same time, the large thickness of the interface does not affect the performance of LET-PF because the $\phi=\frac12$ level set, which specifies the actual interface in the elasticity problem, only weakly depends on the thickness of the diffuse interface. 
            %
            %
            %
            %
            
\begin{figure}[htbp]
                \centerline{\scriptsize
                    \begin{tabular}{ccc}
                        \includegraphics[height=0.3\textwidth]{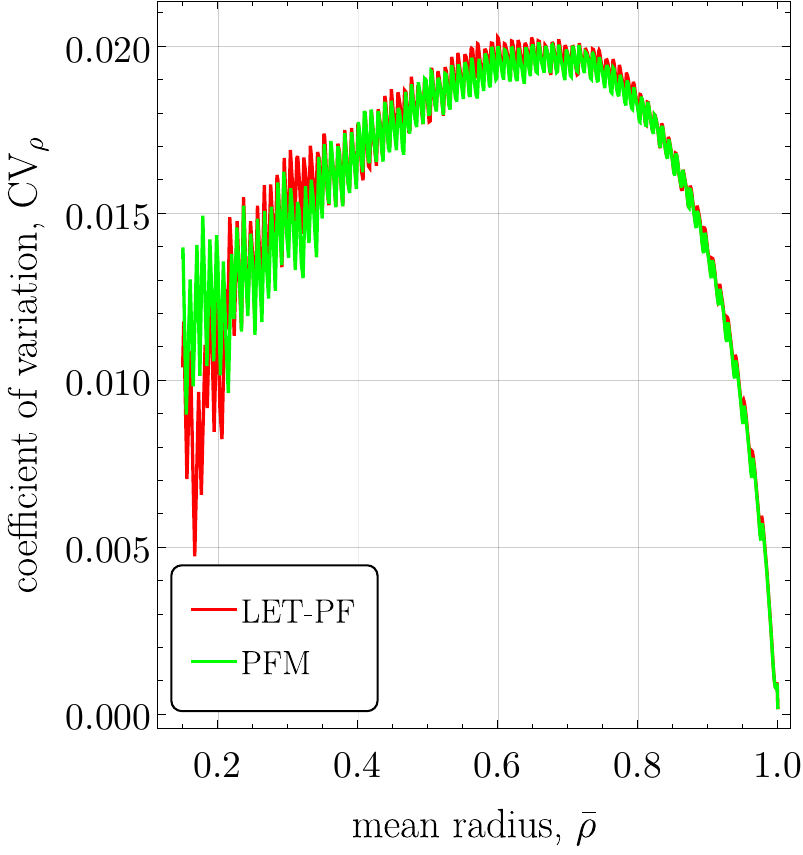} &
                        \includegraphics[height=0.3\textwidth]{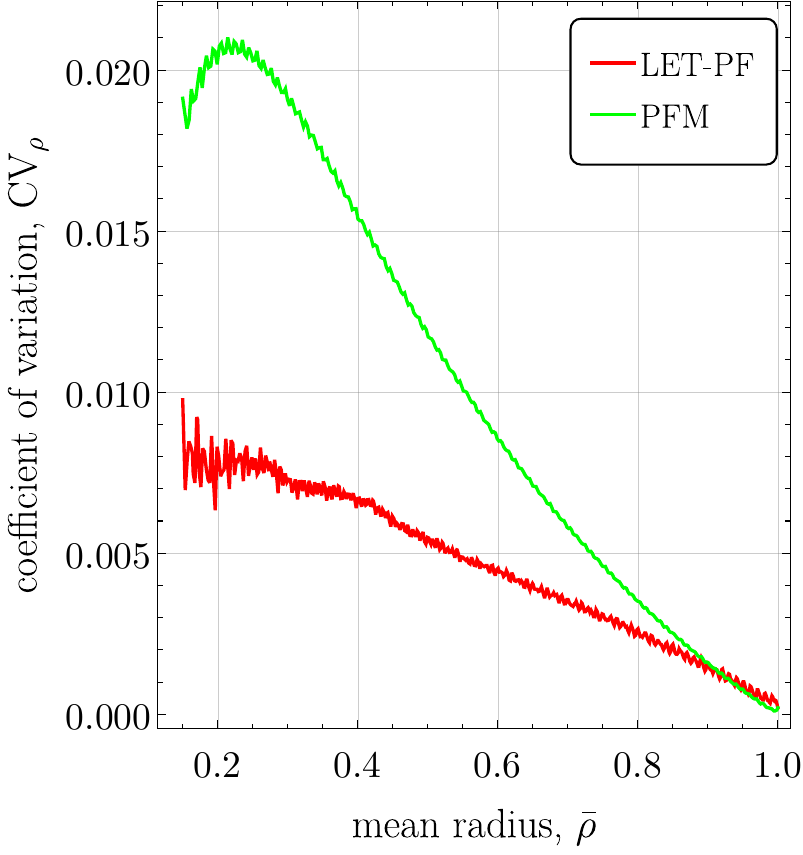} &
                        \includegraphics[height=0.3\textwidth]{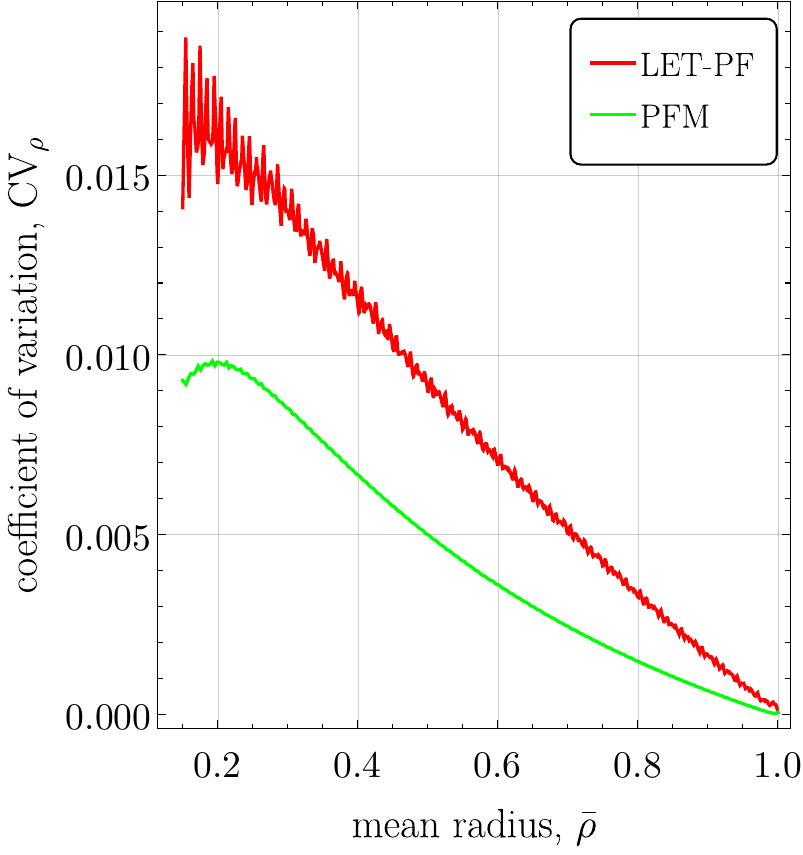} \\[1ex]
                       \hspace*{3em}(a) & \hspace*{3em}(b) & \hspace*{3em}(c)
                    \end{tabular}
                    }
                \caption{Deviation of the inclusion shape from a circular one, quantified by the coefficient of variation of the radius, $\mathrm{CV}_\rho$, Eq.~\eqref{eq:CVrho}, as a function of the mean inclusion radius $\bar{\rho}$ for three representative cases: (a) $\gamma=0.003$, $\ell=0.75h$; (b) $\gamma=0.0001$, $\ell=0.75h$; (c) $\gamma=0.0001$, $\ell=1.5h$ ($h=0.01$, $\Delta t_{\rm max}=T_{\rm exact}/500$).
                }
                \label{fig:RoundnessExamples}
            \end{figure}

            Since the reference sharp-interface problem is axisymmetric and the interface $\Gamma$ is then a circle, any deviation from the circular shape can be treated as a measure of the numerical solution error. 
            As a quantitative measure of this deviation, the coefficient of variation of the radius can be calculated at each time instant,
            \begin{equation}\label{eq:CVrho}
                \mathrm{CV}_\rho=\frac{\sigma_\rho}{\bar{\rho}},
            \end{equation}
            where $\sigma_\rho$ is the standard deviation of the orientation-dependent inclusion radius $\rho(\theta)$. 
            Some representative results presenting $\mathrm{CV}_\rho$ as a function of the mean inclusion radius $\bar{\rho}$ are shown in Fig.~\ref{fig:RoundnessExamples}. 
            It follows that the deviation from the circular shape depends on the model parameters and varies during the evolution process, including small oscillations resulting from the interface traversing the finite-element mesh. 
            Anyway, the deviation is relatively small in all cases (typically $\mathrm{CV}_\rho<0.02$).
            %
            %
            %
            %
    
\begin{figure}[htbp]
                \centerline{\scriptsize
                    \begin{tabular}{ccc}
                        \includegraphics[width=0.45\textwidth]{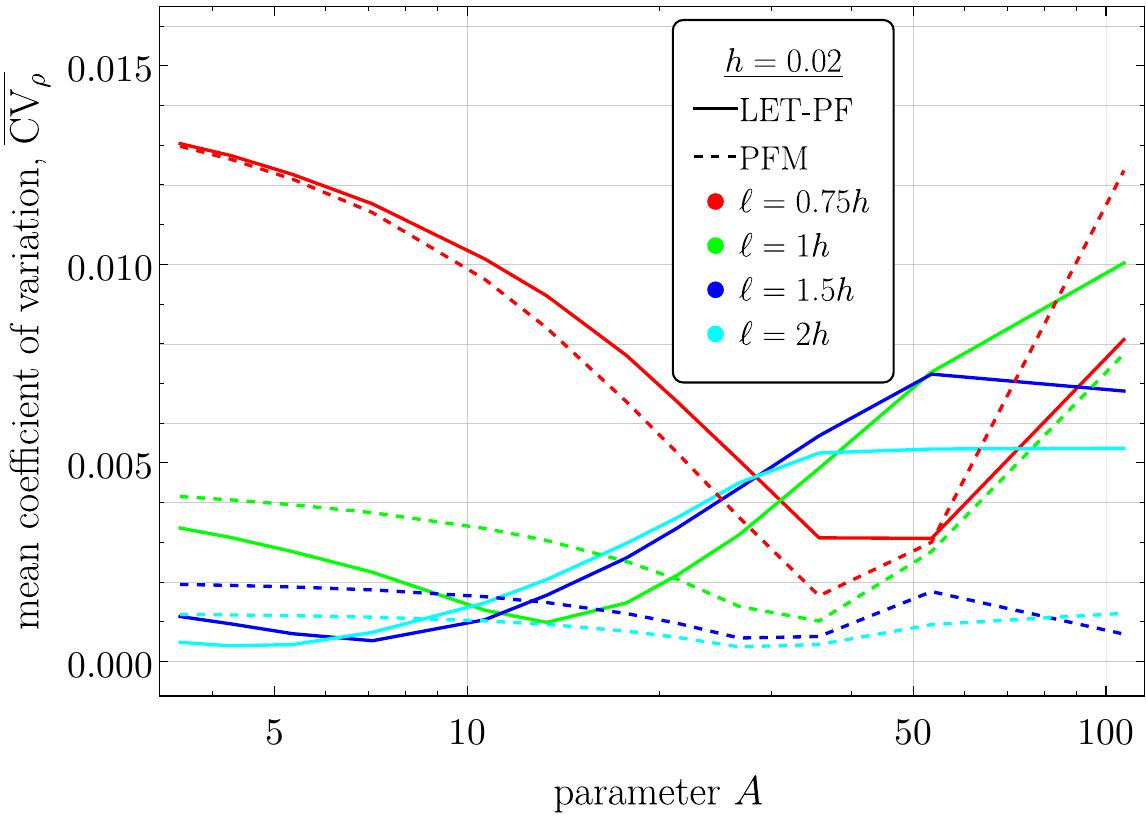} &&
                        \includegraphics[width=0.45\textwidth]{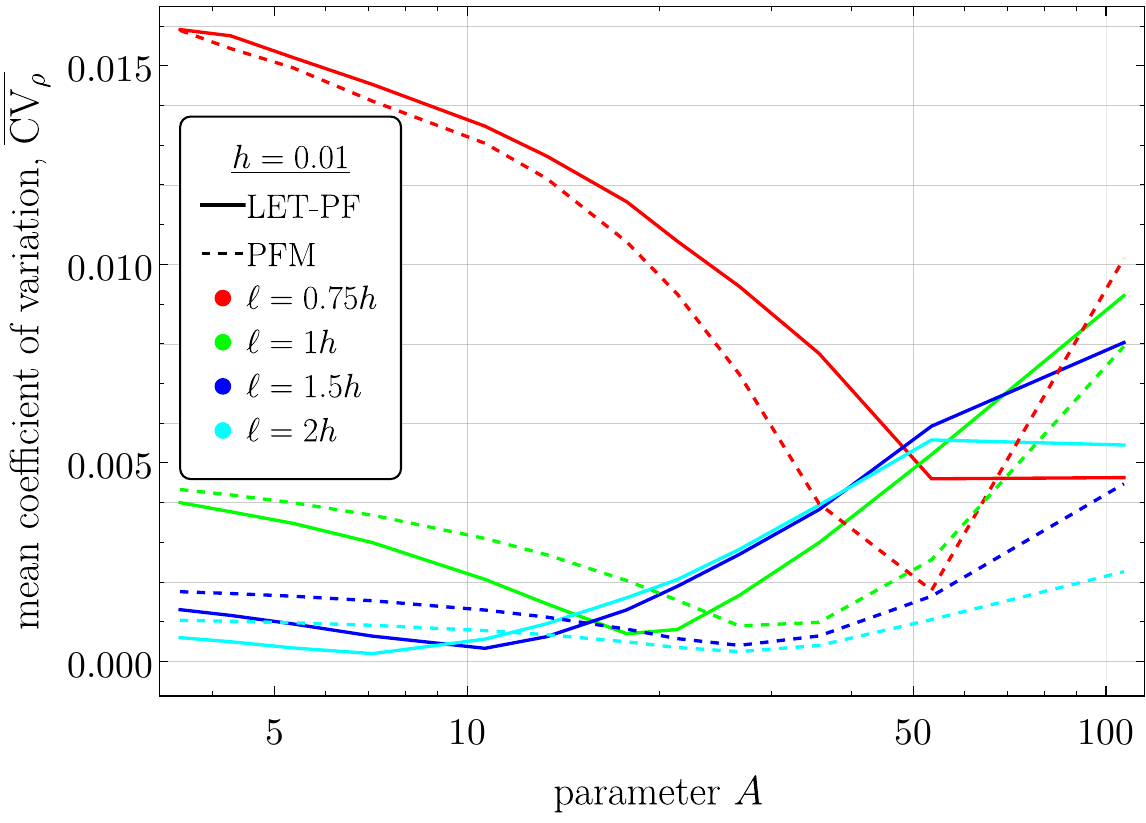} \\[1ex]
                        \hspace*{4em}(a) && \hspace*{4em}(b)
                    \end{tabular}
                    }
                \caption{Mean coefficient of variation of the inclusion radius $\overline{\mathrm{CV}_\rho}$, Eq.~\eqref{eq:meanCVrho}, for (a) $h=0.02$ and (b) $h=0.01$ ($\Delta t_{\rm max}=T_{\rm exact}/500$).
                }
                \label{fig:RoundnessMean}
            \end{figure}

            In order to comprehensively investigate this characteristic, the mean value of $\mathrm{CV}_\rho$ has been calculated for each case,
            \begin{equation}\label{eq:meanCVrho}
                \overline{\mathrm{CV}_\rho}=\frac{1}{0.85\rho_0}\Int_{0.15\rho_0}^{\rho_0}\mathrm{CV}_\rho \,\mathrm{d}\rho,
            \end{equation}
            and the aggregated results, separately for the two mesh sizes, are shown in Fig.~\ref{fig:RoundnessMean}. 
            This analysis reveals that both LET-PF and PFM result in a low mean value of the coefficient of variation of the radius with the maximum value less than 2\%. 
            However, it can be seen that, for the interface thickness $\ell=0.75h$, for both methods the deviation from the circular shape is higher than for larger interface thicknesses, which can be attributed to the interface pinning which is expected to diminish with the interface thickness increasing relative to the element size. 
            %
            %
            %
            %
    
            \FloatBarrier

    \subsection{Single inclusion in a constrained domain}\label{sec:single}
        In this example, a 2D square domain with a single evolving inclusion is considered. 
        To ensure a non-trivial steady-state solution, the displacements are fully constrained on the boundary and a volumetric eigenstrain of the same value but with opposite signs is introduced in both phases, namely $\bm{\varepsilon}^{\rm t}_{2}=-\bm{\varepsilon}^{\rm t}_{1}=\epsilon\bm{I}$ with $\epsilon=0.1$. 
        Accordingly, considering a viscous evolution towards the minimum of the elastic strain energy, a steady-state solution is expected with the overall volume fraction of both phases close to 0.5, because this leads to a null overall eigenstrain in the domain, thus effectively minimizing the total elastic strain energy (the chemical energy is not considered).

        The dimensions of the domain are $1\times 1$, the inclusion is initially located in the domain centre and its initial radius is $\rho_0=0.1$. 
        Homogeneous elastic properties are assumed, specifically, the Young's modulus $E=1$ and the Poisson's ratio $\nu=0.25$. Further, the following interface parameters are assumed: the interfacial energy density $\gamma=0.0003$, the effective mobility $\hat{m}=1$, and the interface thickness parameter $\ell=1.5h$.

\begin{figure}[htbp]
            \centerline{\scriptsize
                    \includegraphics[width=0.98\textwidth]{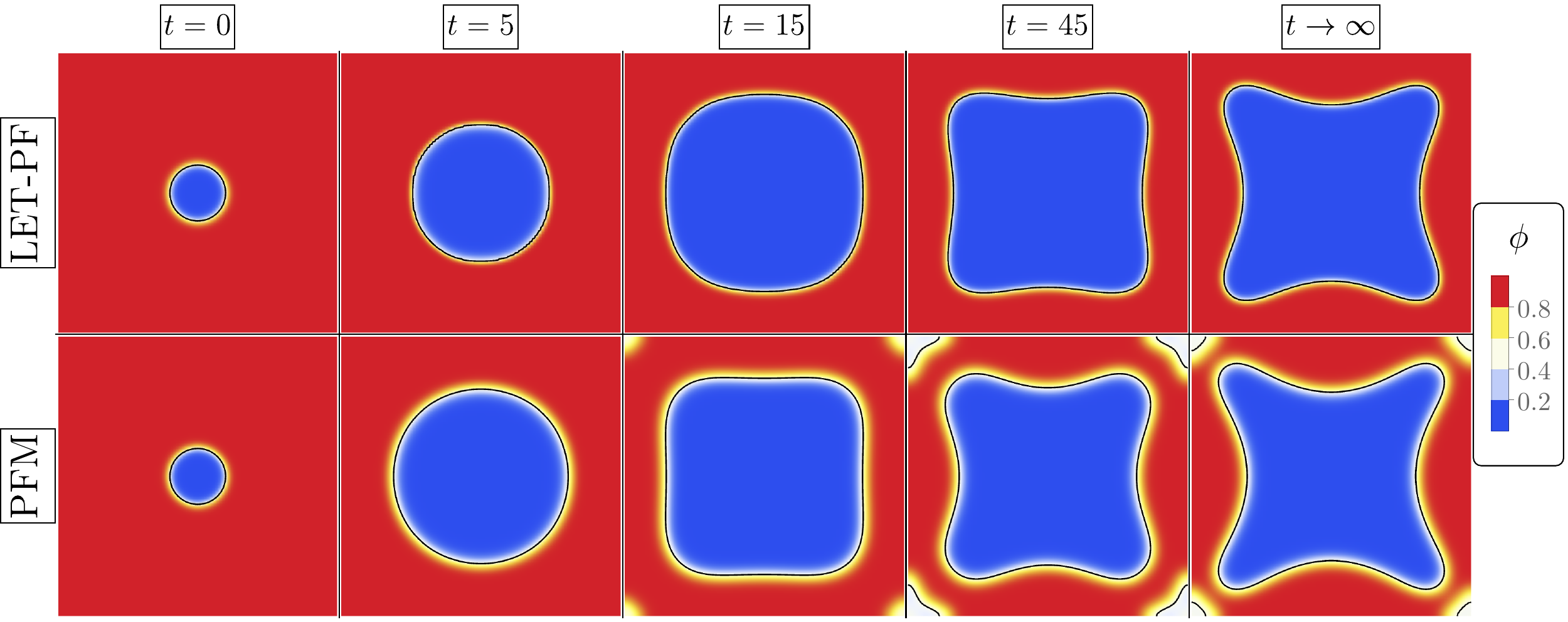} 
                }
            \caption{Evolution of the inclusion within a constrained domain for LET-PF (top) and PFM (bottom) obtained for the mesh of $100\times 100$ elements. The colour maps depict the order parameter $\phi$. 
            The black contours indicate the position of the $\phi=\frac12$ level set.
            }
            \label{fig:fullEvolution}
        \end{figure}

        Preliminary computations have been performed on the mesh of $100\times 100$ elements. Fig.~\ref{fig:fullEvolution} illustrates the evolution of the inclusion shape for both methods, starting from $t=0$ and ending in a steady state at $t\rightarrow\infty$ (a sufficiently large final time has been adopted).
        Discernible differences between LET-PF and PFM become evident almost instantaneously. 
        The final shape of the inclusion is different, with the corners being more acute for PFM. 
        Also, the evolution towards a square-like shape proceeds faster for PFM, compare the snapshots corresponding to $t=15$.
        Moreover, for PFM, artifacts manifest in the corners of the domain at a relatively early stage, a phenomenon not observed for LET-PF.
        Finally, in the case of PFM, the interface thickness increases with respect to the theoretical one, which is prescribed at $t=0$, while no visible difference in the interface thickness is observed for LET-PF. 
        %
        %
        %
        %

\begin{figure}[htbp]
            \centerline{\scriptsize
                    \includegraphics[width=0.98\textwidth]{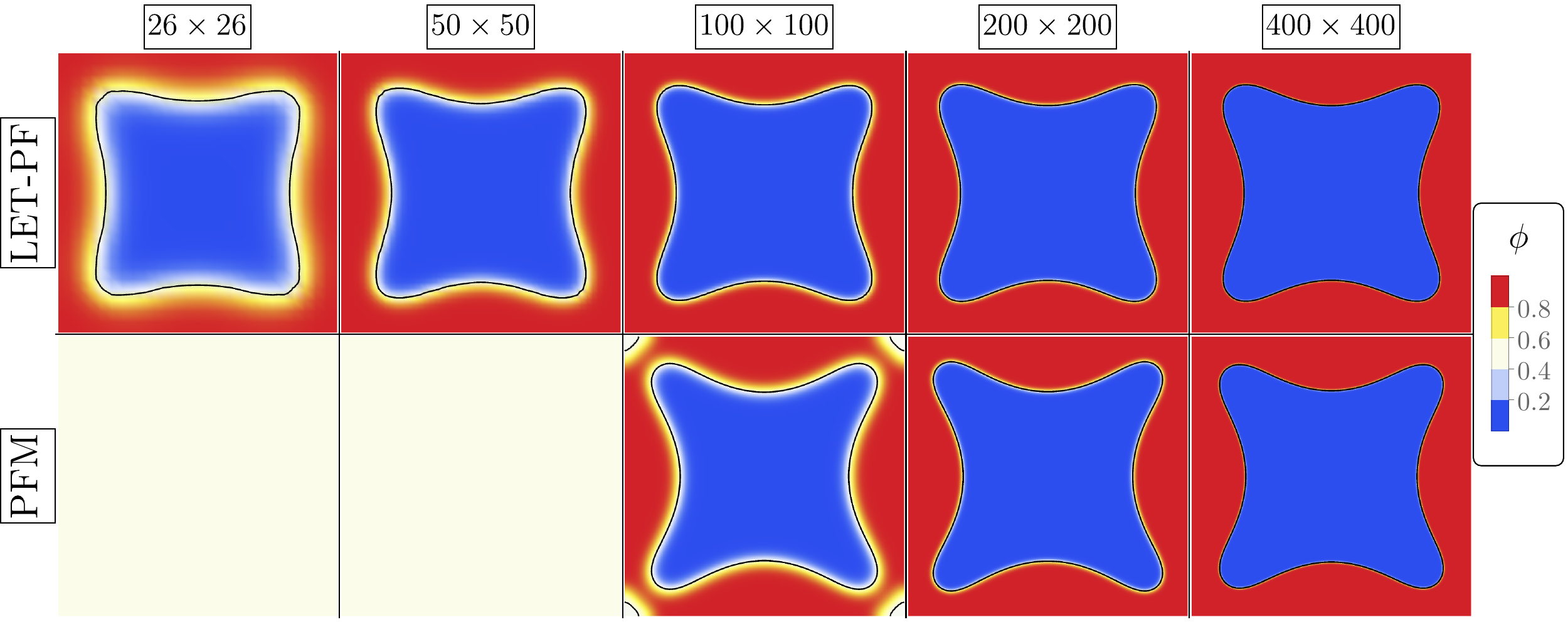} 
                }
            \caption{Steady-state shapes for LET-PF (top) and PFM (bottom) as a function of mesh density (increasing from left to right).}
            \label{fig:steadyStates}
        \end{figure}

\begin{figure}[htbp]
            \centerline{\scriptsize
                \begin{tabular}{cc}
                    \includegraphics[width=0.45\textwidth]{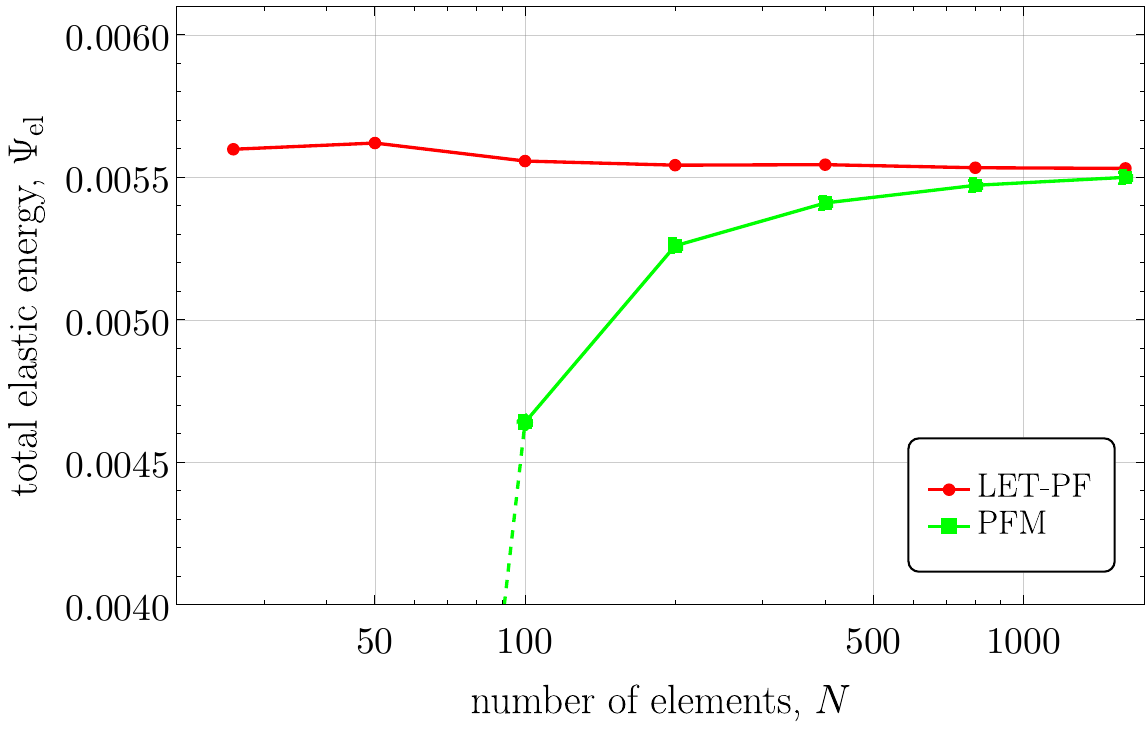} &
                    \includegraphics[width=0.45\textwidth]{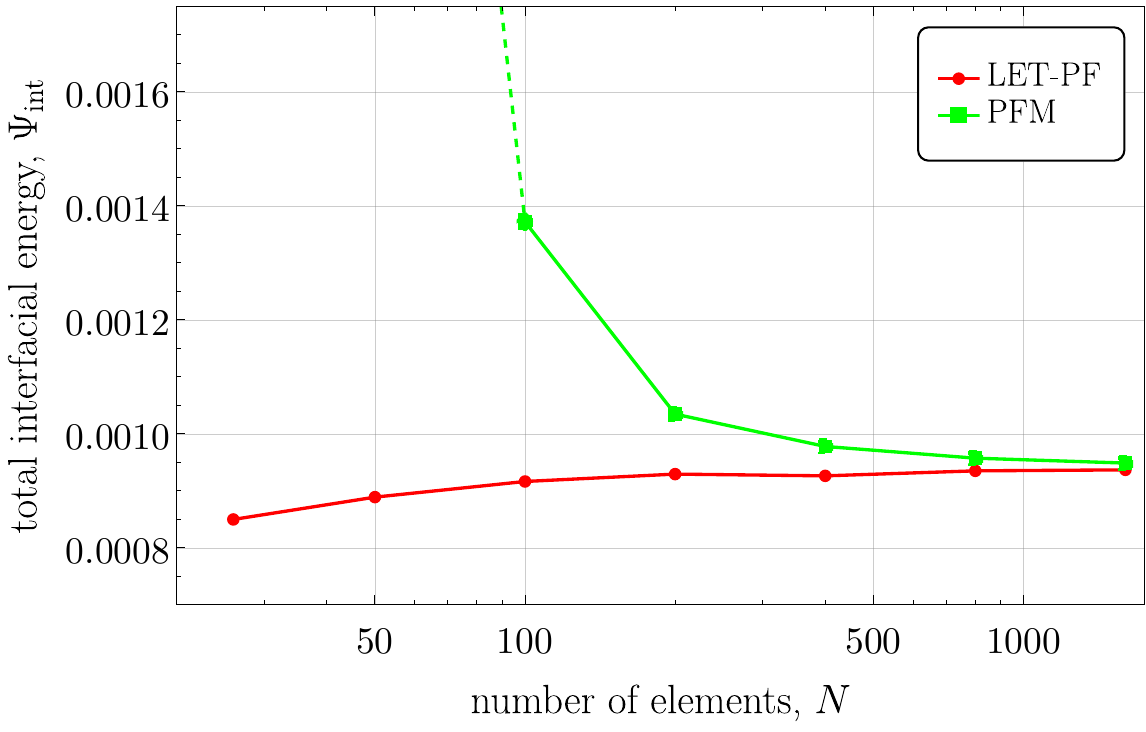}  \\[1ex]
                    \hspace*{4em}(a) & \hspace*{4em}(b) 
                \end{tabular}
                }
            \caption{Convergence of the total elastic strain energy $\Psi_{\rm el}$ (a) and of the total interfacial energy $\Psi_{\rm int}$ (b) with mesh refinement.
            The points corresponding to the trivial solution obtained for PFM for $N=26$ and $N=50$ are not included, the dashed lines show that these points lie well outside the plot range.
            }
            \label{fig:ConvRates}
        \end{figure}

        To get a better insight,
        the computations have been performed for a range of mesh densities with $N\times N$ elements where $N\in\{26,50,100,200,400,800,1600\}$. Fig.~\ref{fig:steadyStates} shows the steady-state shapes for selected mesh densities (the results for two finest meshes are omitted, as no differences are visible with further mesh refinement). Based on the results obtained, a clear advantage of LET-PF over PFM can be observed. For coarse meshes, PFM fails, yielding a trivial solution with $\phi=0.5$ in the entire domain. This stands in contrast to LET-PF which produces qualitatively correct results already for the coarsest mesh. 
        Secondly, while the final shape seems to converge to the same solution for both methods, the convergence is faster for LET-PF. 
        This is also apparent in Fig.~\ref{fig:ConvRates} which shows the total elastic strain energy $\Psi_{\rm el}$ and the total interfacial energy $\Psi_{\rm int}$ in the steady state as a function of mesh density.
        %
        %
        %
        %
        %
        %
        %
        %

\begin{figure}[htbp]
            \centerline{\scriptsize
                \begin{tabular}{cc}
                    \includegraphics[width=0.45\textwidth]{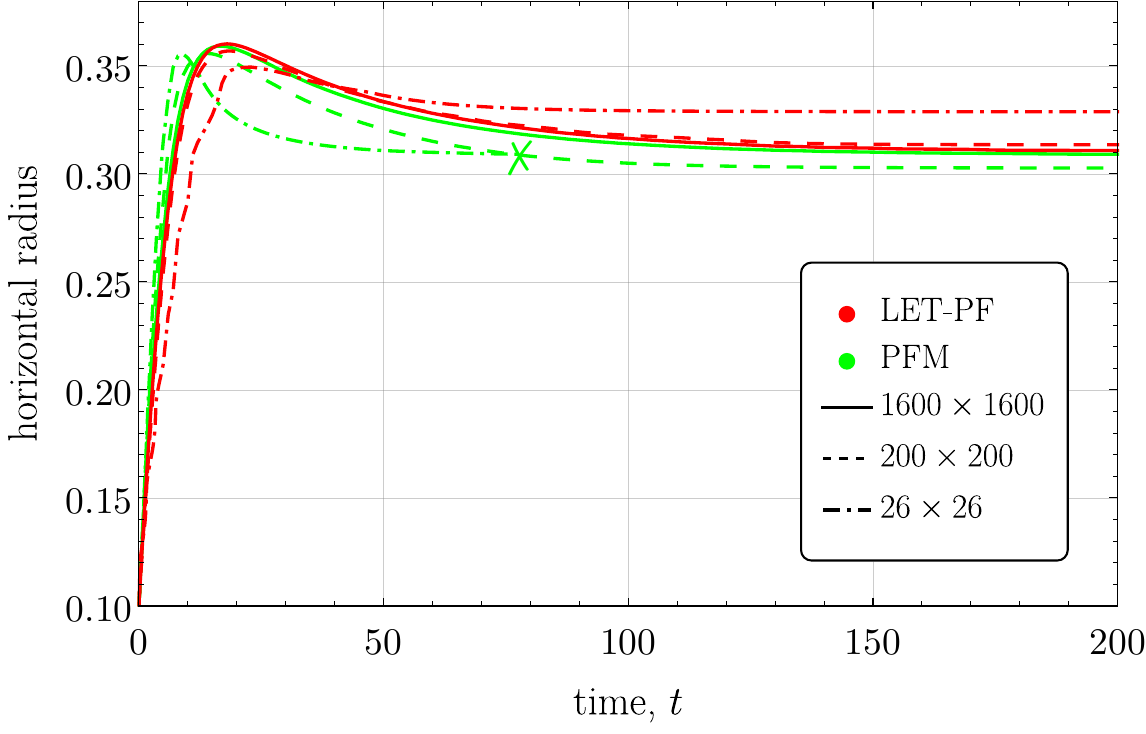} &
                    \includegraphics[width=0.45\textwidth]{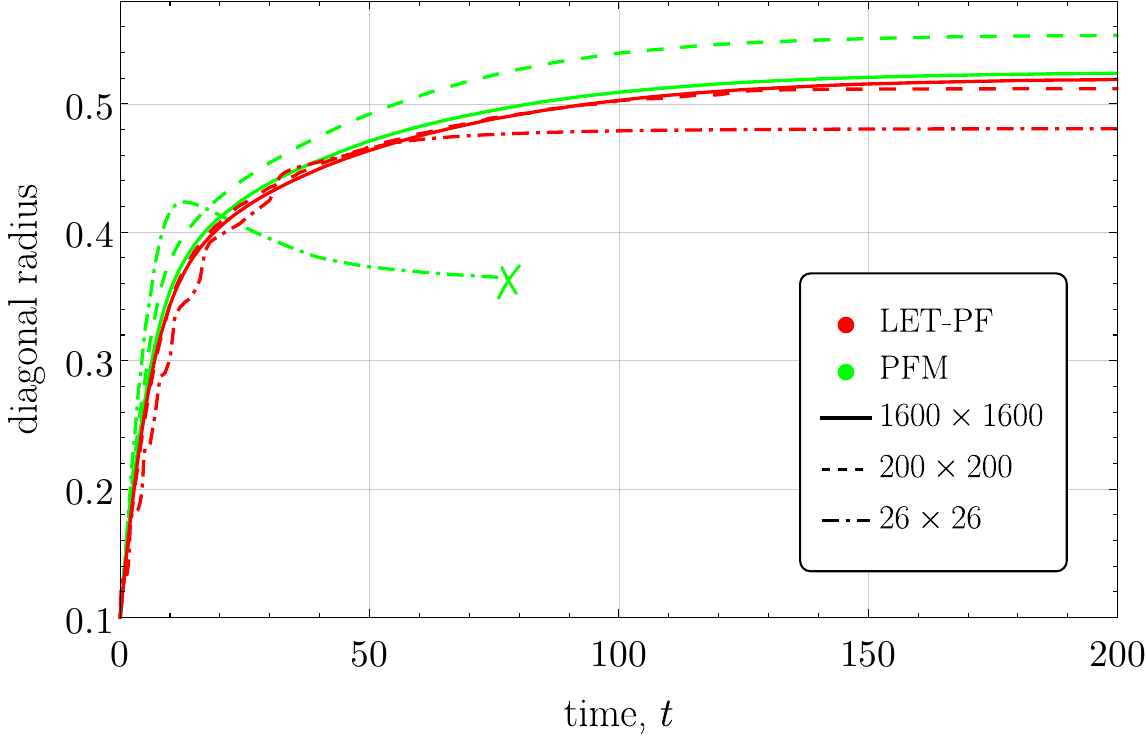}  \\[1ex]
                    \hspace*{3em}(a) & \hspace*{3em}(b) 
                \end{tabular}
                }
            \caption{Evolution of the inclusion size expressed by the distance from the centre to the propagating interface measured in the horizontal (a) and diagonal (b) directions, shown for selected mesh densities. 
            Cross markers indicate the instant when the trivial solution is achieved for PFM with $\phi=0.5$ in the entire domain.
            }
            \label{fig:radiiEvolutions}
        \end{figure}

        To examine the evolution of shape quantitatively, the distance from the centre to the interface (level set $\phi=0.5$), measured in the horizontal direction and along the diagonal, is shown in Fig.~\ref{fig:radiiEvolutions} for three mesh densities. For the finest mesh ($N=1600$), the results obtained for LET-PF and PFM are very close one to the other, which again suggests that both methods converge to the same solution. However, consistent with Fig.~\ref{fig:ConvRates}, the convergence rate is visibly different. For PFM, the difference between the results obtained for $N=200$ and $N=1600$ is significantly larger than in the case of LEF-PF. For the coarsest mesh ($N=26$), a trivial solution is obtained for PFM, 
        while LET-PF delivers a reasonably accurate solution.
        %
        %
        %
        %
        
        
        \FloatBarrier

    \subsection{Three inclusions in a constrained domain}\label{sec:three}
        As the last example, a problem similar to that of Section~\ref{sec:single} is studied with the only difference that three inclusions are considered in the initial conditions. The initial inclusion radii are 0.1, 0.15 and 0.2, and the corresponding positions of the inclusion centres are (0.25,0.25), (0.75,0.30) and (0.35,0.75), respectively, for the origin of the coordinate system located in the lower-left corner of the domain. The remaining parameters of the problem are the same as in Section~\ref{sec:single}.

\begin{figure}[htbp]
            \vspace{-5eM}
            \centerline{\scriptsize
                \begin{tabular}{c}
                    \hspace*{-0.5em}\includegraphics[width=1.1\textwidth]{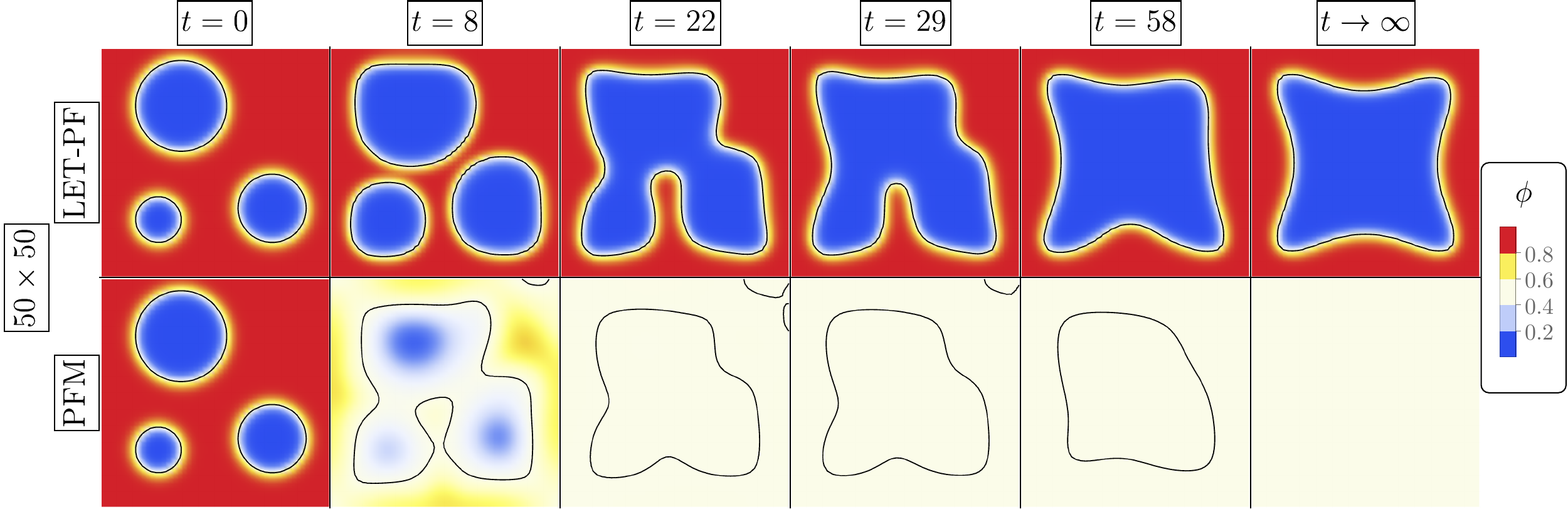}\hspace*{-0.5em}
                    \\[-0.75ex]\hline \\[-2.65ex]
                    \hspace*{-0.5em}\includegraphics[width=1.1\textwidth]{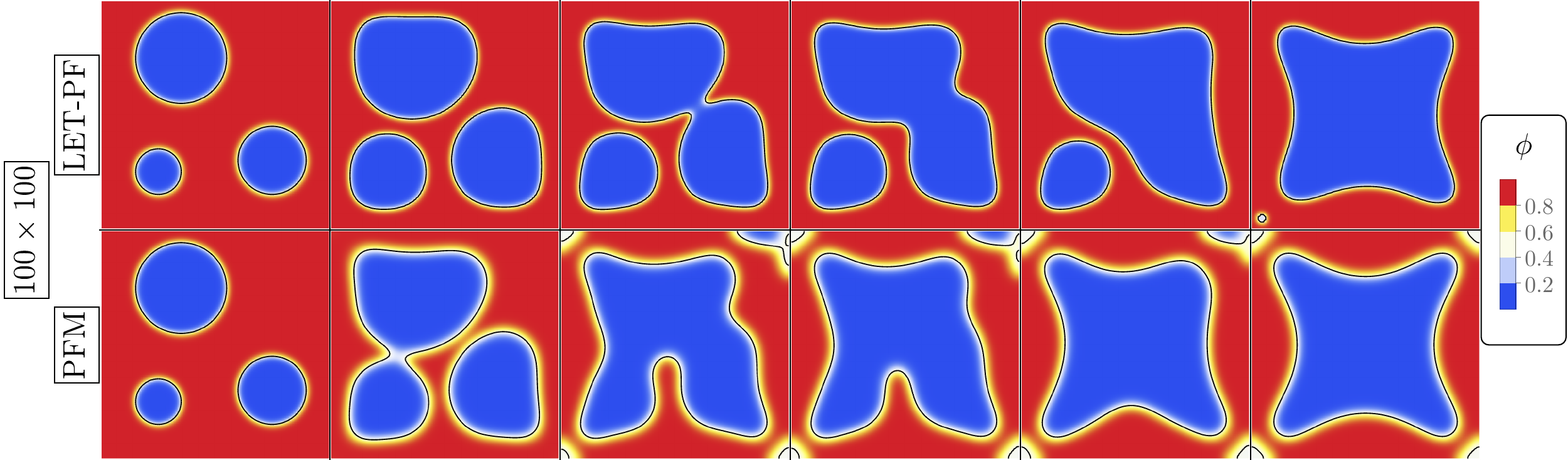}\hspace*{-0.5em}
                    \\[-0.75ex]\hline \\[-2.65ex]
                    \hspace*{-0.5em}\includegraphics[width=1.1\textwidth]{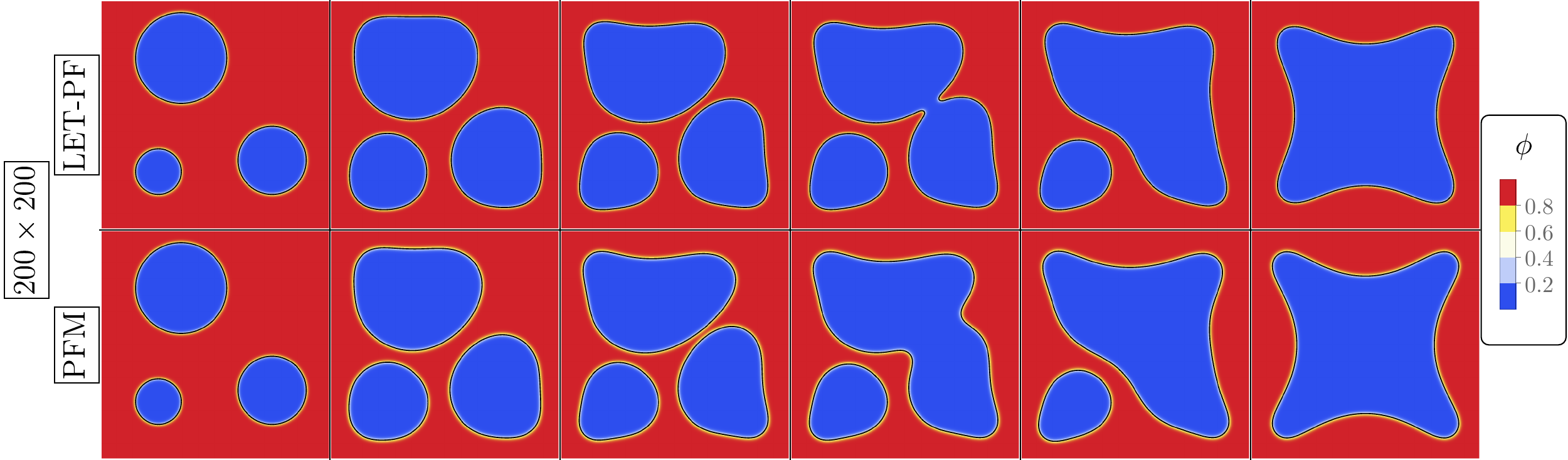}\hspace*{-0.5em}
                    \\[-0.75ex]\hline \\[-2.65ex]
                    \hspace*{-0.5em}\includegraphics[width=1.1\textwidth]{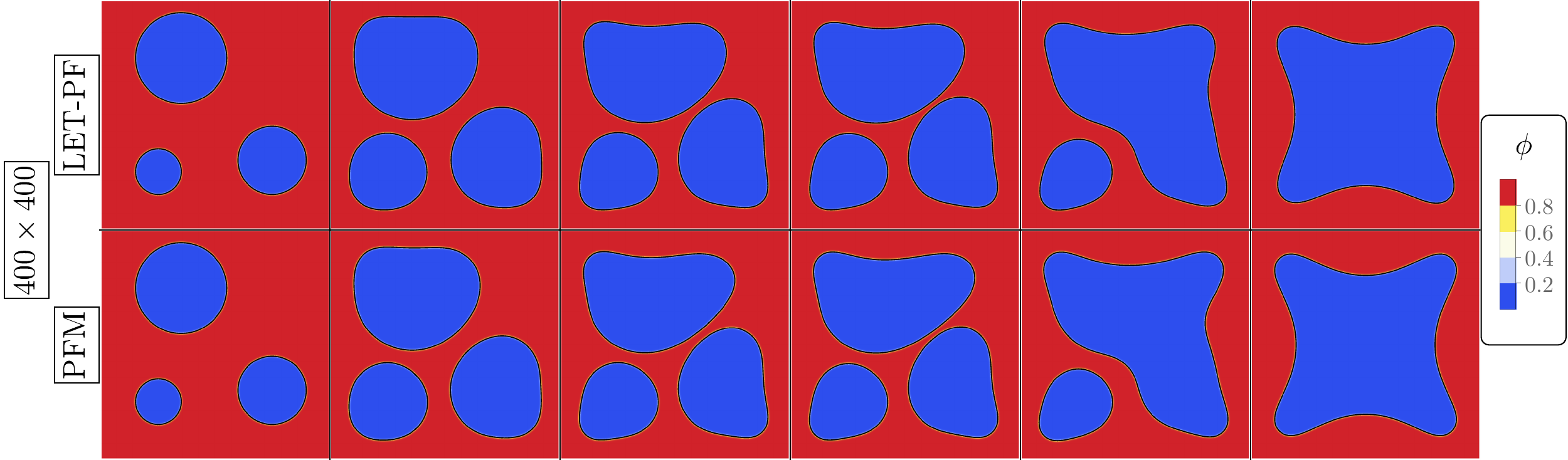}\hspace*{-0.5em}
                \end{tabular}
                }
            \caption{Evolution of three inclusions within a constrained domain obtained for four mesh densities and for LET-PF (upper snapshots) and PFM (lower snapshots). The colour maps depict the order parameter $\phi$. The black contours indicate the
            position of the $\phi=0.5$ level set.}
            \label{fig:3inclsMeshes}
        \end{figure}

        Fig.~\ref{fig:3inclsMeshes} shows selected snapshots of the evolution process for four mesh densities.
        Consider first the case of the finest mesh of $400\times400$ elements. Here, the evolution is very similar for LET-PF and PFM. At the beginning, all inclusions start to grow. Then two inclusions coalesce and grow further, while the third inclusion (that in the lower-left corner) shrinks and ultimately vanishes. 
        In the case of LET-PF, this scenario is obtained also for two coarser meshes, although the two inclusions coalesce at an earlier stage (at $t=29$ for $N=200$ and at $t=22$ for $N=100$), and only for the coarsest mesh ($N=50$) a different scenario is obtained in which all inclusions coalesce. 
        However, in the case of PFM, the correct scenario is obtained for $N=200$, even if the evolution proceeds then significantly faster (see the snapshot at $t=29$). At the same time, an incorrect scenario is obtained for $N=100$ (three inclusions coalesce), while for $N=50$ the microstructure evolves towards the trivial solution, as in Section~\ref{sec:single}. 
        %
        %
        %
        %

        The final shapes correspond, in general, to those obtained for the single-inclusion problem in Section~\ref{sec:single}. Note that, for LET-PF and $N=100$, the shrinking inclusion does not vanish completely, rather a small remnant inclusion persists in the steady state, see the lower-left corner of the domain. 
        Also, the final shape obtained for LET=PF for $N=50$ is not perfectly symmetric with the upper-right corner of the inclusion being somewhat more rounded than the other ones. 
        Both effects suggest that LET-PF may be prone to mesh pinning, and this effect may require further attention in the future.

        Nevertheless, the present numerical example confirms that LET-PF performs significantly better than PFM for the class of problems considered here, in particular, it is capable of capturing the correct evolution scenario for a coarser mesh.

\section{Conclusion}
        The laminated element technique (LET), which has been recently developed for an efficient treatment of weak discontinuities on non-matching meshes~\citep{Dobrzanski2024}, has been combined with the phase-field method and thus extended to problems involving moving interfaces and, specifically, to microstructure evolution problems. The essence of LET is in treating an element that is cut by an interface as a simple laminate with the lamination orientation and phase volume fraction depending on the orientation and position of the interface within the element. The position of a propagating interface is then determined using a continuous order parameter which plays the role of a level set function and, at the same time, carries the information about the interfacial energy. Evolution of the phase field order parameter is governed by a Ginzburg--Landau type evolution equation.

        The resulting LET-PF method has been implemented in a finite-element code, and its performance has been illustrated by a comprehensive set of numerical examples, including a comparison to the conventional phase-field method (PFM). The implementation and the computations are limited to 2D problems in the framework of small-strain elasticity, but generalization to 3D and finite deformations is immediate. The results show that, when the interfacial energy has a significant contribution to the total thermodynamic driving force (interfacial-energy-driven evolution), the two methods behave similarly in terms of both accuracy and robustness. However, when the elastic strain energy contribution dominates (elasticity-driven evolution), LET-PF is found significantly more accurate than PFM. As a result, a coarser mesh (along with a larger interface thickness) can be used in LET-PF to get results of similar accuracy, thus reducing the computational cost. 

        Overall, LET-PF seems to be a promising approach for problems involving propagating interfaces. Further work is needed to improve its robustness and to examine its performance for a wider class of problems. While generalization to a multi-phase framework is an open problem, it can be readily applied to a wide range of two-phase evolution problems.
        

    \subsection*{Acknowledgement}
        This work has been partially supported by the National Science Centre (NCN) in Poland through Grant No.\ 2018/29/B/ST8/00729.
\appendix
\section{Micro-macro transition for a simple laminate}\label{app:laminate}
    In this section, a 
    laminate composed of two phases is considered in the framework of linear elasticity with eigenstrain, and the micro-macro transition is discussed for this microstructure. 
    A simple laminate is a microstructure in which layers of two phases are separated by parallel planar interfaces, and the microstructure is defined by the interface normal $\bm{n}$ and volume fraction $\eta$. 
    Under the assumption of separation of scales, strains and stresses are homogeneous within each layer and are identical in the layers of the same phase. 
    Accordingly, the compatibility conditions formulated in Eq.~\eqref{eq:compat} for the interface jumps hold for the strains $\bm{\varepsilon}_i$ and stresses $\bm{\sigma}_i$ within individual phases,
    \begin{equation}
    \label{eq:app:compat}
    \Delta\bm{\varepsilon} = \frac12 (\bm{c}\otimes\bm{n} + \bm{n}\otimes\bm{c}) , \qquad
    \Delta\bm{\sigma}\cdot\bm{n} = \bm{0} ,
    \end{equation}
    with $\Delta(\mysquare)=(\mysquare)_2-(\mysquare)_1$, as defined in Section~\ref{sec:PFM}. 
    
    The overall quantities, denoted by a superimposed bar, are obtained by simple averaging,
    \begin{equation}
    \label{eq:averaging}
    \bar\psi=\{\psi\} , \qquad
    \bar{\bm{\varepsilon}}=\{\bm{\varepsilon}\} , \qquad
    \bar{\bm{\sigma}}=\{\bm{\sigma}\} ,
    \end{equation}
    where $\{\mysquare\}=(1-\eta)(\mysquare)_1+\eta(\mysquare)_2$. 
    Note that the overall free energy $\bar\psi$ of the laminate constitutes the bulk free energy in the LET setting, Eqs.~\eqref{eq:sigmaBar} and~\eqref{eq:GL:weak:LET}, hence $\bar\psi_{\rm bulk}=\bar\psi$, and the subscript `bulk' is used in Section~\ref{sec:hybrid} to differentiate $\bar\psi_{\rm bulk}$ from the interfacial energy contribution, as in the conventional phase-field model of Section~\ref{sec:diffuse}.
    
    Using the compatibility condition~\eqref{eq:app:compat}${}_1$ and the averaging rule~\eqref{eq:averaging}${}_2$, the local strains $\bm{\varepsilon}_i$ can be expressed as a function of the overall strain $\bar{\bm{\varepsilon}}$ and vector $\bm{c}$,
    \begin{equation}
    \label{eq:epsilon12}
    \bm{\varepsilon}_1 = \bar{\bm{\varepsilon}} - \eta \Delta\bm{\varepsilon} , \qquad
    \bm{\varepsilon}_2 = \bar{\bm{\varepsilon}} + (1-\eta) \Delta\bm{\varepsilon} .
    \end{equation}
    The unknown vector $\bm{c}$ is then determined by solving the compatibility condition~\eqref{eq:app:compat}${}_2$ with the local stresses governed by the constitutive equation $\bm{\sigma}_i=\partial\psi_i/\partial\bm{\varepsilon}_i=\mathbb{L}_i:(\bm{\varepsilon}_i-\bm{\varepsilon}^{\rm t}_i)$. 
    In the present case of linear elasticity, this yields a linear equation for vector $\bm{c}$ (note that the transformation strains $\bm{\varepsilon}^{\rm t}_i$ are here constant and known).
    
    
    In view of linearity, the complete micro-macro transition scheme can be solved in closed form. In particular, the overall free energy $\bar\psi$ can be expressed in the following form,
    \begin{equation}
    \label{eq:psiBar}
    \bar\psi = \bar\psi(\bar{\bm{\varepsilon}},\eta,\bm{n}) = \bar\psi{}^0
      + \frac12 (\bar{\bm{\varepsilon}}-\bar{\bm{\varepsilon}}^{\rm t}) : \bar{\mathbb{L}} : (\bar{\bm{\varepsilon}}-\bar{\bm{\varepsilon}}^{\rm t}) ,
    \end{equation}
    where
    \begin{equation}
    \bar\psi{}^0 = (1-\eta)\psi^0_1 + \eta \psi^0_2 + \frac12 \eta(1-\eta) \Delta\bm{\varepsilon}^{\rm t} : \mathbb{S} : \Delta\bm{\varepsilon}^{\rm t} , \qquad
    \bar{\bm{\varepsilon}}^{\rm t} = \{ \mathbb{B}^{\rm T} : \bm{\varepsilon}^{\rm t} \} .
    \end{equation}
    Here, $\bar{\mathbb{L}}$ is the overall elastic stiffness tensor, $\mathbb{B}_1$ and $\mathbb{B}_2$ are the fourth-order stress concentration tensors, and $\mathbb{S}$ is a fourth-order tensor.
    Closed-form formulae for $\bar{\mathbb{L}}$, $\mathbb{B}_1$, $\mathbb{B}_2$ and $\mathbb{S}$ in the general case of distinct anisotropy of the phases can be found, for instance, in \citep{StupkiewiczPetryk2002,Stupkiewicz2007}. 
    In the case of homogeneous elastic properties, $\mathbb{L}=\mathbb{L}_1=\mathbb{L}_2$, those formulae simplify significantly, namely $\bar{\mathbb{L}}=\mathbb{L}$ and $\mathbb{B}_1=\mathbb{B}_2=\mathbb{I}$, where $\mathbb{I}$ denotes the fourth-order identity tensor, so that $\bar{\bm{\varepsilon}}^{\rm t}=\{\bm{\varepsilon}^{\rm t}\}$, while $\mathbb{S}$ is then given by,
    \begin{equation}
    \mathbb{S} = \mathbb{L} - \mathbb{L} : (\bm{n} \otimes (\bm{n} \cdot \mathbb{L} \cdot \bm{n})^{-1} \otimes \bm{n}) : \mathbb{L} .
    \end{equation}
    
    The explicit expression~\eqref{eq:psiBar} for the the overall free energy $\bar\psi$ is provided above in order to illustrate the general structure of the micro-macro transition for a simple laminate, and the details are omitted here. 
    In the actual computer implementation of LET-PF, these closed-form formulae are not used. Rather, the micro-macro transition equations~\eqref{eq:app:compat}--\eqref{eq:epsilon12} along with the local constitutive equations are solved directly on demand. 
    The computational scheme relies on the use of the automatic differentiation (AD) technique implemented in \emph{AceGen} \cite{Korelc2009,KorelcWriggers2016}, which facilitates consistent linearization and computation of the derivatives with respect to $\eta$ and $\bm{n}$ that are needed in LET-PF, cf.\ Eq.~\eqref{eq:GL:weak:LET}. 
    The computational scheme 
    is a special case of a more general one that is discussed in Appendix~A in \citep{Dobrzanski2024} in the context of elastoplastic laminates at finite deformation. 
    The reader is referred to \citep{Dobrzanski2024} for details.
    
    
\section{Analytical solution for the evolving circular inclusion}\label{app:BimaterialLame}
    For a prescribed inclusion radius $\rho$, the elasticity problem corresponding to the setup defined in Section~\ref{sec:circular}, see Fig.~\ref{fig:inclusionscheme}(a), can be solved in a closed form by using the classical Lam\'{e} solution of a thick-walled cylinder under internal pressure. 
    Enforcing the continuity of the displacement and traction at the interface between the elastic cylinder and the elastic inclusion with eigenstrain $\bm{\varepsilon}_1^\text{t}=\epsilon\bm{I}$, the radial displacement $u_r$ can be determined as
    \begin{equation}
        \label{eq:uBiLameEigen}
        u_r(r)=
        \begin{cases}
             (1-\epsilon^\ast) \epsilon \, r
               & 0\leqslant r\leqslant \rho , \\[2ex]
             \displaystyle
             \frac{\left(\lambda_1+\mu_1\right) \left(\mu_2 r^2 + \left(\lambda_2+\mu_2\right) R^2\right) \rho^2 \epsilon^\ast \epsilon}
             {\mu_2 \left(\lambda_2+\mu_2\right) \left(R^2-\rho^2\right) r} \quad
               & \rho < r \leqslant R ,
        \end{cases}
    \end{equation}
    while the circumferential displacement is zero, $u_\theta=0$. 
    Here, $\epsilon^\ast$ denotes the fraction of the eigenstrain that is recovered elastically in the inclusion,
    \begin{equation}
        \epsilon^\ast = \frac{\mu_2 \left(\lambda_2+\mu_2\right) \left(R^2-\rho^2\right)}
            {\mu_2 \left(\lambda_1-\lambda_2+\mu_1-\mu_2\right) \rho^2 + \left(\lambda_2+\mu_2\right) \left(\lambda_1+\mu_1+\mu_2\right) R^2} ,
    \end{equation}
    so that $(1-\epsilon^\ast) \epsilon$ is the total strain in the inclusion, and $\lambda_i$ and $\mu_i$ are the Lam\'{e} constants of the phases, $i=1,2$.
    The total elastic strain energy of the system is then given by
    \begin{equation}
        \hat{\Psi}_{\rm el} = 
        2 \pi \left(\lambda_1+\mu_1\right) \rho^2 \epsilon^\ast \epsilon^2 .
    \end{equation}

    The total free energy is a sum of the elastic strain energy $\hat{\Psi}_{\rm el}$ and interfacial energy $\hat{\Psi}_{\rm int}$ (the chemical energy contribution is not considered in this example),
    \begin{equation}
        \hat{\Psi} = \hat{\Psi}_{\rm el} + \hat{\Psi}_{\rm int} , \qquad
        \hat{\Psi}_{\rm int} = \int_{\Gamma} \gamma \, \rd S = 2 \pi \rho \gamma ,
    \end{equation}
    where $\gamma$ is the interfacial energy density, assumed isotropic.

    Now consider that the interface $\Gamma$ is evolving, i.e., the inclusion radius is a function of time, $\rho=\rho(t)$. The global rate-potential $\hat{\Pi}$ is thus formulated as a sum of the rate of the total free energy, $\dot{\hat{\Psi}}=\rd\hat{\Psi}/\rd t$, and the global dissipation potential $\hat{\mathcal{D}}$,
    \begin{equation}
        \hat{\Pi} = \dot{\hat{\Psi}} + \hat{\mathcal{D}}.
    \end{equation}
    Note that no external loading is applied so that the potential energy of the loading is equal to zero and does not contribute to $\hat{\Pi}$. The dissipation is here assumed to be of purely viscous nature, i.e., the local dissipation potential $\hat{D}$ is quadratic in terms of the interface speed $\hat{v}_n=-\dot\rho$,
    \begin{equation}
        \hat{\mathcal{D}} = \int_{\Gamma}\hat{D}\left(\dot{\rho}\right)\, \rd S = 2 \pi \rho \hat{D}\left(\dot{\rho}\right) = \frac{\pi\rho\dot{\rho}^2}{\hat{m}} , \qquad
        \hat{D}(\hat{v}_n)=\frac{\hat{v}_n^2}{2\hat{m}} ,
    \end{equation}
    where $\hat{m}$ is the local interface mobility parameter.

    The evolution problem can now be formulated as a minimization problem for the global rate-potential $\hat{\Pi}$,
    \begin{equation}
        \dot{\rho} = \arg \min_{\dot{\rho}}\hat{\Pi}(\dot{\rho}) .
    \end{equation}
    Since $\hat{\Pi}$ is a quadratic function of $\dot\rho$, the condition of stationarity of $\hat{\Pi}$, $\rd\hat{\Pi}/\rd\dot\rho=0$, delivers a linear equation for $\dot\rho$ that can be easily solved. In the general case of distinct elastic properties, the resulting expression for $\dot\rho$ is rather lengthy and is not provided here for brevity.


    Assuming identical elastic properties of the two phases, $\lambda_1=\lambda_2=\lambda$ and $\mu_1=\mu_2=\mu$, the formula for $\hat{\Psi}_{\rm el}$ simplifies considerably,
    \begin{equation}
        \hat{\Psi}_{\rm el} = \frac{2 \pi \mu (\lambda +\mu) (R^2-\rho^2) \rho^2 \epsilon ^2}{(\lambda +2 \mu) R^2} ,
    \end{equation}
    and so does the resulting evolution equation,
    \begin{equation}\label{eq:rhodot:app}
        \dot{\rho} = - \hat{m}\left( \left( \frac{1}{2} - \left( \frac{\rho}{R} \right)^2 \right)\frac{E \epsilon^2}{1-\nu^2} + \frac{\gamma}{\rho}\right) ,
    \end{equation}
    where the Lam{\'e} constants have been converted into the Young's modulus $E$ and Poisson's ratio $\nu$. 
    The above equation can be rewritten in a more concise format,
    \begin{equation}\label{eq:rhodot:forces}
        \dot{\rho} = - \hat{m}\left( \hat{f}_{\rm bulk}+\hat{f}_{\rm int} \right),
        \qquad
        \hat{f}_{\rm bulk} = \frac{E \epsilon^2}{1-\nu^2} \left( \frac{1}{2} - \left( \frac{\rho}{R} \right)^2 \right) , 
        \qquad
        \hat{f}_{\rm int} = \frac{\gamma}{\rho} ,
    \end{equation}
    where $\hat{f}_{\rm bulk}$ and $\hat{f}_{\rm int}$ denote the local thermodynamic driving forces related to the elastic and interfacial energy, respectively. It can be checked that the thermodynamic driving force $\hat{f}=\hat{f}_{\rm bulk}+\hat{f}_{\rm int}$ derived above from the total energy balance can be equivalently obtained from the local definition, Eq.~\eqref{eq:kinetic}${}_2$ (to have a consistent sign of the driving force, the interface normal must point into the inclusion so that the inclusion is treated as phase 2 in the notation of Section~\ref{sec:sharp}). 


\section{Regularization of the phase volume fraction in laminated elements}\label{sec:regEffect}
\setcounter{figure}{0}
    As discussed in Section~\ref{sec:FE}, the formula~\eqref{eq:eta:el} for the phase volume fraction $\eta^{(\omega)}$ in the laminated elements is regularized according to Eq.~\eqref{eq:eta:el:reg}. In this appendix, we illustrate the effect of the regularization parameter $\phi_{\rm reg}$ on the efficiency and accuracy of LET-PF.

\begin{figure}[htbp]
        \centerline{\scriptsize
            \begin{tabular}{rcl}
                \includegraphics[height=0.3\textwidth]{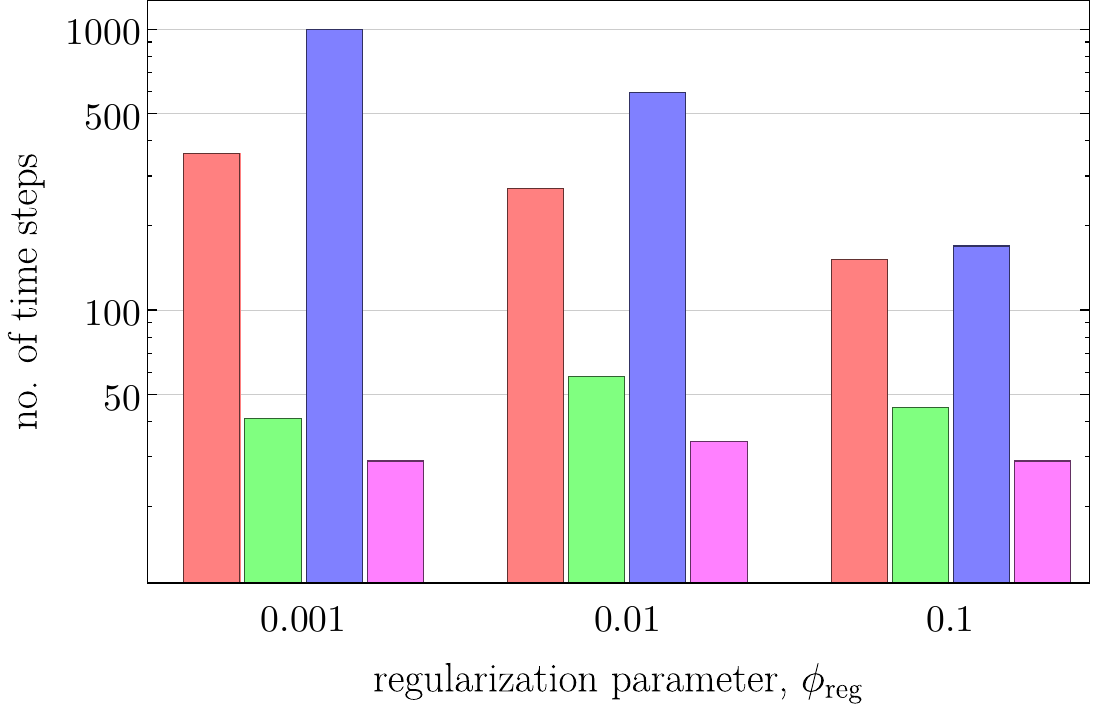} &
                \includegraphics[height=0.3\textwidth]{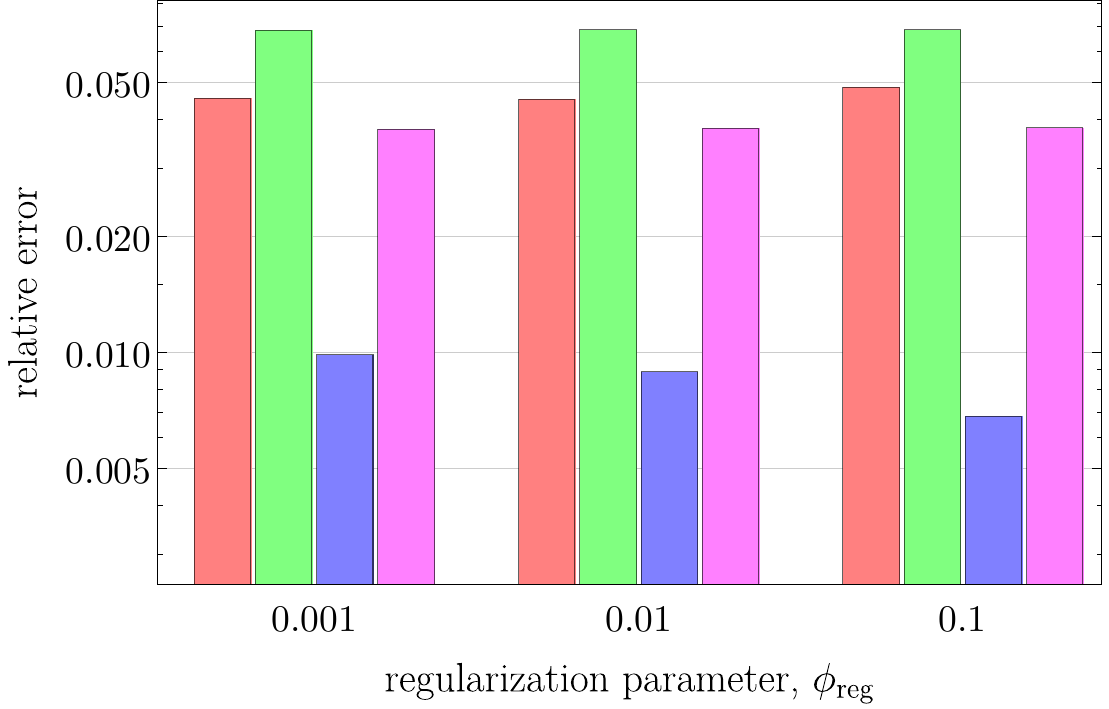} &
                \raisebox{15ex}{\includegraphics[height=0.12\textwidth]{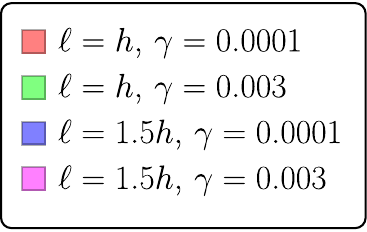}}  \\[1ex]
                (a)\hspace*{10.2eM} & \hspace*{3.4eM} (b)  & 
            \end{tabular}
            }
        \caption{The effect of regularization parameter $\phi_{\rm reg}$: (a) on the total number of time steps (for $\Delta t_{\rm max}=T_{\rm exact}$) and (b) on the relative errer (for $\Delta t_{\rm max}=T_{\rm exact}/500$).
        }
        \label{fig:regular}
    \end{figure}

    Additional simulations of the problem studied in Section~\ref{sec:circular} have thus been performed for three values of $\phi_{\rm reg}\in\{0.001,0.01,0.1\}$ and for a wide range of model parameters. 
    Fig.~\ref{fig:regular} presents representative results concerning the efficiency and robustness (quantified by the total number of time steps when $\Delta t_{\rm max}=T_{\rm exact}$) and concerning the accuracy (expressed by the relative error, Eq.~\eqref{eq:relerror}). 
    In the latter case, small time increments have been used ($\Delta t_{\rm max}=T_{\rm exact}/500$) to isolate the error related to the regularization from that resulting from time integration. 
    It has been observed that the dependence of both indicators on the interfacial energy $\gamma$ is (approximately) monotonic, hence only the extreme values $\gamma=0.0001$ and $\gamma=0.003$ are included in Fig.~\ref{fig:regular}. 
    The computations have been performed for the coarse mesh ($h=0.02$) and for two representative values of the interface thickness parameter $\ell\in\{h,1.5h\}$.
    %
    %
    %
    %

    Fig.~\ref{fig:regular}(a) shows that the number of time steps needed to complete the simulation decreases with increasing $\phi_{\rm reg}$ for $\gamma=0.0001$, and it only weakly depends on $\phi_{\rm reg}$ for $\gamma=0.003$. This confirms that the robustness and efficiency of the method increases with increasing $\phi_{\rm reg}$. Note that, in many cases, the simulations could not be completed without regularization, i.e., for $\phi_{\rm reg}=0$.
    Concerning the relative error, Fig.~\ref{fig:regular}(b), the effect of $\phi_{\rm reg}$ is weak, except one case ($\ell=1.5h$, $\gamma=0.0001$) when the error visibly increases with increasing $\phi_{\rm reg}$. 
    

    Concluding, the results presented above, as well as other results not reported here, do not give a definite answer on the choice of the regularization parameter. Higher values of $\phi_{\rm reg}$ are, in general, beneficial for robustness and computational efficiency. At the same time, the effect on accuracy is in general weak.
\bibliography{bibliografia}

\begin{thebibliography}{50}
\providecommand{\natexlab}[1]{#1}
\providecommand{\url}[1]{\texttt{#1}}
\expandafter\ifx\csname urlstyle\endcsname\relax
  \providecommand{\doi}[1]{doi: #1}\else
  \providecommand{\doi}{doi: \begingroup \urlstyle{rm}\Url}\fi

\bibitem[Chen(2002)]{Chen2002}
L.-Q. Chen.
\newblock Phase-field models for microstructure evolution.
\newblock \emph{Annual Review of Materials Research}, 32\penalty0 (1):\penalty0 113--140, 2002.
\newblock \doi{10.1146/annurev.matsci.32.112001.132041}.

\bibitem[Steinbach(2009)]{Steinbach2009}
I.~Steinbach.
\newblock Phase-field models in materials science.
\newblock \emph{Modelling and Simulation in Materials Science and Engineering}, 17\penalty0 (7):\penalty0 073001, 2009.
\newblock \doi{10.1088/0965-0393/17/7/073001}.

\bibitem[Provatas and Elder(2010)]{ProvatasElder2010}
N.~Provatas and K.~Elder.
\newblock \emph{Phase-Field Methods in Materials Science and Engineering}.
\newblock Wiley, 2010.
\newblock \doi{10.1002/9783527631520}.

\bibitem[Wang and Li(2010)]{WangLi2010}
Y.~Wang and J.~Li.
\newblock Phase field modeling of defects and deformation.
\newblock \emph{Acta Materialia}, 58\penalty0 (4):\penalty0 1212--1235, 2010.
\newblock \doi{10.1016/j.actamat.2009.10.041}.

\bibitem[Tourret et~al.(2022)Tourret, Liu, and LLorca]{Tourret2022}
D.~Tourret, H.~Liu, and J.~LLorca.
\newblock Phase-field modeling of microstructure evolution: Recent applications, perspectives and challenges.
\newblock \emph{Progress in Materials Science}, 123:\penalty0 100810, 2022.
\newblock \doi{10.1016/j.pmatsci.2021.100810}.

\bibitem[Ode et~al.(2001)Ode, Kim, and Suzuki]{Ode2001}
M.~Ode, S.~G. Kim, and T.~Suzuki.
\newblock Mathematical modeling of iron and steel making processes. recent advances in the phase-field model for solidification.
\newblock \emph{{ISIJ} International}, 41\penalty0 (10):\penalty0 1076--1082, 2001.
\newblock \doi{10.2355/isijinternational.41.1076}.

\bibitem[Ubachs et~al.(2004)Ubachs, Schreurs, and Geers]{Ubachs2004}
R.~Ubachs, P.~Schreurs, and M.~Geers.
\newblock A nonlocal diffuse interface model for microstructure evolution of tin{\textendash}lead solder.
\newblock \emph{Journal of the Mechanics and Physics of Solids}, 52\penalty0 (8):\penalty0 1763--1792, 2004.
\newblock \doi{10.1016/j.jmps.2004.02.002}.

\bibitem[Guin and Kochmann(2023)]{GuinKochmann2023}
L.~Guin and D.~M. Kochmann.
\newblock A phase-field model for ferroelectrics with general kinetics, part i: Model formulation.
\newblock \emph{Journal of the Mechanics and Physics of Solids}, 176:\penalty0 105301, 2023.
\newblock \doi{10.1016/j.jmps.2023.105301}.

\bibitem[Wang and Khachaturyan(1997)]{WangKhachaturyan1997}
Y.~Wang and A.~Khachaturyan.
\newblock Three-dimensional field model and computer modeling of martensitic transformations.
\newblock \emph{Acta Materialia}, 45\penalty0 (2):\penalty0 759--773, 1997.
\newblock \doi{10.1016/s1359-6454(96)00180-2}.

\bibitem[Levitas and Preston(2002)]{LevitasPreston2002part1}
V.~I. Levitas and D.~L. Preston.
\newblock Three-dimensional {L}andau theory for multivariant stress-induced martensitic phase transformations. {I}. {A}ustenite$\ensuremath{\leftrightarrow}$martensite.
\newblock \emph{Phys. Rev. B}, 66:\penalty0 134206, 2002.
\newblock \doi{10.1103/PhysRevB.66.134206}.

\bibitem[Xu et~al.(2020)Xu, Kang, Kan, Yu, and Xie]{Bo2020}
B.~Xu, G.~Kang, Q.~Kan, C.~Yu, and X.~Xie.
\newblock Phase field simulation on the cyclic degeneration of one-way shape memory effect of {NiTi} shape memory alloy single crystal.
\newblock \emph{International Journal of Mechanical Sciences}, 168:\penalty0 105303, 2020.
\newblock \doi{10.1016/j.ijmecsci.2019.105303}.

\bibitem[T{\r{u}}ma et~al.(2021)T{\r{u}}ma, Rezaee-Hajidehi, Hron, Farrell, and Stupkiewicz]{Tuma2021}
K.~T{\r{u}}ma, M.~Rezaee-Hajidehi, J.~Hron, P.~Farrell, and S.~Stupkiewicz.
\newblock Phase-field modeling of multivariant martensitic transformation at finite-strain: Computational aspects and large-scale finite-element simulations.
\newblock \emph{Computer Methods in Applied Mechanics and Engineering}, 377:\penalty0 113705, 2021.
\newblock \doi{10.1016/j.cma.2021.113705}.

\bibitem[Clayton and Knap(2011)]{ClaytonKnap2011}
J.~Clayton and J.~Knap.
\newblock A phase field model of deformation twinning: Nonlinear theory and numerical simulations.
\newblock \emph{Physica D: Nonlinear Phenomena}, 240\penalty0 (9-10):\penalty0 841--858, 2011.
\newblock \doi{10.1016/j.physd.2010.12.012}.

\bibitem[Liu et~al.(2018)Liu, Shanthraj, Diehl, Roters, Dong, Dong, Ding, and Raabe]{Liu2018}
C.~Liu, P.~Shanthraj, M.~Diehl, F.~Roters, S.~Dong, J.~Dong, W.~Ding, and D.~Raabe.
\newblock An integrated crystal plasticity{\textendash}phase field model for spatially resolved twin nucleation, propagation, and growth in hexagonal materials.
\newblock \emph{International Journal of Plasticity}, 106:\penalty0 203--227, 2018.
\newblock \doi{10.1016/j.ijplas.2018.03.009}.

\bibitem[Rezaee-Hajidehi et~al.(2022)Rezaee-Hajidehi, Sadowski, and Stupkiewicz]{Rezaee2022}
M.~Rezaee-Hajidehi, P.~Sadowski, and S.~Stupkiewicz.
\newblock Deformation twinning as a displacive transformation: Finite-strain phase-field model of coupled twinning and crystal plasticity.
\newblock \emph{Journal of the Mechanics and Physics of Solids}, 163:\penalty0 104855, 2022.
\newblock \doi{10.1016/j.jmps.2022.104855}.

\bibitem[Bourdin et~al.(2000)Bourdin, Francfort, and Marigo]{Bourdin2000}
B.~Bourdin, G.~Francfort, and J.-J. Marigo.
\newblock Numerical experiments in revisited brittle fracture.
\newblock \emph{Journal of the Mechanics and Physics of Solids}, 48\penalty0 (4):\penalty0 797--826, 2000.
\newblock \doi{10.1016/s0022-5096(99)00028-9}.

\bibitem[Ambati et~al.(2014)Ambati, Gerasimov, and Lorenzis]{Ambati2014}
M.~Ambati, T.~Gerasimov, and L.~D. Lorenzis.
\newblock A review on phase-field models of brittle fracture and a new fast hybrid formulation.
\newblock \emph{Computational Mechanics}, 55\penalty0 (2):\penalty0 383--405, 2014.
\newblock \doi{10.1007/s00466-014-1109-y}.

\bibitem[Cui et~al.(2021)Cui, Ma, and Mart{\'{\i}}nez-Pa{\~{n}}eda]{Cui2021}
C.~Cui, R.~Ma, and E.~Mart{\'{\i}}nez-Pa{\~{n}}eda.
\newblock A phase field formulation for dissolution-driven stress corrosion cracking.
\newblock \emph{Journal of the Mechanics and Physics of Solids}, 147:\penalty0 104254, 2021.
\newblock \doi{10.1016/j.jmps.2020.104254}.

\bibitem[T{\r{u}}ma et~al.(2016)T{\r{u}}ma, Stupkiewicz, and Petryk]{Tuma2016}
K.~T{\r{u}}ma, S.~Stupkiewicz, and H.~Petryk.
\newblock Size effects in martensitic microstructures: Finite-strain phase field model versus sharp-interface approach.
\newblock \emph{Journal of the Mechanics and Physics of Solids}, 95:\penalty0 284--307, 2016.
\newblock \doi{10.1016/j.jmps.2016.04.013}.

\bibitem[Yeddu(2018)]{Yeddu2018}
H.~K. Yeddu.
\newblock Phase-field modeling of austenite grain size effect on martensitic transformation in stainless steels.
\newblock \emph{Computational Materials Science}, 154:\penalty0 75--83, 2018.
\newblock \doi{10.1016/j.commatsci.2018.07.040}.

\bibitem[Rezaee-Hajidehi and Stupkiewicz(2020)]{Rezaee2020}
M.~Rezaee-Hajidehi and S.~Stupkiewicz.
\newblock Phase-field modeling of multivariant martensitic microstructures and size effects in nano-indentation.
\newblock \emph{Mechanics of Materials}, 141:\penalty0 103267, 2020.
\newblock \doi{10.1016/j.mechmat.2019.103267}.

\bibitem[Finel et~al.(2018)Finel, Bouar, Dabas, Appolaire, Yamada, and Mohri]{FinelPRL2018}
A.~Finel, Y.~L. Bouar, B.~Dabas, B.~Appolaire, Y.~Yamada, and T.~Mohri.
\newblock Sharp phase field method.
\newblock \emph{Physical Review Letters}, 121\penalty0 (2), 2018.
\newblock \doi{10.1103/physrevlett.121.025501}.

\bibitem[Dimokrati et~al.(2020)Dimokrati, Bouar, Benyoucef, and Finel]{Dimokrati2020}
A.~Dimokrati, Y.~L. Bouar, M.~Benyoucef, and A.~Finel.
\newblock S-{PFM} model for ideal grain growth.
\newblock \emph{Acta Materialia}, 201:\penalty0 147--157, 2020.
\newblock \doi{10.1016/j.actamat.2020.09.073}.

\bibitem[Fleck et~al.(2022)Fleck, Schleifer, and Zimbrod]{Fleck2022}
M.~Fleck, F.~Schleifer, and P.~Zimbrod.
\newblock Frictionless motion of diffuse interfaces by sharp phase-field modeling.
\newblock \emph{Crystals}, 12\penalty0 (10):\penalty0 1496, 2022.
\newblock \doi{10.3390/cryst12101496}.

\bibitem[Fleck and Schleifer(2022)]{Fleck2023}
M.~Fleck and F.~Schleifer.
\newblock Sharp phase-field modeling of isotropic solidification with a super efficient spatial resolution.
\newblock \emph{Engineering with Computers}, 39\penalty0 (3):\penalty0 1699--1709, 2022.
\newblock \doi{10.1007/s00366-022-01729-z}.

\bibitem[Dobrza{\'{n}}ski et~al.(2024)Dobrza{\'{n}}ski, Wojtacki, and Stupkiewicz]{Dobrzanski2024}
J.~Dobrza{\'{n}}ski, K.~Wojtacki, and S.~Stupkiewicz.
\newblock Lamination-based efficient treatment of weak discontinuities for non-conforming finite element meshes.
\newblock \emph{Computers {\&} Structures}, 291:\penalty0 107209, 2024.
\newblock \doi{10.1016/j.compstruc.2023.107209}.

\bibitem[G{\'{e}}l{\'{e}}bart and Ouaki(2015)]{Gelebart2015}
L.~G{\'{e}}l{\'{e}}bart and F.~Ouaki.
\newblock Filtering material properties to improve {FFT}-based methods for numerical homogenization.
\newblock \emph{Journal of Computational Physics}, 294:\penalty0 90--95, 2015.
\newblock \doi{10.1016/j.jcp.2015.03.048}.

\bibitem[Kabel et~al.(2015)Kabel, Merkert, and Schneider]{Kabel2015}
M.~Kabel, D.~Merkert, and M.~Schneider.
\newblock Use of composite voxels in {FFT}-based homogenization.
\newblock \emph{Computer Methods in Applied Mechanics and Engineering}, 294:\penalty0 168--188, 2015.
\newblock \doi{10.1016/j.cma.2015.06.003}.

\bibitem[Kabel et~al.(2017)Kabel, Fink, and Schneider]{Kabel2017}
M.~Kabel, A.~Fink, and M.~Schneider.
\newblock The composite voxel technique for inelastic problems.
\newblock \emph{Computer Methods in Applied Mechanics and Engineering}, 322:\penalty0 396--418, 2017.
\newblock \doi{10.1016/j.cma.2017.04.025}.

\bibitem[Mareau and Robert(2017)]{Mareau2017}
C.~Mareau and C.~Robert.
\newblock Different composite voxel methods for the numerical homogenization of heterogeneous inelastic materials with {FFT}-based techniques.
\newblock \emph{Mechanics of Materials}, 105:\penalty0 157--165, 2017.
\newblock \doi{10.1016/j.mechmat.2016.12.002}.

\bibitem[Keshav et~al.(2022)Keshav, Fritzen, and Kabel]{Keshav2022}
S.~Keshav, F.~Fritzen, and M.~Kabel.
\newblock {FFT}-based homogenization at finite strains using composite boxels ({ComBo}).
\newblock \emph{Computational Mechanics}, 2022.
\newblock \doi{10.1007/s00466-022-02232-4}.

\bibitem[Durga et~al.(2013)Durga, Wollants, and Moelans]{Durga2013}
A.~Durga, P.~Wollants, and N.~Moelans.
\newblock Evaluation of interfacial excess contributions in different phase-field models for elastically inhomogeneous systems.
\newblock \emph{Modelling and Simulation in Materials Science and Engineering}, 21\penalty0 (5):\penalty0 055018, 2013.
\newblock \doi{10.1088/0965-0393/21/5/055018}.

\bibitem[Mosler et~al.(2014)Mosler, Shchyglo, and Hojjat]{Mosler2014}
J.~Mosler, O.~Shchyglo, and H.~M. Hojjat.
\newblock A novel homogenization method for phase field approaches based on partial rank-one relaxation.
\newblock \emph{Journal of the Mechanics and Physics of Solids}, 68:\penalty0 251--266, 2014.
\newblock \doi{10.1016/j.jmps.2014.04.002}.

\bibitem[Schneider et~al.(2015)Schneider, Tschukin, Choudhury, Selzer, B\"{o}hlke, and Nestler]{Schneider2015}
D.~Schneider, O.~Tschukin, A.~Choudhury, M.~Selzer, T.~B\"{o}hlke, and B.~Nestler.
\newblock Phase-field elasticity model based on mechanical jump conditions.
\newblock \emph{Computational Mechanics}, 55\penalty0 (5):\penalty0 887--901, 2015.
\newblock \doi{10.1007/s00466-015-1141-6}.

\bibitem[Bartels and Mosler(2017)]{Bartels2017}
A.~Bartels and J.~Mosler.
\newblock Efficient variational constitutive updates for {A}llen{\textendash}{C}ahn-type phase field theory coupled to continuum mechanics.
\newblock \emph{Computer Methods in Applied Mechanics and Engineering}, 317:\penalty0 55--83, 2017.
\newblock \doi{10.1016/j.cma.2016.11.024}.

\bibitem[Sethian(1999)]{Sethian-book1999}
J.~A. Sethian.
\newblock \emph{Level {S}et {M}ethods and {F}ast {M}arching {M}ethods: Evolving {I}nterfaces in {C}omputational {G}eometry, {F}luid {M}echanics, {C}omputer {V}ision, and {M}aterials {S}cience}.
\newblock Cambridge {M}onographs on {A}pplied and {C}omputational {M}athematics ({N}o. 3). Cambridge University Press, Cambridge, England, 2nd edition, 1999.

\bibitem[Moës et~al.(1999)Moës, Dolbow, and Belytschko]{Moes1999}
N.~Moës, J.~Dolbow, and T.~Belytschko.
\newblock A finite element method for crack growth without remeshing.
\newblock \emph{International Journal for Numerical Methods in Engineering}, 46\penalty0 (1):\penalty0 131--150, 1999.
\newblock \doi{10.1002/(sici)1097-0207(19990910)46:1<131::aid-nme726>3.0.co;2-j}.

\bibitem[Sukumar et~al.(2001)Sukumar, Chopp, Moës, and Belytschko]{Sukumar2001}
N.~Sukumar, D.~Chopp, N.~Moës, and T.~Belytschko.
\newblock Modeling holes and inclusions by level sets in the extended finite-element method.
\newblock \emph{Computer Methods in Applied Mechanics and Engineering}, 190\penalty0 (46-47):\penalty0 6183--6200, 2001.
\newblock \doi{10.1016/s0045-7825(01)00215-8}.

\bibitem[Moës et~al.(2003)Moës, Cloirec, Cartraud, and Remacle]{Moes2003}
N.~Moës, M.~Cloirec, P.~Cartraud, and J.-F. Remacle.
\newblock A computational approach to handle complex microstructure geometries.
\newblock \emph{Computer Methods in Applied Mechanics and Engineering}, 192\penalty0 (28-30):\penalty0 3163--3177, 2003.
\newblock \doi{10.1016/s0045-7825(03)00346-3}.

\bibitem[Ji et~al.(2002)Ji, Chopp, and Dolbow]{Ji2002}
H.~Ji, D.~Chopp, and J.~E. Dolbow.
\newblock A hybrid extended finite element/level set method for modeling phase transformations.
\newblock \emph{International Journal for Numerical Methods in Engineering}, 54\penalty0 (8):\penalty0 1209--1233, 2002.
\newblock \doi{10.1002/nme.468}.

\bibitem[Duddu et~al.(2011)Duddu, Chopp, Voorhees, and Moran]{Duddu2011}
R.~Duddu, D.~L. Chopp, P.~Voorhees, and B.~Moran.
\newblock Diffusional evolution of precipitates in elastic media using the extended finite element and the level set methods.
\newblock \emph{Journal of Computational Physics}, 230\penalty0 (4):\penalty0 1249--1264, 2011.
\newblock \doi{10.1016/j.jcp.2010.11.002}.

\bibitem[Munk et~al.(2022)Munk, Reschka, Maier, Wriggers, and L\"{o}hnert]{Munk2022}
L.~Munk, S.~Reschka, H.~J. Maier, P.~Wriggers, and S.~L\"{o}hnert.
\newblock A sharp-interface model of the diffusive phase transformation in a nickel-based superalloy.
\newblock \emph{Metals}, 12\penalty0 (8):\penalty0 1261, 2022.
\newblock \doi{10.3390/met12081261}.

\bibitem[Gurtin(2000)]{Gurtin2000}
M.~E. Gurtin.
\newblock \emph{Configurational Forces as Basic Concepts of Continuum Physics}.
\newblock Springer New York, 2000.
\newblock \doi{10.1007/b97847}.

\bibitem[Moës et~al.(2023)Moës, Remacle, Lambrechts, Lé, and Chevaugeon]{Moes2023}
N.~Moës, J.-F. Remacle, J.~Lambrechts, B.~Lé, and N.~Chevaugeon.
\newblock The e{X}treme {M}esh deformation approach ({X-MESH}) for the {S}tefan phase change model.
\newblock \emph{Journal of Computational Physics}, 477:\penalty0 111878, 2023.
\newblock ISSN 0021-9991.
\newblock \doi{10.1016/j.jcp.2022.111878}.

\bibitem[Hildebrand and Miehe(2012)]{HildebrandMiehe2012}
F.~Hildebrand and C.~Miehe.
\newblock A phase field model for the formation and evolution of martensitic laminate microstructure at finite strains.
\newblock \emph{Philosophical Magazine}, 92\penalty0 (34):\penalty0 4250--4290, 2012.
\newblock \doi{10.1080/14786435.2012.705039}.

\bibitem[Penrose and Fife(1990)]{PenroseFife1990}
O.~Penrose and P.~C. Fife.
\newblock Thermodynamically consistent models of phase-field type for the kinetic of phase transitions.
\newblock \emph{Physica D: Nonlinear Phenomena}, 43\penalty0 (1):\penalty0 44--62, 1990.
\newblock \doi{10.1016/0167-2789(90)90015-h}.

\bibitem[Korelc(2009)]{Korelc2009}
J.~Korelc.
\newblock Automation of primal and sensitivity analysis of transient coupled problems.
\newblock \emph{Computational Mechanics}, 44\penalty0 (5):\penalty0 631--649, 2009.
\newblock \doi{10.1007/s00466-009-0395-2}.

\bibitem[Korelc and Wriggers(2016)]{KorelcWriggers2016}
J.~Korelc and P.~Wriggers.
\newblock \emph{Automation of Finite Element Methods}.
\newblock Springer International Publishing, 2016.
\newblock \doi{10.1007/978-3-319-39005-5}.

\bibitem[Stupkiewicz and Petryk(2002)]{StupkiewiczPetryk2002}
S.~Stupkiewicz and H.~Petryk.
\newblock Modelling of laminated microstructures in stress-induced martensitic transformations.
\newblock \emph{Journal of the Mechanics and Physics of Solids}, 50\penalty0 (11):\penalty0 2303--2331, 2002.
\newblock \doi{10.1016/s0022-5096(02)00029-7}.

\bibitem[Stupkiewicz(2007)]{Stupkiewicz2007}
S.~Stupkiewicz.
\newblock \emph{Micromechanics of Contact and Interphase Layers}.
\newblock Springer Berlin Heidelberg, 2007.
\newblock \doi{10.1007/978-3-540-49717-2}.

\end{thebibliography}
\end{document}